\documentclass[a4paper]{amsart}

\usepackage{graphicx}%
\usepackage{multirow}%
\usepackage{amsmath,amssymb,amsfonts}%
\usepackage{amsthm}%
\usepackage{mathrsfs}%
\usepackage[title]{appendix}%
\usepackage{xcolor}%
\usepackage{textcomp}%
\usepackage{manyfoot}%
\usepackage{booktabs}%
\usepackage{algorithm}%
\usepackage{algorithmicx}%
\usepackage{algpseudocode}%
\usepackage{listings}%
\usepackage{todonotes}

\usepackage{epstopdf}
\usepackage{tikz-cd}

\newtheorem{theorem}{Theorem}[section]
\newtheorem{proposition}[theorem]{Proposition}%
\newtheorem{definition}[theorem]{Definition}%
\newtheorem{remark}[theorem]{Remark}%

\newtheorem{corollary}[theorem]{Corollary}
\newtheorem{Lemma}[theorem]{Lemma}

\newcommand{\lieg}{{\tilde{\mathfrak g}}} 
\newcommand{\Lieg}{{\mathfrak g}} 
\newcommand{\norm}[1]{\lVert#1\rVert}

\newcommand{\R}{\mathbb{R}}

\newcommand{\N}{\mathbb{N}}

\newcommand{\eps}{\varepsilon}
\newcommand{\beq}{\begin{equation}}
\newcommand{\eeq}{\end{equation}}
\newcommand{\bqn}{\begin{eqnarray}}
\newcommand{\eqn}{\end{eqnarray}}

\definecolor{lgray}{rgb}{0.7,0.8,0.9}
\definecolor{sina}{rgb}{1.0,.0,1.0}         

\definecolor{sofya}{rgb}{0.0, 0.7, 0.3}

\definecolor{kathrin}{rgb}{0.3,.0,1.0}     

\newcommand{\diff}{\,{\rm d}}

\graphicspath{{figures/}}

\title[Trim Turnpikes for OCPs with Symmetries]{Trim Turnpikes for Optimal Control Problems with Symmetries}

\author{K.\ Flaßkamp}
\address{Kathrin Flaßkamp \\
	Systems Modeling and Simulation, Saarland University, Germany}
\email{kathrin.flasskamp@uni-saarland.de}

\author{S.\ Maslovskaya }
\address{Sofya Maslovskaya \\
		Department of Mathematics, Paderborn University, Paderborn, Germany}
\email{sofyam@math.upb.de}

\author{S.\ Ober-Blöbaum }
\address{Sina Ober-Blöbaum \\
	Department of Mathematics, Paderborn University, Paderborn, Germany}
\email{sinaober@math.upb.de}

\author{B.\ Wembe}
\address{Boris Wembe \\Department of Mathematics, Paderborn University, Paderborn, Germany}
\email{wboris@math.upb.de}

\begin{document}


\begin{abstract} 
Motivated by mechanical systems with symmetries, we focus on optimal control problems possessing certain symmetries. Following recent works \cite{FFOSW22, Tre:23}, which generalized the classical concept of \emph{static turnpike to \emph{manifold turnpike}} 
we extend the \emph{exponential turnpike property} to the \emph{exponential trim turnpike} for control systems with symmetries induced by abelian or non-abelian groups. 
Our analysis is mainly based on the geometric reduction of control systems with symmetries. More concretely, we first reduce the control system on the quotient space and state the turnpike theorem for the reduced problem. Then we use the group properties to obtain the \emph{trim turnpike theorem} for the full problem.
Finally, we illustrate our results on the Kepler problem and the Rigid body problem.
\end{abstract}


\maketitle

\section{Introduction}
\label{sec1}

The turnpike property in optimal control problems (OCPs) refers to a quasi-static asymptotic behavior of optimal solutions.
Typically, three characteristic parts can be distinguished in every optimal trajectory on a sufficiently large time horizon, which satisfies the turnpike property. 
In the first part, the solution quickly converges from the initial conditions to a certain steady-state. In the second part it stays for most of the time in some neighborhood of the steady-state. Finally, in the third part, it quickly leaves the neighborhood to satisfy the finial conditions. Such a behavior has been observed and analyzed in various classes of OCPs. 
Thus, turnpikes have become a central topic in optimal control theory.
Originating from the analysis of problems arising in economics, they have subsequently appeared ubiquitous in many application areas of optimal control. Initially, the term turnpike has been introduced by \cite{DSS:58} and popularized by \cite{McKenzie:76} and \cite{CHL:91}. Earlier reports on turnpike phenomena in optimal control theory date back to \cite{Neumann:38}. 

The vast majority of research focuses on steady-state turnpikes \cite{FKJB17, GTZ:16, Tre_Zua:2014}. The turnpike phenomenon in this case is characterized by the convergence of the solutions to a steady-state, as it was described initially. In this setting,
one of the most important results is the \emph{exponential turnpike theorem}, which gives precision on the exponential type of the convergence.

However, for other important classes of OCPs, optimal solutions do not converge to a static point but toward a more complicated set in state space. 
Several recent works have introduced different notions of non-static turnpikes, as for instance  \emph{periodic turnpike} \cite{Samuelson:76}, \emph{linear turnpike} \cite{Tre:23}, \emph{trim turnpike} \cite{FFOSW22} etc. 
In \cite{FFOSW22}, the turnpike was set into relation with 
symmetries of the optimal control problem in the context of mechanical systems and for a particular class of symmetries. In such case, the turnpikes are given by trim primitives, or shortly trims.
Trim primitives have been introduced by Frazzoli et al.~\cite{frazzoli2005} as an important class of constantly controlled motion primitives that can be used for motion planning, among others, in autonomous driving applications \cite{pedrosa2021, scheffe2022}. 
Library-based motion planning greatly benefits from system symmetries, because symmetric motions do not have to be included into the library more than with one representative.
In classical mechanics, symmetry is known to give rise to the existence of relative equilibria, i.e.\ purely symmetry-induced motions of the mechanical system. Trim primitives can be seen as the extension of relative equilibria to the control system case by
restricting to constant control values \cite{frazzoli2005, FOK12}.

The main goal of this work is to show \emph{trim-exponential turnpike property}, i.e., the exponential convergence of solutions toward trims for a general class of optimal control problems with symmetries.
To this aim, we consider an optimal control problem for a dynamical control system on a smooth manifold $M$ with control space $\R^m$, given by
\begin{equation} 
\label{eq:OCP_original}
\tag{OCP}
\begin{aligned}
\min_{u} &~~ J(u) = \int_{0}^{T} f^0(x(t),u(t)) d t, \\
\text{s.t.} & ~~ 
\begin{cases} 
\dot{x}(t) &= f(x(t),u(t)), ~~ x(0) = x_0,\\
0 &= \varphi_T(x(T)),
\end{cases}
\end{aligned}
\end{equation}
where $f: M \times \R^m \rightarrow TM$, $f^0: M \times \R^m \rightarrow \R$ and $\varphi_T : M \mapsto \R^{k_T}$ ($k_T \in \mathbb{N}_0$), are smooth. 
In this paper, we assume that \eqref{eq:OCP_original} admits a symmetry, i.e.\ 
there exists an action $\Phi : G \times M \rightarrow M$ of a Lie group $G$ on $M$ such that $f$ is equivariant w.r.t.\ 
$\Phi$ and $f^0$, $\varphi_T$ are invariant w.r.t.\ $\Phi$ \cite{Ohsawa13}. These conditions insure that the optimal control is invariant w.r.t.\ 
the group action.

Symmetry in terms of $G$-equivariant vector fields $f$ often arises in mechanical systems, where symmetries are related to conservation laws. Classical examples of
mechanical systems with symmetry, that have been studied extensively in the literature, are the double spherical pendulum \cite{MS93}, the Kepler (two-body) problem \cite{FFOSW22}, rigid body system \cite{BL23} or the snake-board system \cite{Ohsawa13}. 

In some cases, such as for the spherical pendulum or the Kepler problem, the symmetry group $G$ is abelian and one can find coordinates for $M$ such that the state $x$ can be subdivided into shape variable $y \in M/G$ and cyclic variable $\theta \in G$.
In this case, \eqref{eq:OCP_original} takes the form
\begin{equation} 
\label{eq:OCP_bis}
\begin{aligned}
    \min_{u}  &~~ J(u) = \int_{0}^{T} f^0(y(t),u(t)) d t, \\
    \text{s.t.} & ~~ 
    \begin{cases}
    \dot{y}(t) &= f_1(y(t),u(t)), ~~ y(0) = y_0 \\
    \dot{\theta}(t) &= f_2(y(t),u(t)), ~~ \theta(0) = \theta_0 \\
    0 &= \varphi_T (y(T)),
    \end{cases} 
\end{aligned}
\end{equation}
keeping in mind the assumptions that $f^0$ and $\varphi_T$ shall share the invariance of the system solutions given by the Lie group symmetry. Then, we can roughly sketch the concept of \emph{trim turnpikes} for the abelian case by analyzing the symmetry-reduced OCP for existence of \emph{static turnpikes}.
Since $\theta$ is a cyclic variable with a free-final position, problem~\eqref{eq:OCP_bis} is independent of $\theta$, the turnpike can then be obtained by first considering the reduced system in which we disregard the dynamics in $\theta$.
The solution of this system exhibiting the turnpike property are characterised by the following behaviour
\[
y(t) \sim \bar y, \quad u(t) \sim \bar u, \quad \text{for} ~~ t \in [\eta, T - \eta],
\]
for some small $\eta > 0$ and $(\bar y, \bar u)$ solution of the static problem
\[
\min \, \{ f^0(\bar y, \bar u) \, | \, (\bar y, \bar u) \in M/G \times \R^m, ~~ f_1(\bar y,\bar u) = 0 \}.
\]
Then, taking into account the dynamics of the cyclic variable, we deduce the trim turnpike for the full problem \eqref{eq:OCP_bis} given by: 
\[
\bar y, \ \bar u, \ \bar \theta(t) = h
f_2(\bar y, \bar u) t + \theta_0.
\]
These calculations show that a static turnpike of the reduced problem defined on the quotient space $M/G$ leads to a non-static turnpike in the group variable. For an optimal control problem, it means that under certain conditions on the reduced problem, the solutions are quasi-static in the shape variable and converge to a trajectory on a group orbit.

However, in the case of non-abelian groups, the splitting in configuration and cyclic variables as in problem \eqref{eq:OCP_bis} is not straightforward and more work has to be done to establish turnpike property.
Thus, the first contribution of this work is
to generalize this procedure to  more general Lie group symmetry cases and completely characterize the trajectory to which the optimal solution converges.
The second contribution which constitutes the main result, is to establish from the resulting reduced system a \emph{trim exponential turnpike property}.

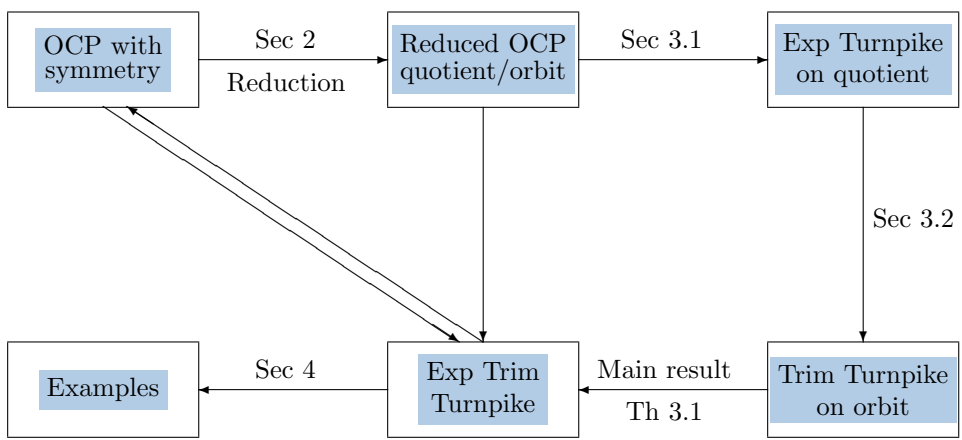
\begin{figure}[ht]
\centering
\vspace{4.cm}
\unitlength=6.25mm
\begin{picture}(28,12)

\put(0,16){\framebox(4,2){\colorbox{lgray}{\shortstack{OCP with \\ symmetry}}}}
\put(5.2,17.3){\parbox{\textwidth}{Sec~2}}
\put(4.6,16.4){\parbox{\textwidth}{Reduction}}
\put(4,17){\vector(1,0){4}}
\put(10,16){\vector(0,-1){5}}
\put(10,11){\vector(-3,2){7.5}}
\put( 2,16){\vector(3,-2){7.5}}

\put(8,16){\framebox(4,2){\colorbox{lgray}{\shortstack{Reduced OCP \\ quotient/orbit}}}}
\put(12.9,17.3){\parbox{\textwidth}{Sec~3.1}}
\put(12,17){\vector(1,0){4}}

\put(16.,16){\framebox(4,2){\colorbox{lgray}{\shortstack{Exp Turnpike \\ on quotient}}}}
\put(18.2,13.5){\parbox{\textwidth}{Sec~3.2}}
\put(18,16){\vector(0,-1){5}}

\put(16,9){\framebox(4,2){\colorbox{lgray}{\shortstack{Trim Turnpike \\ on orbit}}}}
\put(12.4,10.3){\parbox{\textwidth}{Main result}}
\put(13,9.4){\parbox{\textwidth}{Th~3.1}}

\put(8,9){\framebox(4,2){\colorbox{lgray}{\shortstack{Exp Trim \\ Turnpike}}}}
\put(16,10){\vector(-1,0){4}}

\put(8,10){\vector(-1,0){4}}
\put(5.2,10.3){\parbox{\textwidth}{Sec~4}}
\put(0,9){\framebox(4,2){\colorbox{lgray}{\shortstack{Examples}}}}

\end{picture}
\vspace{-6.cm}
\caption{Overview of the remainder of the paper.}
\label{fig:overview}
\end{figure}

The rest of the paper is organized as follows, see also the flowchart in Figure~\ref{fig:overview}: In Section~\ref{Sec:symm.reduction}, we recall some useful geometric tools as well as the reduction procedures for general optimal control problems \eqref{eq:OCP_original} and for mechanical systems 
from literature. We show the equivalence of the two reductions for controlled mechanical systems. We analyze from the obtained reduced system the turnpike property -- first on the quotient space, then on the group -- in order to set our main result which is the \emph{Trim exponential turnpike} in Section~\ref{Sec:turnpike}. Then, Section~\ref{sec:examples} is devoted to examples. We first consider the Kepler problem and then the Rigid body problem and also verify our analytic results by numerical computations. We end the paper by a conclusion, where we discuss some open research questions.

\section{Reduction of systems with symmetries}  \label{Sec:symm.reduction}


\subsection{Definitions and notations} 
\label{Sec:def.symm}

We will extensively use notions related to Lie group actions on manifolds along the paper, the most important notions are introduced in this subsection and for more details we refer to~\cite{MarsdenRatiu2013,CMR:2001b}. Consider a left action of a Lie group $G$ on a smooth manifold $M$ defined by a smooth map $\Phi_g : G\times M \rightarrow M$. For any $g\in G$ and $x \in M$, we denote $\Phi_g(x) = \Phi(g,x) = g x$. 

\begin{definition}
 We call \emph{$G$-orbit} (or just orbit if no confusion) at $x \in M$ the set $Gx \subseteq M$ defined by $$Gx = \{ z \in M \ : \ z = gx, \ g \in G\}.$$
\end{definition}

The Lie algebra of $G$ is denoted by $\mathfrak g$. We assume that the group action is free and proper, in this case the group action induces a \emph{principal bundle structure} on $M$, i.e.\ there exists a projection map $\pi :  M \rightarrow M/G$ with $M/G$ smooth manifold and $\pi$ submersion. Elements of $M/G$ are equivalence classes denoted by $[x]_G$, with equivalence relationship on $M$ being to belong to the same orbit $Gx$.
Principal bundle admits local trivialization, that is, locally it takes the form $M = M/G \times G$. If a principal bundle admits the product form globally, then it is called \emph{trivial bundle}. The group action on $M$ admits a natural lift to the action on $TM$ and $T^*M$ 
which endows $TM$ and $T^*M$ with the principal bundle structure. In a local trivialization $(x,v)$ of $TM$, the action on $TM$ is defined by $d\Phi_{g}(x, v) = (\Phi_g(x), d\Phi_g(x)(v))$. The action on $T^*M$ is defined in local trivialization $(x,\lambda)$ of $T^*M$ by $(d\Phi_{g^{-1}})^*(x, \lambda) = (\Phi_g(x), (d{\Phi_g(x)})^*\Phi_{g^{-1}}(\lambda))$, where for any $w \in T_{\Phi_g(x)}M, \ \langle (d{\Phi_g(x)})^*\Phi_{g^{-1}}(\lambda) \ | \ w \rangle = \langle \lambda \ | \ d{\Phi_g(x)}\Phi_{g^{-1}}(w) \rangle$ with $\langle \cdot \ | \ \cdot \rangle: T^*M \times TM \to \mathbb{R}$ the canonical pairing 
of co-vectors and vectors.

\begin{definition}
Let  $\xi \in {\mathfrak g}$. We call a \emph{flow of $\xi$} from $x\in M$ a trajectory on $M$ defined by $\phi^\xi(t) = \Phi(\exp(\xi t), x)$ with $\exp(\xi t)$ the Lie group exponential. The corresponding vector field $\xi_x = \left. \frac{d}{dt}\phi^\xi (t)\right|_{t=0}$ is called infinitesimal generator. 
\end{definition}

The notion of a flow of $\xi \in {\mathfrak g}$ is closely related to relative equilibria in dynamical systems \cite{MS93} and the analogous notion in control systems is given by trim primitives \cite{frazzoli2005} and can be defined as follows.

\begin{definition}
Solution $x(\cdot)$ of $\dot{x} = f(x,u)$, with constant control $u$, is called \emph{trim} if $x(\cdot) = \phi^\xi(\cdot)$ for some $\xi \in \mathfrak{g}$. 
\end{definition}

\begin{definition}
\emph{A principal connection} $\mathcal{A}$ on $M$ is a one form $\mathcal{A} : TM \rightarrow \mathfrak{g}$ with the following properties
\begin{enumerate}
\item[i)] $\mathcal{A}(\xi_x) = \xi$, for all $\xi \in \mathfrak{g}$;
\item[ii)] $\mathcal{A}(T_x \Phi_g \cdot v) = \mathrm{Ad}_g(\mathcal{A}(v))$, with $\mathrm{Ad}_g$ the adjoint action of $G$ on $\mathfrak{g}$.
\end{enumerate}
Restriction of a connection $\mathcal{A}$ on $T_x M$ is denoted by $\mathcal{A}_x$.
\end{definition} 

\begin{figure}[ht!]
\centering  
\includegraphics[width=0.4\textwidth, height = 0.37\textwidth]{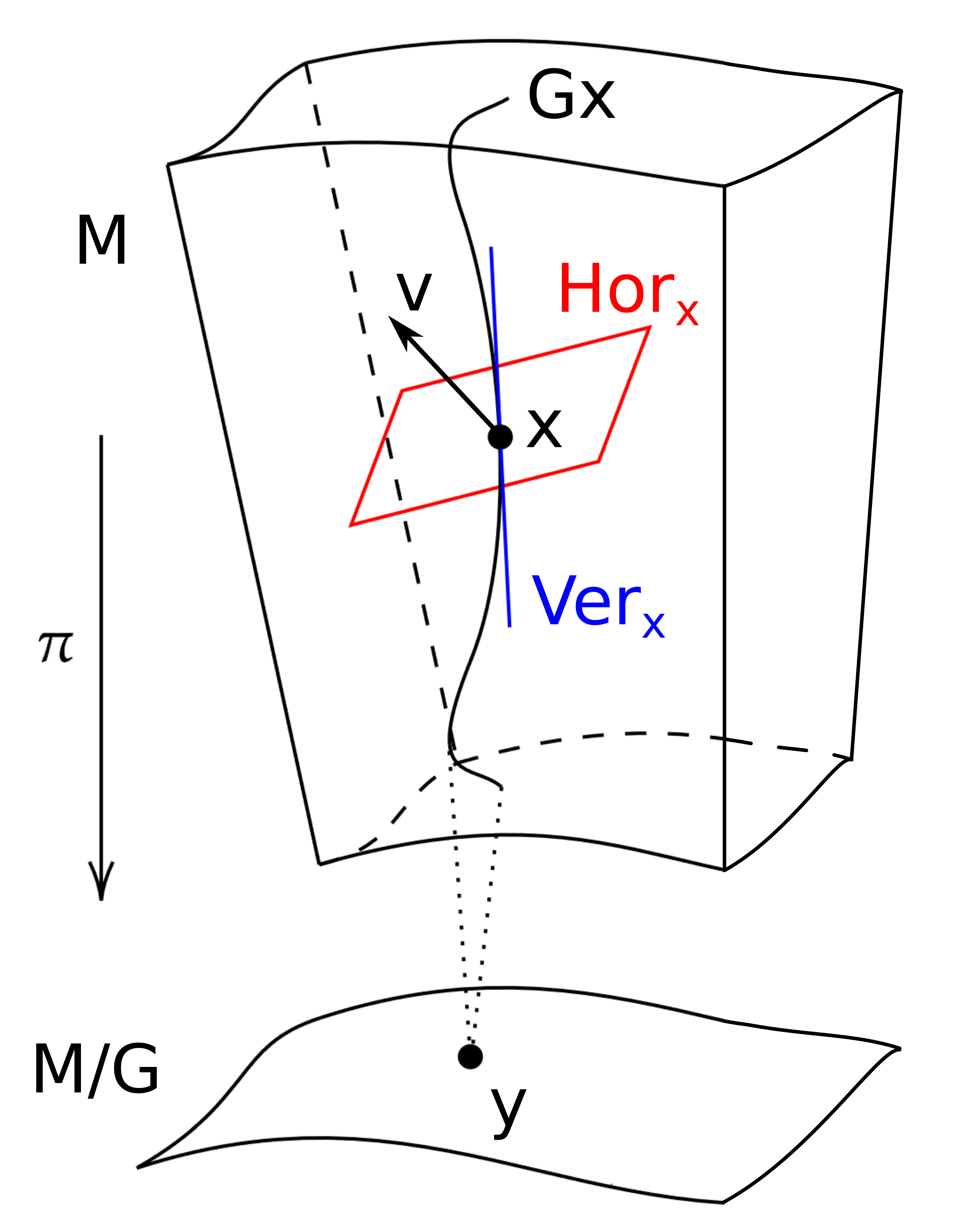}
\vspace{0.3cm}
\caption{Illustration of the principal bundle geometry including decomposition of $T_x M$ into horizontal $\mathrm{Hor}_x$ and vertical $\mathrm{Ver_x}$ spaces defined by a connection $\mathcal{A}_x$.}
\label{fig:bundle}
\end{figure}

Any connection can be equivalently characterized by its \emph{vertical and horizontal spaces}, namely $\mathrm{Ver_x} = \mathrm{Ker}T_x \pi$ and $\mathrm{Hor}_x = \mathrm{Ker}\mathcal{A}_x$. This induces the definition of a vertical and a horizontal components of a vector $v \in T_x M$, $\mathrm{Ver}(v) = \mathcal{A}_x(v)$ and $\mathrm{Hor}(v) = v - \mathcal{A}_x(v)$. As a result, we have a decomposition into sub-bundles $TM = \mathrm{Hor}(TM) \oplus \mathrm{Ver}(TM)$, see Figure~\ref{fig:bundle} for illustration. Notice that  $d_x \pi : \mathrm{Hor}_x \rightarrow T_{\pi(x)} (Q/G)$ is an isomorphism. \emph{The curvature} of $\mathcal{A}$ is denoted by $B$ and defined as $B_x: T_xM \oplus T_xM \rightarrow \mathfrak g$, $B_x(v,w) = d \mathcal{A}_x(\mathrm{Hor}(v), \mathrm{Ver}(w))$. 

\begin{remark}
In case of trivial bundle $P\times G \rightarrow P $, it is always possible to choose $\mathcal{A}$ as the trivial connection. The \emph{trivial connection} is defined by $\mathcal{A}_{(p,g)} =g^{-1}\circ dg$ and has zero curvature. 
\end{remark}

The \emph{adjoint bundle} to $M \rightarrow M/G$ is the vector bundle $\tilde{\mathfrak g} = (M \times \mathfrak g)/G$ over the base manifold $M/G$. The action of $G$ on $\mathfrak g$ is the adjoint action $\xi \mapsto g \xi g^{-1}$. Elements of $\tilde{\mathfrak g}$ are denoted by $[x,\xi]_G$ for $x\in M$ and $\xi \in \mathfrak{g}$. The connection $\mathcal{A}$ induces the following vector bundle isomorphisms
\begin{equation} \label{eq:identification.TM}
    \begin{aligned}
        \alpha_\mathcal{A} : TM/G &\rightarrow T(M/G) \oplus \tilde{\mathfrak g} \\
        [(x,v)]_G &\mapsto d_x \pi(x,v) \oplus [x,\mathcal{A}_x(v)]_G
    \end{aligned}
\end{equation}
and its pull back
\begin{equation} \label{eq:identification.T*M}
    \begin{aligned}
        \alpha^*_\mathcal{A} : T(M/G)^* \oplus \tilde{\mathfrak g}^* &\rightarrow T^*M/G.
    \end{aligned}
\end{equation}
where $\oplus$ is a Whitney sum and $\tilde{\mathfrak g}^*$ its adjoint space of $\tilde{\mathfrak g}$.

\subsection{Reduction of control systems with symmetries} \label{Sec:reduction}

Symmetry of \eqref{eq:OCP_original} induces \emph{invariance of the optimal control}. This permits to perform a reduction of the control system and obtain a reduced optimal control problem equivalent to the original one. Reduction procedure present in the following is mainly based on \cite{Ohsawa13}. First, let us introduce equivariant control systems.

\begin{definition}
Control system $\dot{x} = f(x,u)$ is \emph{equivariant} with respect to $G$-action on $M$ if 
$$d\Phi_g \circ f(x,u) = f(\Phi_g(x), u).$$
\end{definition}

Assume that the control system $\dot{x} = f(x,u)$ is equivariant with respect to the $G$-action. 
The equivariance of a control system induces the reduced control system defined on the quotient space $\bar{f} :  M/G\times \R^m \rightarrow TM/G$.
Let $\mathcal{A}$ be a principal connection of $M \rightarrow M/G$ and $\alpha_A$ defines the isomorphism between $TM/G$ and $T(M/G) \oplus \tilde{\mathfrak g}$ given by \eqref{eq:identification.TM}. 
The \emph{reduced system} on $TM/G$ can be decomposed in two parts, on $T(M/G)$ and on $\tilde{\mathfrak g}$ by identification
$$\alpha_A \circ \bar f = \bar f_{M/G} \oplus \bar f_{\tilde{\mathfrak g}}.
$$
The \emph{reduced dynamics} can then be obtained in the following way
$$\bar f_{M/G}([x]_G, u) = d \pi_M(x) \circ f(x,u) \quad \text{and} \quad f_{\tilde{\mathfrak g}}([x]_G,u) =[x,\mathcal{A}(x, f(x,u))]_G.
$$
The \emph{reduced equations} are defined in variables  $y = 
\pi_M (x) \in M/G$ and $\zeta = [x,\mathcal{A}(x,\dot x)]_G \in \tilde{\mathfrak g}$ by
\begin{equation} 
\label{eq:reduced.contro.syst}
    \dot {y}= f_{M/G}(y, u) \quad \text{and} \quad \zeta = f_{\tilde{\mathfrak g}}(y,u).  
\end{equation}

\begin{definition}
We say that an optimal control problem admits \emph{$G$-symmetry} if $f^0$ and $\varphi_T$ are invariant with respect to $G$, i.e.\ $f^0(\Phi_g(x),u) = f^0(x,u), \ \varphi_T(\Phi_g(x))=   \varphi_T(x)$, and $f$ is equivariant with respect to $G$.
\end{definition}

Consider \eqref{eq:OCP_original} which admits $G$-symmetry.
In this case we say that \eqref{eq:OCP_original} is $G$-symmetric. 
The symmetry induces the corresponding reduced problem on the quotient space $M/G$ defined as follows. 

\begin{definition}[Reduced optimal control problem] 
\label{def:reduced.OCP}

Find a control $u(t) \in \R^m$ and trajectories $y(t) \in M/G, \ \zeta(t) \in \tilde{\mathfrak g} = (M \times \mathfrak g)/G $ solutions of
\begin{equation} 
\label{eq:reduced_OCP}
\begin{aligned}
    &\min_{u}  & J(u) = & \int_{0}^{T}  f^0_{M/G}(y(t),u(t)) d t  \\
    &  \text{s.t.} & \dot{y}(t) = &~ f_{M/G}(y(t),u(t)), ~~ y(0) = y_0 \\
    &   &   \zeta(t) = &~ f_{\tilde{\mathfrak g}}(y(t),u(t)), \\
    &   &	0 =&~ \bar \varphi_T(y(T)).
\end{aligned}
\end{equation}
\end{definition}

Reconstruction of the solution of \eqref{eq:OCP_original} from solution of \eqref{eq:reduced_OCP} can be done in local trivialization $(y,g)$ of $M \rightarrow M/G$. using the reconstruction equation for $g$
$$ \dot{g}= g \,\xi(y,u), \qquad \xi = \zeta - \mathcal{A}_{(y, e)} (\dot{y}, 0). $$

\begin{theorem}[\cite{Ohsawa13}]  
\label{th:reduced.OCP.equiv}
If optimal control problem \eqref{eq:OCP_original} admits $G$-symmetry, then the reduced problem \eqref{eq:reduced_OCP} together with reconstruction equation is equivalent to  \eqref{eq:OCP_original} and leads to the same solutions.
\end{theorem}

\begin{remark} 
\label{rem:trivial.bundle}
In case of trivial bundle $\pi : P\times G \rightarrow P $, 
the reduced equations are defined similarly as in \eqref{eq:reduced_OCP} but on $TP \oplus {\mathfrak g}$. In this case, one can particularly distinguish the cyclic case which corresponds to the case where the group coordinates do not appear in the control system and in the cost function. This naturally makes the \eqref{eq:OCP_original} G-symmetric and happens especially when $G$ is abelian.
The reduced OCP is defined on $P$ and obtained by ignoring the dynamics in $g$, which will only be useful for the reconstruction. In this case, reduced OCP simply write 
\begin{equation} 
\label{eq:OCP_cyclic}
\begin{aligned}
&\min_{u}  & J(u) &= \int_{0}^{T} f^0(y(t),u(t)) d t, \\
&  \text{s.t.} & \dot{y}(t) &= f(y(t),u(t)), ~~ y(0) = y_0 \\
&   &	0 &= \bar \varphi_T(y(T)), 
\end{aligned}
\end{equation}
with the reconstruction equation $\dot g(t) = h(y(t),u(t)), \ g(0) = g_0$ for $g \in G$. The turnpike property for optimal control problems of this form have been studied in \cite{Tre:23}.
\end{remark}

\subsection{Mechanical systems with symmetries} \label{Sec:mechanical.syst}

In the mechanical system the state is given by a position $q$ in \emph{configuration space $Q$} and velocity $v = \dot q$ in the \emph{tangent space $T_q Q$}, so that $M = TQ$. For many mechanical systems the Lagrangian can be defined on $TQ$ by
\begin{equation} 
\label{eq:Lagrangian}
L(q,\dot q) = \frac{1}{2} \dot q^\top M(q) \dot q - V(q),    
\end{equation}
where the first term represent the kinetic energy and the second term stands for potential energy. In addition we assume that $M(q)$ is a positive definite matrix for any $q \in Q$. The movement of a mechanical system in presence of external controlled force $F : TQ \times \R^m \rightarrow T^*Q$ can be described by Lagrange-d'Alembert principle
\begin{equation} 
\label{eq:d'Alembert}
\delta \int_0^T L(q(t),\dot q(t))dt + \int_0^T F(q(t),\dot q(t),u(t)) \, \delta q \, dt = 0,  
\end{equation}
resulting in controlled Euler-Lagrange equations defined by
\begin{equation}
\label{eq:E-L}
\frac{d}{dt} \left(\frac{\partial}{\partial \dot q} L(q,\dot q) \right)  -    \frac{\partial}{\partial q} L(q,\dot q) = F(q(t),\dot q(t),u(t)).
\end{equation}
For a Lagrangian given by \eqref{eq:Lagrangian}, second derivative $\ddot{q}$ is explicitly given by
\begin{equation}
\label{eq:E-L.second.ord}
\ddot{q} = M(q)^{-1} \left( \frac{\partial}{\partial q}V(q) - \frac{1}{2}\left(\frac{\partial}{\partial q}M(q)\circ\dot{q} \right)\dot{q}\right) + M(q)^{-1} F(q(t),\dot q(t),u(t)).
\end{equation}
Controlled Euler-Lagrange equations can then be written as a first order control system with state variables $(q,v) \in TQ$ by setting $v = \dot{q}$
\begin{equation} 
\label{eq:E-L.first.order}
\begin{cases}
\dot{q} = v, \\
    \dot{v} = M(q)^{-1} \left( \frac{\partial}{\partial q}V(q) - \frac{1}{2}\left(\frac{\partial}{\partial q}M(q)\circ v \right)v\right) + M(q)^{-1} F(q(t),\dot q(t), u(t)).
\end{cases}
\end{equation}
For simplicity, right-hand side of \eqref{eq:E-L.first.order} will be denoted by $(v, l(q,v,u))$ in the following.

In the classical case of mechanical systems without forcing, the reduction procedure leads to a variational formulation on the quotient space and dates back to the seminal works \cite{RalphMarsden2021, CMR:2001, CMR:2001b}. In references \cite{BlochKrishMarsdenRatiu1996,deLeon2021}, reduction %
has been extend to the case of mechanical systems with uncontrolled external forcing, nonholonomic systems \cite{CMR:2001} and to the case of feedback controlled Lagrangian and Hamiltonian systems \cite{CM:2004}. In the latter case the notion of controlled Lagrangian systems takes a different form from \eqref{eq:d'Alembert} and assumes that the control is a function of $(q, \dot q)$. Next, we will show that the reduction appearing in the listed literature have a straightforward generalization to  the systems with controlled forcing in the form \eqref{eq:d'Alembert}.

The main idea behind the variational principle \eqref{eq:d'Alembert} is to consider how the action integral varies on the set of curves $q(t)$ to determine the curve of the least action. The principal bundle structure $Q \rightarrow Q/G$ implies the definition of the reduced set of curves $[q]_G(t)$ on $Q/G$. 
Following \cite{CMR:2001b}, any reduced curve can be split in two part using the identification 
of $TQ/G$ from \eqref{eq:identification.TM},  $[q]_G(t) = y(t) \oplus \bar \Omega(t)$ with $y(t) \in T(Q/G)$ and $\bar \Omega(t) = \left[q(t),\mathcal{A}(q(t),\dot{q}(t))\right]_G\in \tilde{\mathfrak g}$. Notice that the identification is defined using a connection $\mathcal{A}$.  
In local trivialization of $Q \rightarrow Q/G$ given by $Y \times G \rightarrow Y$, where $Y$ is an open subset of $\R^r$ with $r = \dim(Q/G)$ the connection can be defined as follows. For $(y,g,\dot y, \dot g) \in T_{(y,g)}(Y \times G)$, $\mathcal{A}(y,g,\dot y, \dot g)  = \mathrm{Ad}_g (\mathcal{A}_e(y) \dot y + \xi)$, where $\mathcal{A}_e(y) \dot y = \mathcal{A}_{(y,e)}(\dot y, 0)$ and $\xi = g^{-1}\dot g$. The isomorphism \eqref{eq:identification.TM} can be expressed as follows
$$\alpha_{\mathcal{A}}([y,g,\dot y, \dot g]_G) = (y,\dot y)\oplus (y, \Omega), 
\qquad \Omega = \mathcal{A}_e(y) \dot y + \xi.$$ Consider an equivariant forcing term $F(q,\dot q, u)$. It induces the reduced map $f(x, \dot x,  \Omega, u)$ from $\left( T(Q/G) \oplus \tilde{\mathfrak g} \right) \times \R^m$ to $T^*(Q/G) \oplus \tilde{\mathfrak g}^*$ defined using the connection and acting on the reduced space of curves. 

\begin{Lemma}
Let $F : TQ \times \R^m \rightarrow T^*Q$ be equivariant with respect to left action of $G$ on $Q$. Then, there exists a reduced map $f : \left( T(Q/G) \oplus \tilde{\mathfrak g} \right) \times \R^m \rightarrow T^*(Q/G) \oplus \tilde{\mathfrak g}^*$ such that
\begin{equation} 
\label{eq:equiv.forcing}
F(q, \dot q, u) \delta q = f(y, \dot y, \bar \Omega, u) (\delta y \oplus \bar \eta),
\end{equation}
where $(q, \dot q) \in TQ, \ u \in \R^m, \ (y, \dot y) \in T(Q/G), \ \bar \Omega = (y, \Omega) \in \tilde{\mathfrak g}$ and $\delta y$ is a horizontal variation of $q(t)$ and $\bar \eta = (y(t),\delta \Omega(t))$ is a vertical variation of $q(t)$. Moreover, \eqref{eq:equiv.forcing} can be decomposed in a vertical and a horizontal parts
$$ f_{Q/G}(y, \dot y, \Omega, u) \delta y \quad \text{and} \quad  f_{\tilde{\mathfrak g}}(y, \dot y, \Omega, u) \bar \eta.$$
\end{Lemma}  

\begin{proof}
In local trivialization $(y,g)$ of $Q$, the forcing term $F$ can be written $F(y,\dot y, g, \dot g, u) = F(y,\dot y, g, g(\Omega - \mathcal{A}_e(y) \dot y), u)$ and it decomposes into $F = F_{Q/G}dy + F_{G}dg$. By equivariance, it holds
\begin{multline*}
F(y,\dot y, e, \Omega - \mathcal{A}_e(y) \dot y, u)( \delta y \oplus \delta g) = \\ = F_{Q/G}(y,\dot y, e, \Omega - \mathcal{A}_e(y) \dot y, u)\delta y + F_{G}(y,\dot y, e, \Omega - \mathcal{A}_e(y) \dot y, u)(g^{-1}\delta g)
\end{multline*}
Then, using the relationship $\delta \Omega = \mathcal{A}_e(y) \delta y + g^{-1} \delta g$, we get
$$F( \delta y \oplus \delta g) = f_{Q/G}(y,\dot y, \Omega, u)\delta y + f_{\tilde{\mathfrak g}}(y, \dot y, \Omega, u) \delta \Omega,$$
where $f_{Q/G}, f_{\tilde{\mathfrak g}}$ are obtained from $F_{Q/G}, F_{G}$ using the connection $\mathcal{A}$.
\end{proof}
In the case of non controlled mechanical systems, \cite{BlochKrishMarsdenRatiu1996,CMR:2001,CM:2004} provide the formulas for vertical and horizontal reduced Euler-Lagrange operator, see for instance \cite[Theorem 3.4.1]{CMR:2001} or \cite[equation (2.9)]{CM:2004}. Combining these results we state the following proposition on symmetry reduction adapted to our formulation of the controlled Lagrangian system. 

\begin{proposition}
\label{Prop:mech.reduction}
Let $Q$ be a configuration space, $L : TQ \rightarrow \R$ a left invariant Lagrangian and $F: TQ \times \R^m \rightarrow T^*Q$ a force field equivariant with respect to the canonical left action of $G$ on $TQ$ and on $T^*Q$. Let $l : T(Q/G) \oplus \tilde{\mathfrak g} \rightarrow \R $ and $f : \left( T(Q/G) \oplus \tilde{\mathfrak g} \right) \times \R^m \rightarrow T^*(Q/G) \oplus \tilde{\mathfrak g}^*$ be respectively the restrictions of $L$ and $F$ to the corresponding quotient spaces. In local trivialization of $Q/G$, a curve $q(t) = (y(t),g(t)) \in Y \times G$. 
Then the following are equivalent
\begin{itemize}
\item $q(t)$ satisfies the Euler-Lagrange equations \eqref{eq:E-L.first.order} with forcing for $L$ on $Q$
\item the Lagrange-Poincaré equations with forcing, also known as Lagrange-d’Alembert-Poincaré equations, are valid with the vertical equation
\[
\frac{d}{dt} \frac{\partial}{\partial \Omega^b} l - \frac{\partial}{\partial \Omega^a} l \left( c^a_{db} \Omega^d - c^a_{db} \mathcal{A}^d_\alpha \dot y^\alpha \right) = f^b_{\tilde{\mathfrak g}}(y, \dot y, \Omega, u)
\]
and the horizontal equation
\[
\frac{d}{dt} \frac{\partial}{\partial \dot y^\alpha} l -\frac{\partial}{\partial y^\alpha} l + \frac{\partial}{\partial \Omega^a} l \left( B^a_{\beta \alpha} \dot y^\beta + c^a_{d b} \Omega^d \mathcal{A}^b_\alpha\right)  = f^\alpha_{Q/G}(y, \dot y, \Omega, u),
\]
where $B^a_{\beta \alpha}$ is the curvature of $\mathcal{A}_e = (\mathcal{A}^a_\alpha)$, $c^a_{b d}$ are the structure constants of the Lie algebra $\mathfrak{g}$, and $x = (x^\alpha)_{a = 1, \dots, \mathrm{dim} (\mathfrak{g})}, \ \Omega = (\Omega^a)_{\alpha = 1, \dots, \mathrm{dim}(Q/G)}$ are coordinates.
\end{itemize}
\end{proposition}

\begin{remark} \label{rk:for.rigid.body}
In case of trivial bundle structure, i.e.,  $Q = G \times Y$, it is possible to use the trivial connection with $\mathcal{A}_e = 0$. Consider the forcing term in the form $F(y,g,\dot y, \dot g, u) = g  f_G(y,\dot y, \xi,u)dg + f_W(y,\dot y, \xi,u)dy$ with $\xi \in \mathfrak{g}$. It is equivariant by construction. Reduced equations from Proposition~\ref{Prop:mech.reduction} are defined on $TY \times \mathfrak{g}$ and take a simplified form with the vertical part given by
$$ \frac{d}{dt} \frac{\partial}{\partial \xi} l(x, \dot x, \xi) - \mathrm{ad}^*_{\xi} \frac{\partial}{\partial \xi} l(x, \dot x,  \xi) =  f_{G}(x, \dot x, \xi, u)$$
and the horizontal part given by
    $$ \frac{d}{dt} \frac{\partial}{\partial \dot x} l(x, \dot x, \xi) -  \frac{\partial}{\partial x} l(x, \dot x, \xi)  = f_{W}(x, \dot x,  \xi, u),$$
where $(x, \dot x) \in TW$, $u \in \R^m, \ v \in \mathfrak{g}$ and $\mathrm{ad}^*_\xi \frac{\partial l}{\partial \xi}$ stands for the co-adjoint operator, $\mathrm{ad}^*_\xi : \mathfrak g^* \rightarrow \mathfrak g$, $\langle \mathrm{ad}^*_\xi \mu, \eta \rangle = \langle \mu, \mathrm{ad}_\xi \eta \rangle = \langle \mu, [\xi, \eta] \rangle$ with $[\cdot,\cdot]$ the Lie bracket. The co-adjoint operator $\mathrm{ad}^*_\xi \frac{\partial l}{\partial \xi}$ has the coordinate expression $c^a_{d b} \xi^d \frac{\partial l}{\partial \xi^a}$.
\end{remark}

\begin{remark}
\label{rem:cyclic.OCP}
Considering abelian group case with cyclic coordinates $(s,\theta) \in Q = Q/G\times G$ ($s$ and $\theta$ being resp. shape and cyclic variables), if the trivialization is global, we can apply Remark~\eqref{rem:trivial.bundle} and use trivial connection in the reduction procedure. \emph{Measured trim turnpike} in this particular case has been analyzed in \cite{FFOSW22}.
The reduced equations are defined for  $(s, v_s, v_{\theta}) \in T(Q/G)\times {\mathfrak g}$ as follows
\begin{equation} 
 \begin{cases}
\dot{s} = v_s, \\
\dot{v}_s = l_s(s,v_s,v_{\theta}, u), \\
\dot{v}_{\theta} = l_{\theta}(s,v_s,v_{\theta}, u).
\end{cases}
\end{equation}
\end{remark}

\subsection{Relation between optimal control and mechanical reductions}

Let us now establish a connection between reduction of variational equations in mechanics and reduction of optimal control problem. The case of interest is when the control system in \eqref{eq:OCP_original} is given by controlled Euler-Lagrange equations coming from a Lagrangian problem with symmetry. 

\begin{proposition} 
\label{Prop:EL.equivariance}
 If Lagrangian $L : TQ \rightarrow \R$ is invariant and a force field  $F: TQ \times \R^m \rightarrow T^*Q$ is equivariant with respect to the action of $G$ on $TQ$ and $T^*Q$, then \eqref{eq:E-L.first.order} is equivariant with respect to $G$-action.
\end{proposition}

\begin{proof}
Assume that $q(\cdot)$ is a critical point of \eqref{eq:d'Alembert}. Now consider $\Phi_g(q)(\cdot)$ for some $g \in G$. Because of the invariance of the Lagrangian, the first term in \eqref{eq:d'Alembert} is invariant with respect to $G$-action on $q(\cdot)$ and thus, its value at $\Phi_h(q)(\cdot)$ equals to the value at $q(\cdot)$. Let us now prove the invariance of the second term in \eqref{eq:d'Alembert}. From the equivariance of $F$ we have 
$$ F(\Phi_g(q), T_q \Phi_g(\dot q), u) \delta \Phi_g(q) =  T^*\Phi_{g^{-1}}F(q, \dot q, u) \, T\Phi_{g} \delta q = F(q, \dot q, u) \delta q. $$
We conclude that \eqref{eq:d'Alembert} is invariant with respect to $G$. This implies that if $q(\cdot)$ is a critical point of \eqref{eq:d'Alembert}, then  $\Phi_h(q)(\cdot)$ is also a critical point of \eqref{eq:d'Alembert} and satisfies \eqref{eq:E-L.second.ord}. As $q(\cdot)$ and $\Phi_g(q)(\cdot)$ satisfy \eqref{eq:E-L.second.ord} for any $g \in G$, 
then the following also holds
\begin{equation} 
\begin{cases}
\frac{d}{dt}\left( \Phi_g{q}\right) = d_q\Phi_g v, \\
\frac{d}{dt}  \left(d_q\Phi_g v\right) = l(\Phi_g(q),d_q\Phi_g v,u) = d_q\Phi_g f(q,v,u).
\end{cases}
\end{equation}
This means that the control system \eqref{eq:E-L.first.order} is equivariant with respect to $G$-action, conclusion then follows.
\end{proof}

By Proposition~\ref{Prop:EL.equivariance}, under invariance condition on the Lagrangian and force of a mechanical system, Euler-Lagrange equations are given by an equivariant control system, which permit to perform a reduction described in Section~\ref{Sec:reduction}.

\begin{corollary}
 If Lagrangian $L : TQ \rightarrow \R$ is invariant and a force field  $F: TQ \times \R^m \rightarrow T^*Q$ is equivariant with respect to the action of $G$ on $TQ$ and $T^*Q$, then there exists a reduced system \eqref{eq:reduced.contro.syst} defined on $T(TQ)/G$ equivalent to~\eqref{eq:E-L.first.order}. 
\end{corollary}
Notice that the reduction procedure requires a connection. In mechanical systems there is a natural connection which is induced by $M(q)$, which defines a Riemannian metric on $Q$. 
\begin{corollary}
The reduction of \eqref{eq:E-L.first.order} as a control system and the reduction based on the variational principle are equivalent. 
\end{corollary}
\begin{remark}
In the optimal control framework, the control $u$ can be expressed by a necessary optimality condition called Pontryagin's maximum principle as a function of state $x$ and co-state (or adjoint vector) $\lambda$. We refer to Section~\ref{Sec:turnpike} for more details. It was shown in \cite{Ohsawa13} that applying 
Pontryagin's maximum principle commutes with 
applying the
reduction procedure. By the equivalence of the two reductions, we can also conclude that the mechanical reduction in Proposition~\ref{Prop:mech.reduction} holds when $u(t) := u(x(t), \lambda(t))$. 
\end{remark}

\section{Turnpike and Trim turnpike}
\label{Sec:turnpike}

This section is devoted to the analysis of the turnpike property based on the reduction covered in the previous section. First, we consider the reduced problem on the quotient space, which allows us to use the classical \emph{static turnpike}. Then, we prove the convergence of optimal solution to \emph{trim turnpike}. 

\subsection{Static problem}

Based on the reduced form \eqref{eq:OCP_cyclic}, problem \eqref{eq:OCP_original} can be written as follows.
Find a control $u(t) \in \R^m$ and a trajectory $(y(t), \xi(t)) \in M/G \times \lieg $, with $\lieg = (M \times \mathfrak g)/G $, satisfying 
\begin{equation}
\label{eq:OCP'}
\tag{OCP'}
\begin{aligned}
& \min_{u} & J(u) &= \int_{0}^{T}  f^0(y(t),u(t)) d t,   \\
&  ~ \text{s.t.} & \dot{y}(t) &= f(y(t),u(t)), ~~ y(0) = y_0 \\
&   &  \dot g(t) &= g(t) \xi(y(t), u(t)), ~~ g(0) = g_0, \\
&   & 0 &= \varphi(y(T)),	
\end{aligned}
\end{equation}
where $y(t) \in M/G$, $g \in \Lieg$, $\xi(t) \in \tilde{\mathfrak g} = (M \times \mathfrak g)/G $ and where, by a slight abuse of notation, we set $f(y,u) = f_{M/G}(y, u) = d \pi(x) \circ f(x,u) : M/G \times \Omega \rightarrow M/G$ and $f^0(y,u) = f^0_{M/G}(y,u) : M/G \times \Omega \rightarrow \R$. Without loosing any generality, we assume final conditions in the form $y(T) = y_f$, for a given $y_f$. Note that the results presented below remain valid for more general initial and final conditions, as long as the conditions are independent from $g$.
We define the \emph{static problem}  associated with \eqref{eq:OCP'} by
\begin{equation}
\label{eq:SOCP'}
\min \, \{ f^0(\bar y, \bar u) \, | \, (\bar y, \bar u) \in M/G \times \R^m, ~~ f(\bar y,\bar u) = 0, ~~  \bar{g} (t) = \bar g_0 e^{t \xi(\bar y,\bar u)} \}.
\end{equation}
We assume that for any $T > 0$, \eqref{eq:OCP'} has a solution and that the static reduced problem 
%
%
has a unique solution $(\bar y, \bar u, \bar g(\cdot, \bar g_0))$ for any given $\bar g_0 \in G$.
The reduced problem corresponding to \eqref{eq:OCP'} is given by
\begin{equation}
\label{eq:ROCP}
\tag{ROCP}
\begin{aligned}
&\min_{u} & J(u) &= \int_{0}^{T}  f^0(y(t),u(t)) d t, \quad u(t) \in \Omega  \\
&~ \text{s.t.} & \dot{y}(t) &= f(y(t),u(t)),  \\
&	& y(0) &= y_0, ~~  y(T) = y_f.	
\end{aligned}
\end{equation}
Its associated static problem is defined by
\begin{equation}
\label{eq:static_reduced_pb_2}
\min \, \{ f^0(\bar y, \bar u) \, | \, (\bar y, \bar u) \in M/G \times \R^m, ~~ f(\bar y,\bar u) = 0 \}.
\end{equation}
We denote by $\mathcal{H}$ (resp.\ by $H$) the Hamiltonian corresponding to \eqref{eq:OCP'} (resp.\ to~\eqref{eq:ROCP}) defined by
\[
\begin{aligned}
\mathcal{H}(y,g,p_y,p_g,p^0,u) &= \langle p_y, f(y,u) \rangle + \langle p_g, g(t) \xi(y(t), u(t)) \rangle + p^0 f^0(y,u) ~~  \\
\text{and} ~~~ H(y,g,p^0,u) &= \langle p_y, f(y,u) \rangle + p^0 f^0(y,u).
\end{aligned}
\]
According to the Pontryagin's Maximum Principle (PMP) \cite{Pontryagin:1962}, if $(y(t),g(t), u(t))$ form an optimal solution of \eqref{eq:OCP'} for $t \in [0,T]$, then there exists a scalar $p^0 \leq 0$ and an absolutely continuous function $p = (p_y, p_g): [0,T] \rightarrow T^* (M/G \times \lieg) $ (called adjoint vector) satisfying $(p(\cdot), p^0) \neq (0,0)$, such that 
\begin{eqnarray}
\label{eq:pmp_conditions_y}
\dot{y}(t) &=& \nabla_{p_y} \mathcal{H}[t] , \quad \dot{p}_y(t) = - \nabla_y \mathcal{H}[t] , \\[0.5em]
\label{eq:pmp_conditions_g}
\dot{g}(t) &=& \nabla_{p_g} \mathcal{H}[t] , \quad \dot{p}_g(t) = - \nabla_g \mathcal{H}[t] , \\[0.5em]
\label{eq:pmp_conditions_u}
0 &=& \nabla_u \mathcal{H}[t]     
\end{eqnarray}
for almost every $t \in [0,T]$ and with $[t] = (y(t), g(t), p_y(t), p_g(t), p^0, u(t)) $. 
Moreover, transversality conditions for $p_g$ imply $p_g(T) = 0$.
Also, we assume the solution not to be abnormal, i.e.\ $p^0 \neq 0$ and normalize $p^0$ to $p^0 = -1$. 

\begin{Lemma} 
\label{lem:equiv}
Solution of \eqref{eq:SOCP'} projected to $(y,u)$ variables coincides with a solution of \eqref{eq:static_reduced_pb_2}.
\end{Lemma}

\begin{proof}
Equation on $\bar g$ in \eqref{eq:SOCP'} does not constitute a constraint in the optimization problem. Indeed, on the one hand, it defines the value of $\bar g(t)$ for any $t \in \R$ and it is well defined for any pair $(\bar y, \bar u)$ and, on the other hand, both $f^0$ and $f$ are independent from $\bar g$. Therefore, the relationship between \eqref{eq:SOCP'} and \eqref{eq:static_reduced_pb_2} follows from their construction.
\end{proof}

It follows from Lemma~\ref{lem:equiv} that a solution $(\bar y, \bar u, \bar g)$ of \eqref{eq:SOCP'} can be obtained as a solution $(\bar y, \bar u)$ of \eqref{eq:static_reduced_pb_2} and $\bar g$ has the following form
\begin{equation}
\label{eq:anal_expression_of_group_turnpike}
    \bar{g} (t) = \bar g_0 e^{t \xi(\bar y,\bar u)}. 
\end{equation}
%

\subsection{Exponential turnpike for the reduced problem}

We first focus on the reduced problem \eqref{eq:ROCP} on the quotient space with $y \in M/G$ and its associated static problem \eqref{eq:static_reduced_pb_2}.
Recall that the corresponding Hamiltonian is given by
$H(y,p_y,u) = \langle p_y, f(y,u) \rangle - f^0(y,u)$ and from PMP we have
\begin{equation}
\label{eq:Ham-syst}
\dot{y}(t) = \nabla_{p_y} H(y,p_y,u) , \quad \dot{p}_y(t) = - \nabla_y H(y,p_y,u),  \quad \nabla_u H(y,p_y,u) = 0.     
\end{equation}
Next, we introduce the following notations
\[
\begin{aligned}
H_{yy} &= \frac{\partial^2 H}{\partial y^2}(\bar y, \bar p_y, \bar u), \quad H_{yp_y} = \frac{\partial^2 H}{\partial y \partial p_y}(\bar y, \bar p_y, \bar u), \\
H_{yu} &= \frac{\partial^2 H}{\partial y \partial u}(\bar y, \bar p_y, \bar u), \quad H_{p_yu} = \frac{\partial^2 H}{\partial p_y \partial u}(\bar y, \bar p_y, \bar u), \quad H_{uu} = \frac{\partial^2 H}{\partial^2 u}(\bar y, \bar p_y, \bar u),
\end{aligned}
\]
and we denote 
\[
A = H_{p_yy} - H_{p_yu}H_{uu}^{-1} H_{uy}, \quad B = H_{p_yu}, \quad W =-H_{yy} + H_{yu}H_{uu}^{-1} H_{uy}.
\]

\begin{definition}[Local exponential turnpike property] 

Optimal control problem \eqref{eq:ROCP} admits the \emph{local exponential turnpike property} if at least one of its solutions satisfies the following property. There exist positive constants $\eps, \mu, C$ and $T_0$ such that if
\[
\|y(0) - \bar y \| + \| \bar p_y \|   \leq \eps,
\]
then for any $T> T_0$, the following holds
\[
\| y(t) - \bar{y}  \|  + \| p_y(t) - \bar p_y \| + \| u(t) - \bar{u}  \| \leq C \left( e^{-\mu t} + e^{-\mu ( T-t)}\right), \qquad \forall \, t \in  [ 0,  T ].
\]
\end{definition}

Necessary conditions for the local exponential turnpike can be stated as follows (see  \cite{Tre_Zua:2014} for more details and for the proof).

\begin{theorem}[Exponential turnpike theorem \cite{Tre_Zua:2014}]
\label{th:exp-turnpike}
Let \eqref{eq:ROCP} have at least one solution for any $T > 0 $ and the following conditions hold,
\begin{enumerate}
\item the matrix $H_{uu}^{-1}$ is symmetric negative definite,
\item the matrix $W$ is symmetric positive definite,
\item the pair $(A, B)$ satisfies the Kalman condition, that is
\[
rank(B, AB, \ldots, A^{n-1} B) = n,
\]
\end{enumerate} 
then the local exponential turnpike property holds for \eqref{eq:ROCP}.
\end{theorem}

Assumptions 1. to 3. from Theorem~\ref{th:exp-turnpike}  insure the hyperbolicity of the state-adjoint system in PMP. More precisely, the matrix 
\[
M = \begin{pmatrix}
    A & -BH_{uu}^{-1}B^T \\
    W & -A^T
\end{pmatrix},
\]
obtained by linearizing the Hamiltonian system given in \eqref{eq:Ham-syst} is Hamiltonian and hyperbolic. For a Hamiltonian matrix to be hyperbolic means that a half of the eigenvalues is with strictly positive real part and the other half is with strictly negative real part.

\subsection{Symmetry, hyperbolicity and eigenvalues}
\label{sec:HSE}

\emph{The exponential turnpike property} is essentially based on the existence of the static problem solutions and the hyperbolicity property, which are not satisfied in the context of symmetric optimal control problems, as it is shown in Proposition~\ref{prop:symmetry-eigenvalues}. Nevertheless, it is still possible to analyze the turnpikes of such OCPs. This can be done by 
first reducing the problem to the quotient space. The reduced problem may have the necessary hyperbolicity property, and therefore, admit classical turnpikes described in the previous section. Before discussing the relation between symmetry and hyperbolicity, let us first come back to the reduction of the control system.

\begin{definition}
Reduction of \eqref{eq:OCP_original} is said to be \emph{maximal} (or the reduced system, \emph{minimal}) if the reduced system \eqref{eq:ROCP} does not admits a symmetry.
\end{definition}

\begin{Lemma}[Maximal reduction]
\label{lem:max-reduction}
If \eqref{eq:OCP_original} has a $G$-symmetry, then it is always possible to apply a maximal reduction to obtain a minimal reduced system.  
\end{Lemma}

The proof of Lemma~\ref{lem:max-reduction} is based on the fact that it is always possible to define a product group of all the symmetries of the system and perform a reduction with respect to it. 
This simple and intuitive result is important, since it allows us to reduce \eqref{eq:OCP_original} in such a way that the reduced problem \eqref{eq:ROCP} does not admit any symmetries. 

\begin{proposition}
\label{prop:symmetry-eigenvalues}
Let $\dot z = h(z) $ be a Hamiltonian system equivariant with respect to a free $G$-action on $M$, with $dim(G) = n$, then one of the following two situations occurs.
\begin{enumerate}
\item $\dot z = h(z) $ has no critical point,
\item for any critical point $\bar z \in M$, $Dh(\bar z)$ has the eigenvalue $\lambda = 0$ of multiplicity $n$, $n \in \N^*$.
\end{enumerate}
\end{proposition}

\begin{proof}
The proof is reduced to showing that situation 2 holds whenever system $\dot z = h(z) $ admits a critical point. Let  $\bar z \in M$ be a critical point of the system. Assume that $G$ is parameterized by $\alpha = (\alpha_1,\ldots, \alpha_n) \in \R^n$.\footnote{For some isomorphic function $\varphi : \R^n \rightarrow G$, $\varphi(\alpha) = g_\alpha$.} Let $g_\alpha \in G$. Equivariance of $\dot z = h(z) $ gives $h(g_\alpha \bar z) = g_\alpha h(\bar z) = 0.$
Differentiating this equation with respect to $\alpha_i, \, i=1,\ldots,n$, and evaluating it at zero, one obtains
\[
0 = \left.\frac{\diff }{\diff \alpha_i}\right|_{\alpha = 0} h(g_\alpha \bar z) = Dh(g_0 \bar z) \left.\frac{\diff }{\diff \alpha_i}\right|_{\alpha = 0} g_\alpha \bar z = Dh(\bar z) v_i,
\]
with $v_i =\left. \frac{\diff }{\diff \alpha_i} \right|_{\alpha = 0} g_\alpha \bar z \in ker Dh(\bar z) $. Moreover, since $G$ acts freely on M (i.e.~$g x = x \Rightarrow g = id_G$), we can choose $g_\alpha$ such that $v_i, ~ i=1,\ldots,n$ are independent. The conclusion then follows.
\end{proof}

If \eqref{eq:OCP_original} admits a $G$-symmetry, then the corresponding Hamiltonian system obtained from PMP is equivariant with respect to $G$-action \cite{Ohsawa13}. According to Proposition~\ref{prop:symmetry-eigenvalues}, such a system can not have hyperbolicity property. Thus, the turnpike analysis requires the maximal reduction of the optimal control problem.

\begin{remark}
If the hyperbolicity property is satisfied by the reduced  Hamiltonian matrix, then the exponential turnpike theorem will hold for the reduced system. This is the case, for example, in the Kepler problem, which will be discussed in Section~\ref{sec:examples}. However, even the maximal reduction of the system does not guarantee the hyperbolicity property for the reduced system. In fact, the reduced Hamiltonian matrix may have zero eigenvalue not linked to symmetries. In this case, a further analysis is required. The Rigid Body problem, discuss in Section~\ref{sec:examples}, will highlight this aspect. 
\end{remark}

\subsection{Trim Turnpike for the full system with free final conditions}
When the reduced system satisfies Theorem~\ref{th:exp-turnpike}, the exponential turnpike property holds for the reduced variables, namely $y \in M/G$ and $u$. The next step is to analyze the convergence properties of the group component $g \in G$ of the state. Our goal is to prove the exponential convergence of optimal solutions to trims.

\begin{definition}[Local exponential trim turnpike property] \label{def:trim.turnpike}
Optimal control problem~\eqref{eq:OCP'} admits \emph{local exponential trim turnpike property} if at least one of its solutions satisfies the following property. There exist $\bar g_{T/2} \in G$, positive constants $\eps, \mu, C$ and $T_0$ such that if
\[
\|y(0) - \bar y \| + \| \bar p_y \|   \leq \eps
\]
then for any $T > T_0$ the following holds
\begin{eqnarray}
\label{eq:trim-turnpike-1}
\| y(t) - \bar{y}  \| + \| p_y(t) - \bar p_y \| + \| u(t) - \bar{u}  \| & \leq & C \left( e^{-\mu t} + e^{-\mu ( T-t)}\right), ~~ t \in  [0,  T ], \\
\label{eq:trim-turnpike-2}
\| g(t, g_0) - \bar{g} (t, \bar g_{T/2})  \| & \leq &  C  \left( e^{-\mu t} + e^{-\mu ( T-t)}\right), ~ t \in  [0, T],
\end{eqnarray}
where $\bar g(\cdot, \bar g_{T/2})$ is a trim satisfying $\bar g(\frac{T}{2}, \bar g_{T/2}) = \bar g_{T/2}$ and $\dot{\bar g}(t) = \bar g(t) \xi(\bar y, \bar u)$ with $\xi$ defined in \eqref{eq:ROCP}.  
\end{definition}

The construction of the trim $\bar g$ in Definition~\ref{def:trim.turnpike} as well as the conditions for the exponential trim turnpike to hold are given in the following theorem.

\begin{theorem}
\label{th:trim-turnpike}
Assume that the following assumptions hold for the problem \eqref{eq:OCP'}:
\begin{enumerate}
\item \eqref{eq:OCP'} has at least one solution for any $T > 0 $;
\item assumptions 1), 2), and 3) of Theorem~\ref{th:exp-turnpike} are satisfied for \eqref{eq:ROCP};
\item quotient spaces $M/G$ and $\tilde{\mathfrak g} = (M \times \mathfrak g)/G $ are compact;
\item function $\xi$ is Lipschitz on $M/G\times \Omega$.
\end{enumerate}
Then \eqref{eq:OCP'} has the \emph{local exponential trim turnpike property}. 
\end{theorem}

\begin{proof}
The estimate \eqref{eq:trim-turnpike-1} comes from Theorem~\ref{th:exp-turnpike}. We are left to prove the convergence of $g$. 
Let us define
$\delta g(t) = g(t) - \bar{g} (t)$ for $t \in [0, T]$, where $g$ is the group component of the state and $\bar{g}$ is a trim, which satisfies the following equation
$$ 
\dot{\bar g}(t) = \bar g(t) \xi(\bar y, \bar u),
$$
with $(\bar y, \bar u)$ solution of \eqref{eq:static_reduced_pb_2}. From \eqref{eq:OCP'} and \eqref{eq:anal_expression_of_group_turnpike}, we have
\[
\begin{aligned}
\delta \dot{ g}(t) &= \dot{g}(t) - \dot{\bar g}(t) \\
&= g(t) \xi(y(t), u(t)) - \bar g(t) \xi(\bar y, \bar u) \\
&= g(t) \xi(y(t), u(t)) - g(t) \xi(\bar y, \bar u) + g(t) \xi(\bar y, \bar u)- \bar g(t) \xi(\bar y, \bar u) \\
&= \delta g(t) \xi(\bar y, \bar u) + g(t) \left( \xi(y(t), u(t)) - \xi(\bar y, \bar u) \right).
\end{aligned}
\]
Integrating this relation for $t \in [T/2, T]$, we obtain 
\[
\begin{aligned}
\delta g(t) &= \delta g(T/2) e^{t\xi(\bar y, \bar u)} + \int_{T/2}^{T/2 + \tau} g(s) \left(\xi(y(s), u(s)) - \xi(\bar y, \bar u) \right) e^{(t-s)\xi(\bar y, \bar u)} ds, 
\end{aligned}
\]
where $\tau = t - T/2$. This together with the fact that $\tilde{\mathfrak g} = (M \times \mathfrak g)/G $ is compact imply
\[
\begin{aligned}
\norm{\delta g(t)} &\leq \norm{\delta g(T/2)}e^{t\xi(\bar y, \bar u)} +  \int_{T/2}^{T/2 + \tau} \norm{g(s) }\norm{\left( \xi(y(s), u(s)) - \xi(\bar y, \bar u) \right) } e^{\left(t-s \right)\xi(\bar y, \bar u)} ds \\
&\leq C_1 \norm{\delta g(T/2)} + C_2 \int_{T/2}^{T/2 + \tau} \norm{\left( \xi(y(s), u(s)) - \xi(\bar y, \bar u) \right) } ds. 
\end{aligned}
\]
Now we use assumptions 2 and 4 of the Theorem, namely, the property of $\xi(\cdot)$ to be Lipschitz and the local exponential turnpike property of the reduced system. This leads us to
\[
\begin{aligned}
\norm{\delta g(t)} &\leq C_1 \norm{\delta g(T/2)} + C_3 \int_{T/2}^{T/2 + \tau} \norm{\left( \xi(y(s), u(s)) - \xi(\bar y, \bar u) \right) } ds \\
&\leq C \left(\norm{\delta g(T/2)} + \int_{T/2}^{T/2 + \tau} \left(e^{-\mu s} + e^{-\mu(T-s)}\right) ds \right)  \\
&\leq C \left(\norm{\delta g(T/2)} + \frac{1}{\mu}\left(e^{-\mu t} + e^{-\mu(T - t)} \right) \right), 
\end{aligned}
\]
where the last inequality is obtained by substituting back $\tau = t -T/2$. Let us now define the value of $\bar g(t)$ at $t = T/2$ by $\bar g(T/2) =  g(T/2)$. This implies $\delta g(T/2) = 0$ and
$$ 
\norm{\delta g(t)} \leq \tilde C \left( e^{-\mu (\frac{T}{2}-t)} + e^{-\mu(\frac{T}{2}+t)} \right).
$$
The same inequality holds for $t \in [0,T/2]$. This ends the proof.
\end{proof}

\begin{remark}
In \cite{Tre:23}, the author obtained a \emph{linear turnpike} for a class of optimal control problems, which can be interpreted as problems admitting a cyclic symmetry, i.e., in the form \eqref{eq:OCP_bis}, but with additional boundary conditions on the group variable. Notice first, that the class of symmetric optimal control problems that we treat is more general and has a strong motivation in mechanical systems, which usually admit conservation laws. Considering no final conditions on the group variable in our case simplifies the proof of the exponential turnpike property. More general boundary conditions will be considered in the future research.  
\end{remark}

\begin{remark}
Although the exact trim turnpike is only determined a posteriori in Theorem~\ref{th:trim-turnpike}, it is still useful from the numerical point of view. Indeed, for implementation purposes, only the static turnpike is needed to determine the optimal solution. On the other hand, this theorem allows us to better characterize the behavior of the solution and can therefore be used to study the behavior of the solution on an infinity time horizon, as it was done in \cite{GKW2017}, for instance. 
\end{remark}

\section{Examples}
\label{sec:examples}

We focus on three examples with different types of symmetries.
Our first example is the Kepler problem. The symmetry group in this case is abelian. It was already considered in \cite{FFOSW22}, where the turnpike was shown in numerical results. In our consideration, we identify the concrete trim toward which the solution converge with exponential rate following Theorem~\ref{th:trim-turnpike}.

The second example is an optimal control problem for rotations of a rigid body. The problem admits a non-abelian symmetry. We show that it fulfills Theorem~\ref{th:trim-turnpike} and identify the trim turnpike in numerical simulations.

Finally, we consider the example of a rigid body equipped with rotors. This system has a group product as its symmetry. Moreover, the corresponding reduced 
optimal control problem does not satisfy the hyperbolicity assumptions of Theorem~\ref{th:exp-turnpike}, but still admits the trim turnpike (in state and control) as we show in numerical results.

\subsection{Kepler problem}

The two-body problem, also known as Kepler problem, is to determine the motion of two point particles in $\R^2$ with masses $m_1$, $m_2\, > 0$.  
Polar coordinates are used to describe the motion of the second body relatively to the position of the first body via the radius $s \in \R_{>0}$ and angle $\theta \in [0,2\pi)$. We make a simplifying assumption $m_1 = 1$. Other parameters of the problem are gravitational constant $\gamma \in \mathbb{R}$ and $k = \gamma m_2 \in \R$. Lagrangian of the Kepler problem reads
\[
L(q,\dot q) = \frac{1}{2} \dot q^T \dot q + \frac{k}{\sqrt{q_1^2+q_2^2}},
\]
%
%
%
with $q = (s,\theta)^T$, $\dot q = (v_s,v_\theta)^T \in \R^2$.
Optimal control of this system has already been discussed in \cite{FFOSW22} where the authors show the existence of trim turnpikes for different cost functions in numerical simulations. Our goal is to show exponential convergence toward a trim turnpike, as formulated in Theorem~\ref{th:trim-turnpike}.
%
%
%
The Euler-Lagrange equations
\[
\frac{\diff}{\diff t}\frac{\partial}{\partial \dot q}L(q,\dot q) - \frac{\partial}{\partial q}L(q,\dot q) = f(u),
\]
with 
forcing term $f : \R^2 \rightarrow T^*M$, 
read 
\begin{align}
\label{eq:Kepler-reduced-syst-1}
    & \quad \dot \theta ~~~ = v_\theta, \\
    \label{eq:Kepler-reduced-syst}
    & \begin{cases}
    \dot s \;~ = v_s, \\
    \dot v_s = sv_\theta^2 - \frac{k}{s^2} + \frac{1}{m_2}f_s(u), \\
    \dot v_\theta = -\frac{2}{s}v_\theta v_s + \frac{1}{m_2s^2}f_\theta(u),
    \end{cases}
\end{align}
with the control $u=(u_s,u_\theta)^T \in \R^2$ and $f(u) = (f_s(u), f_\theta(u)) = (u_s,u_\theta)$. Let us denote the state of subsystem \eqref{eq:Kepler-reduced-syst} by $y=(s,v_s, v_\theta)$. 
\begin{remark}
In case, for example, where $f_s(u) = u_s = 0$,  \eqref{eq:Kepler-reduced-syst-1}-\eqref{eq:Kepler-reduced-syst} does not admit any static point. This corresponds to Scenario 1 in Proposition~\ref{prop:symmetry-eigenvalues}.
\end{remark}

Now, Let $(\bar y, \bar u)
= ((\bar s, \bar v_s, \bar v_\theta), (\bar u_s, \bar u_\theta))$ be a critical point of system \eqref{eq:Kepler-reduced-syst}. Fixing $\bar u = (0,0)$, one has 
\[
\bar v_s = 0, \quad \bar v_\theta = \sqrt{\frac{k}{\bar s^3}}, \quad \text{with} \quad \bar s \in \R_+^*.
\]
Moreover, we consider the following cost function
\[
f^0(y,u) = \frac{1}{2} \left( \|y-\bar y \|^2_2 + \| u-\bar u\|^2_2\right). 
\]
Solution of the corresponding optimal control problem can be obtained by first solving the reduced problem
\begin{eqnarray}
&& \min_{u} ~~ J(u) = \int_{0}^{T}  f^0(y(t),u(t)) d t,  \nonumber  \\
&& \text{subject to} ~~ \eqref{eq:Kepler-reduced-syst}, \nonumber  \\
&& y(0) = y_0 = \left(\bar s, \, 0, \, \sqrt{\frac{k}{\bar s^3}} \right), \nonumber 	
\end{eqnarray}
and then reconstructing $\theta$ from \eqref{eq:Kepler-reduced-syst-1}.

\begin{figure}[ht] 
\centering
\includegraphics[width=0.6\textwidth]{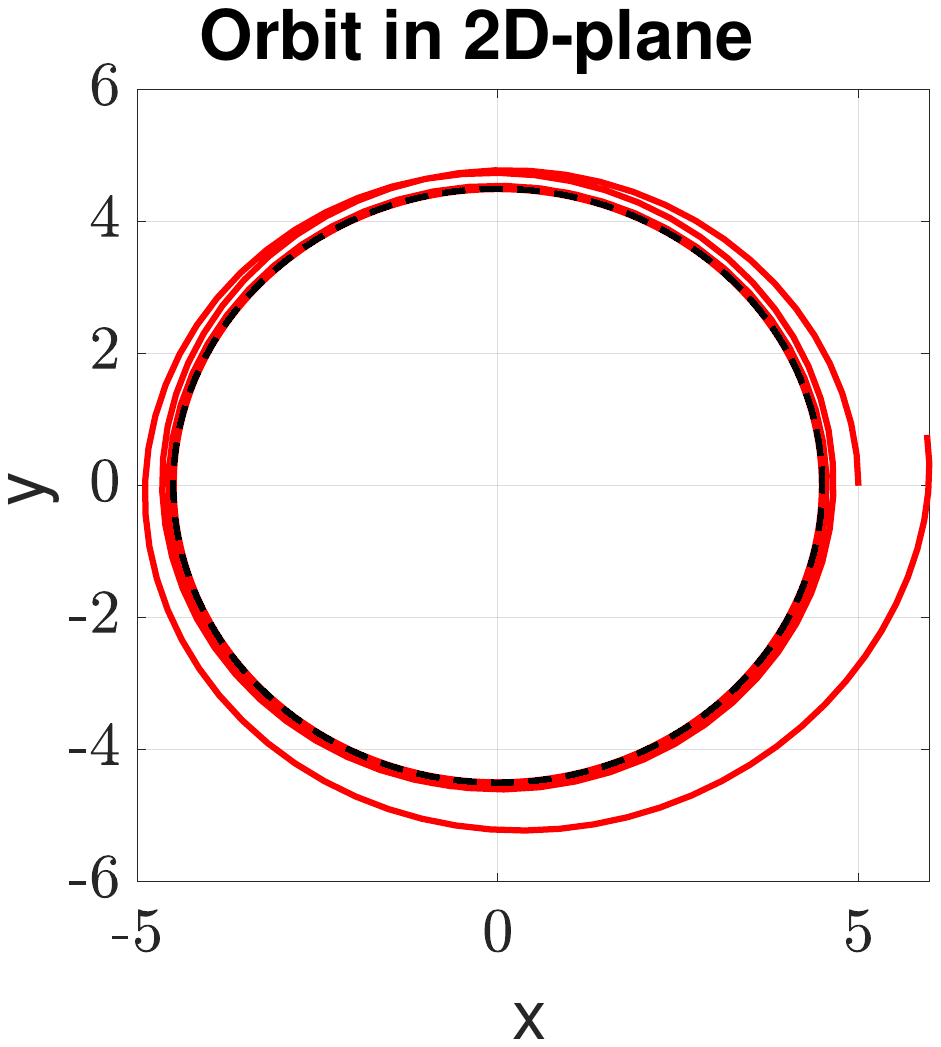}
\vspace{-2cm}
\caption{Illustration of the convergence of a solution toward a trim for the Kepler problem in 2D-plane. Tick red is used for optimal solution with turnpike property, dashed black for associated trim.}
\label{fig:Kepler-example-1}
\end{figure}
\begin{figure}[ht]
\centering
\def\size{0.33}
\def\Size{0.21}
\def\sizeh{-0.8}
\includegraphics[width=\size\textwidth]{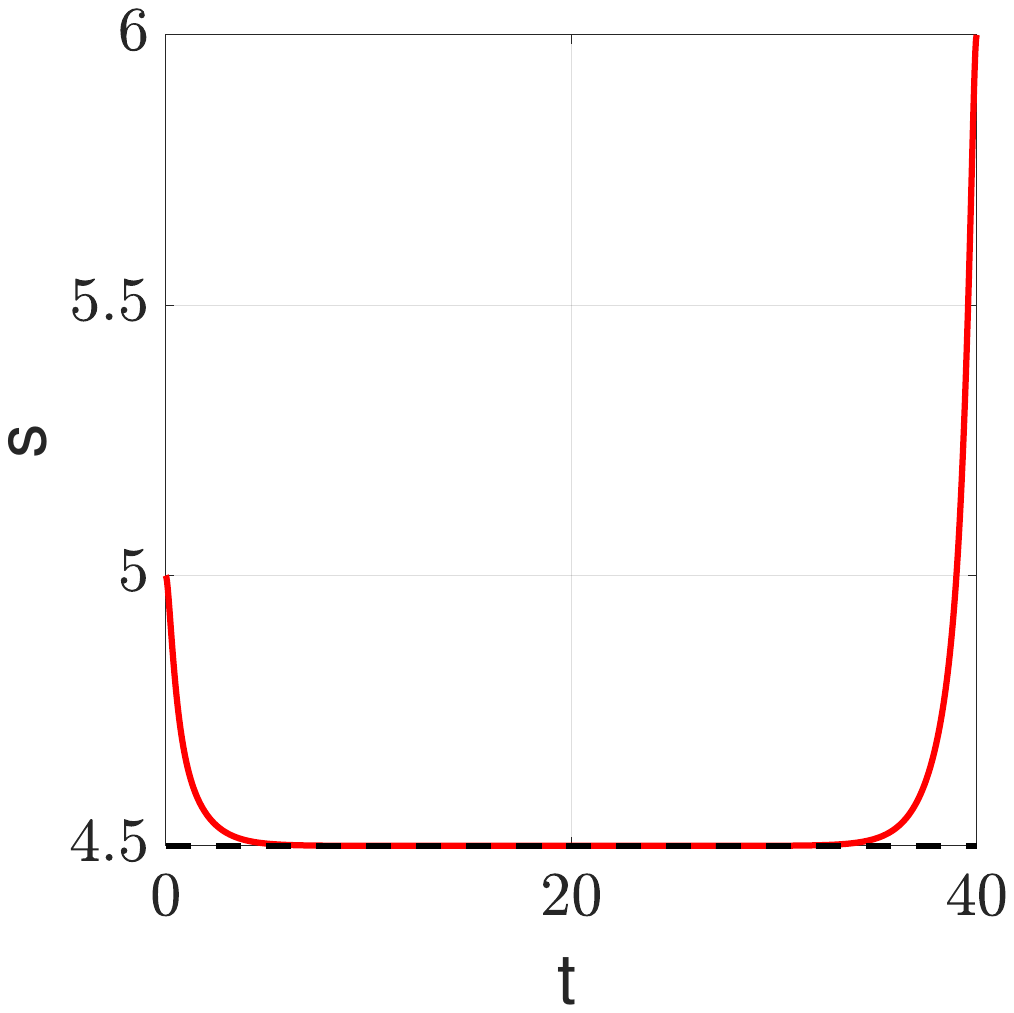}
\hspace{\sizeh cm}
\includegraphics[width=\size\textwidth]{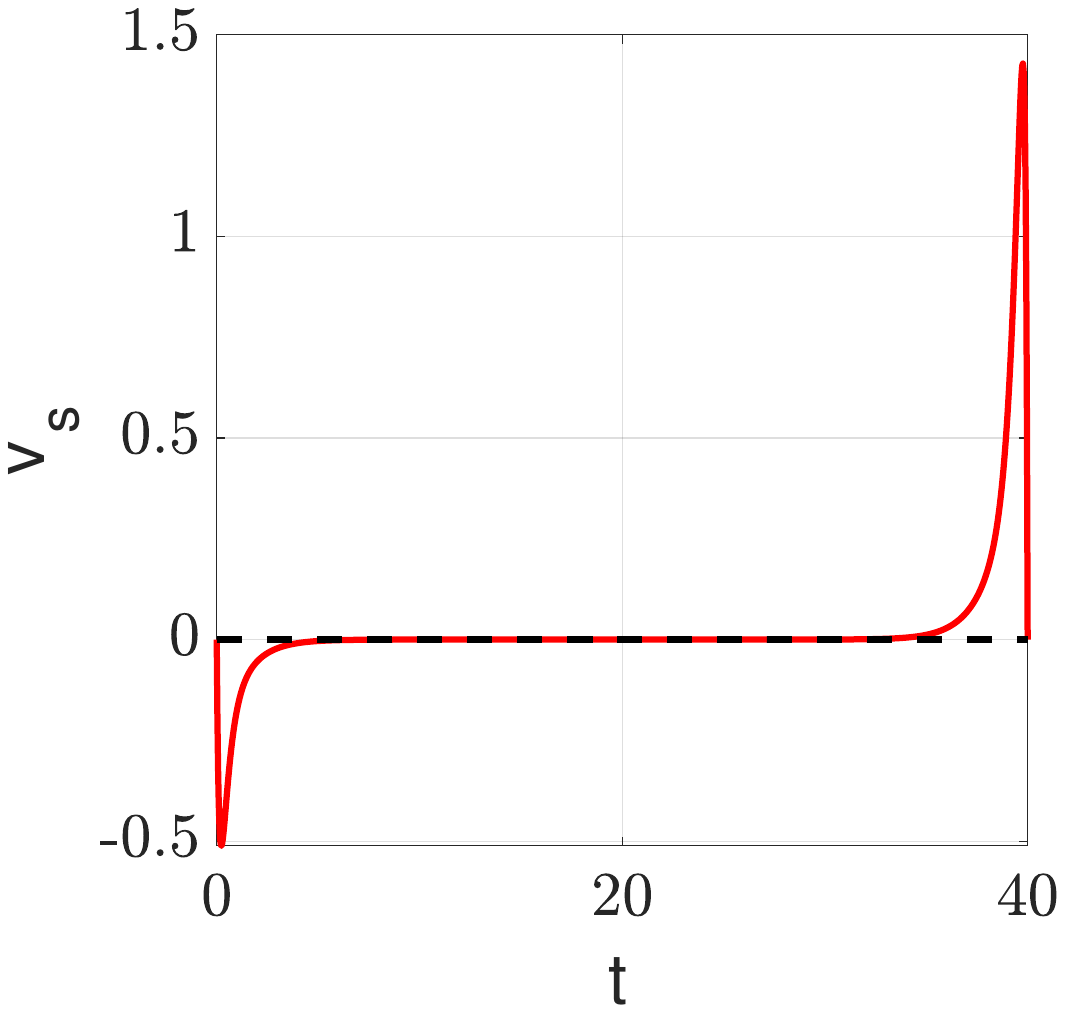}
\hspace{\sizeh cm}
\vspace{-2.cm}
\includegraphics[width=\size\textwidth]{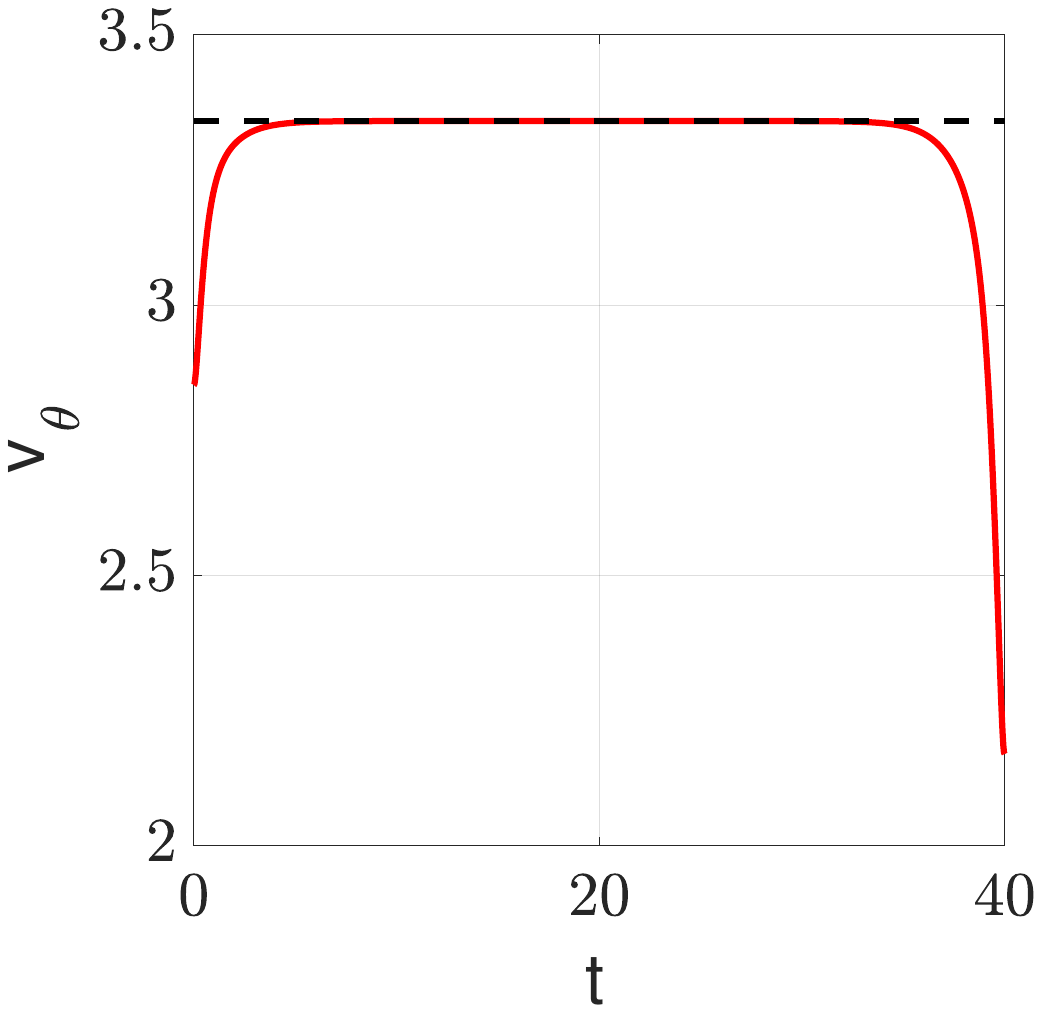}
\vspace{-1cm}
\includegraphics[width=\size\textwidth]{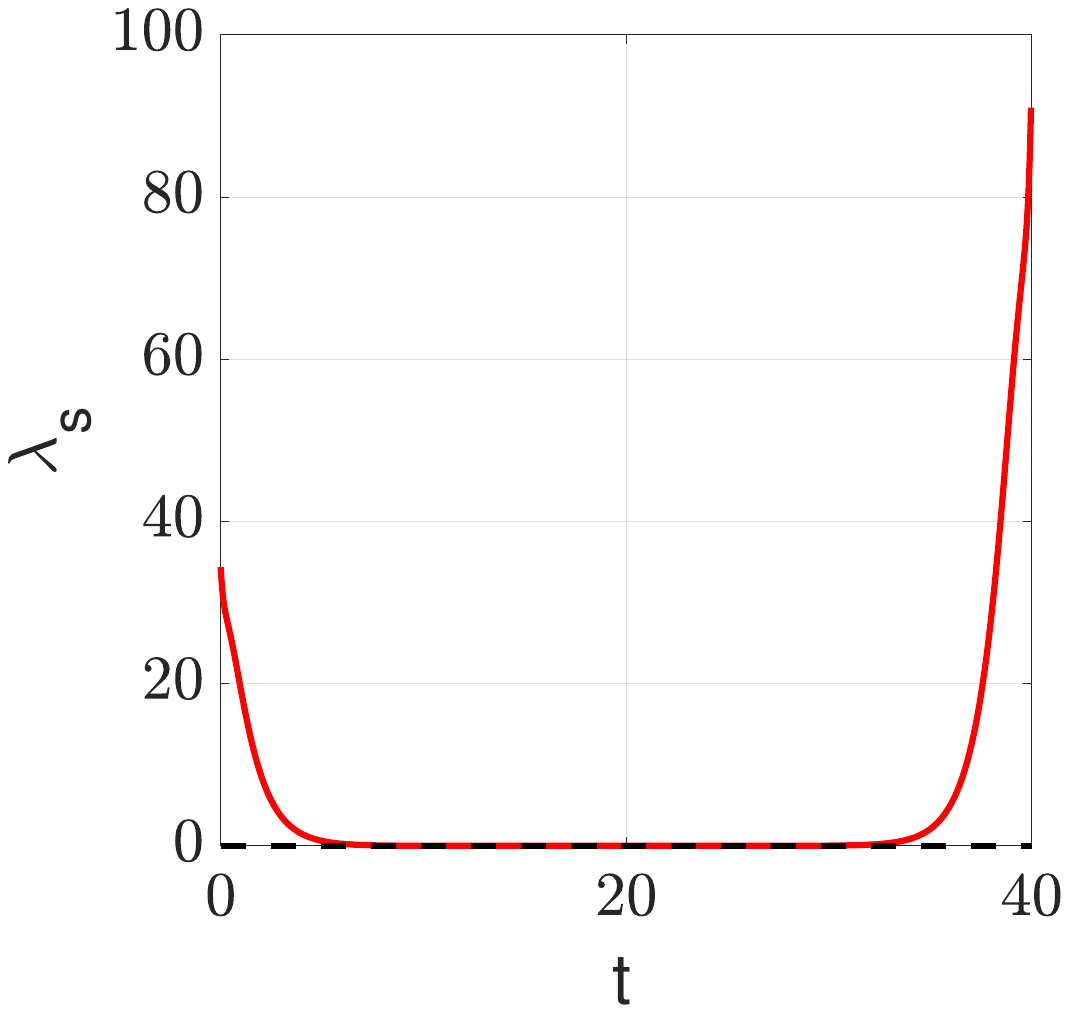}
\hspace{\sizeh cm}
\includegraphics[width=\size\textwidth]{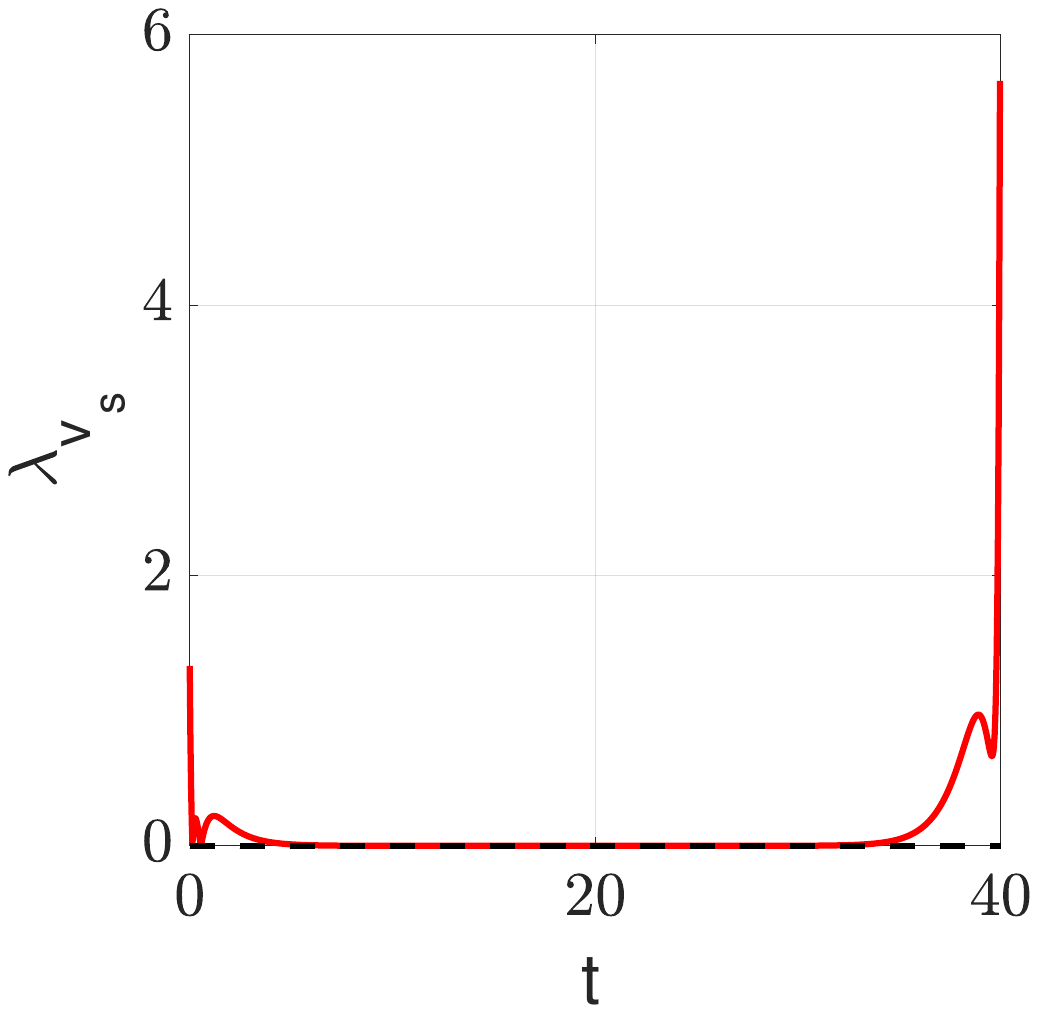}
\hspace{\sizeh cm}
\vspace{-1cm}
\includegraphics[width=\size\textwidth]{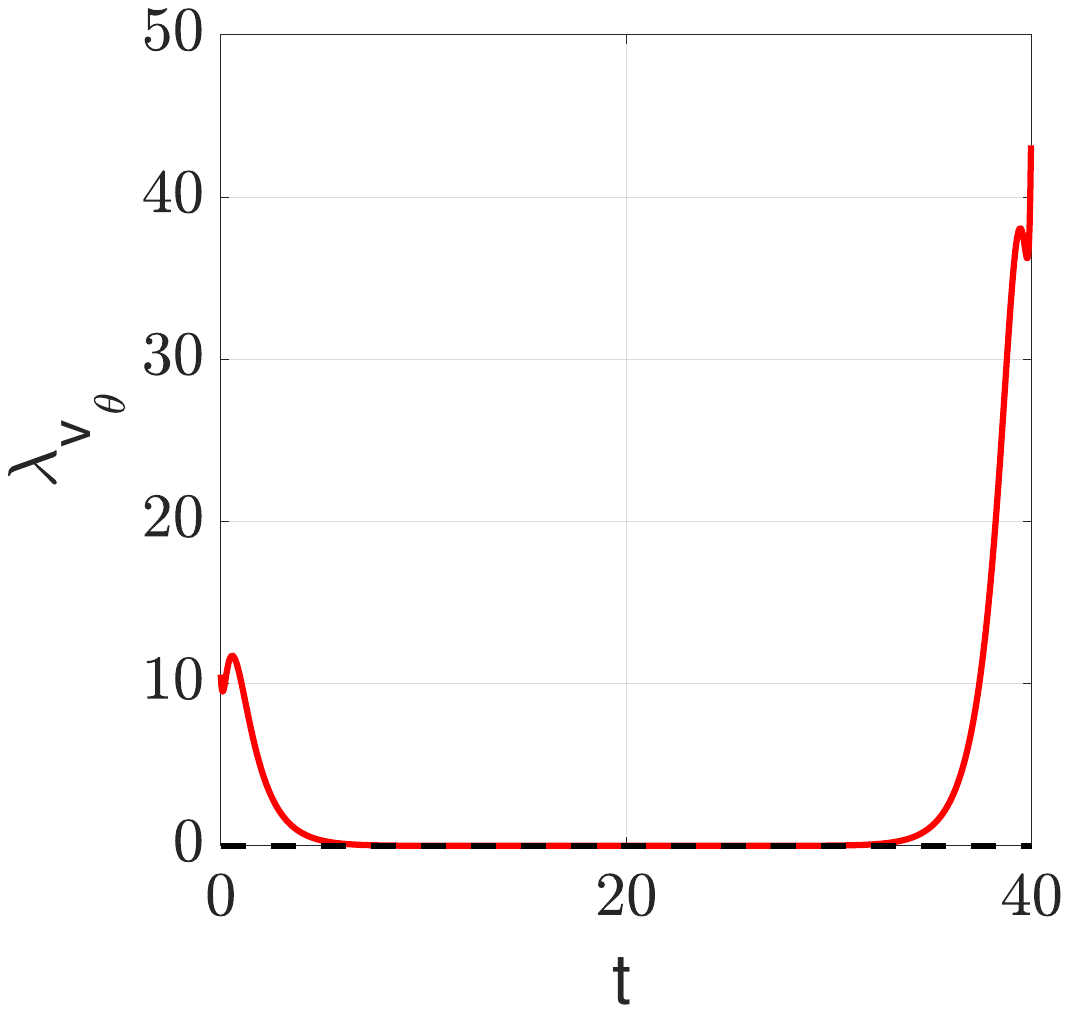}
\vspace{-0cm}
\includegraphics[width=\size\textwidth]{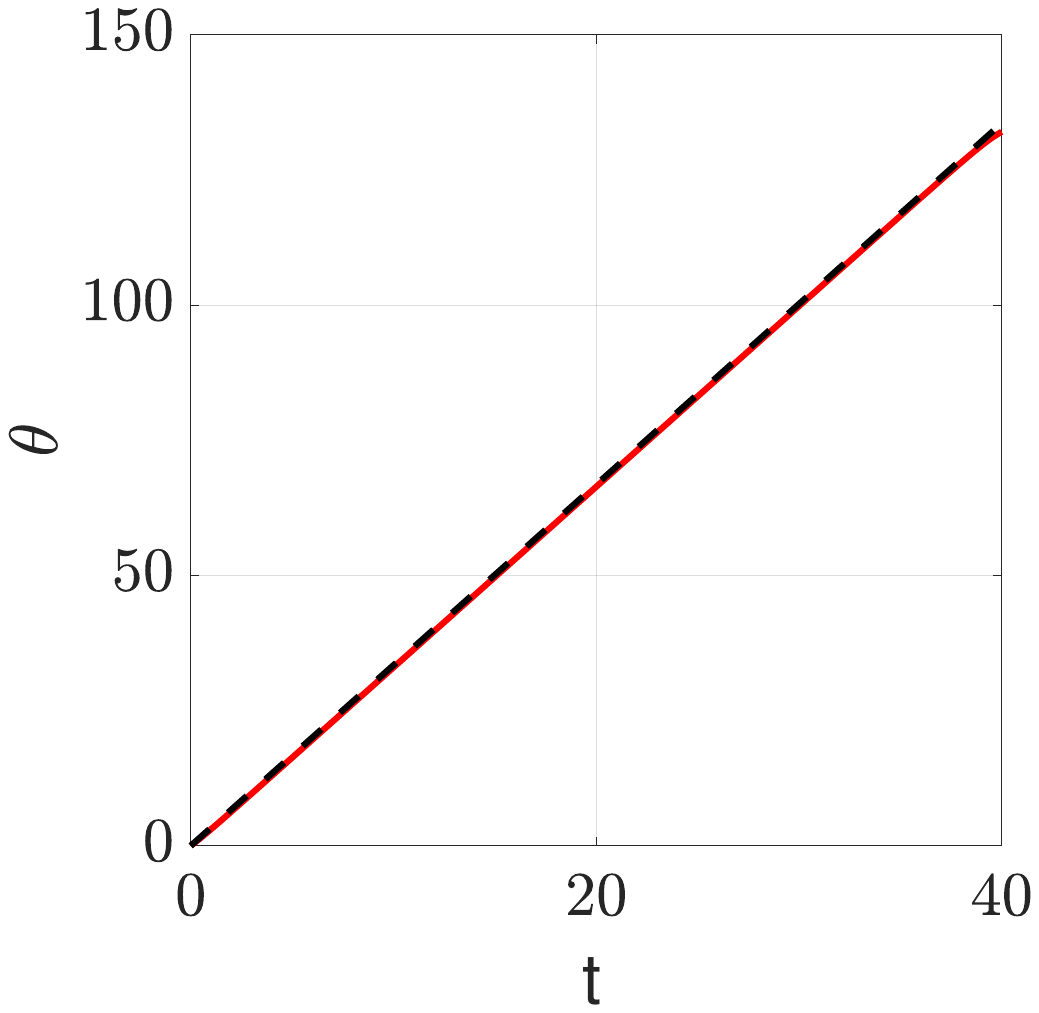}
\hspace{\sizeh cm}
\includegraphics[width=\size\textwidth]{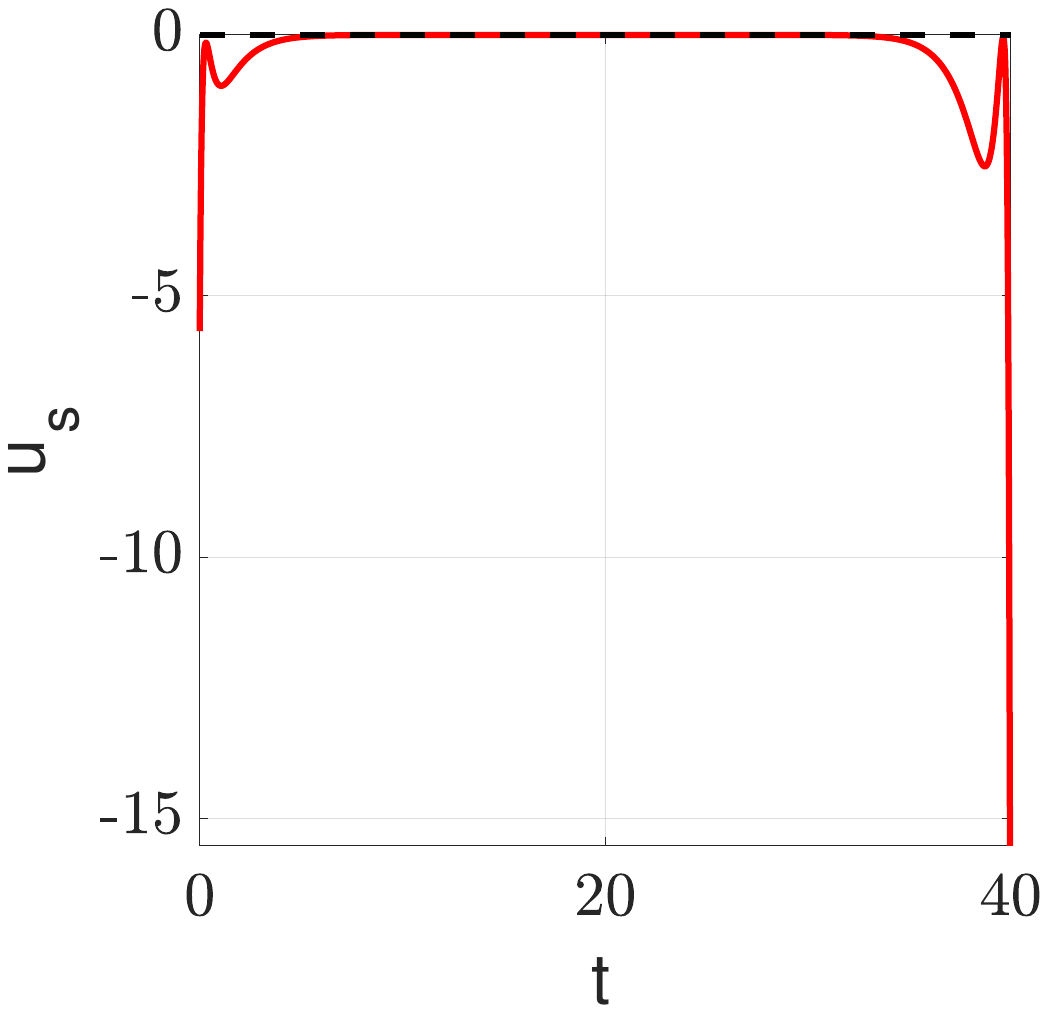}
\hspace{\sizeh cm}
\includegraphics[width=\size\textwidth]{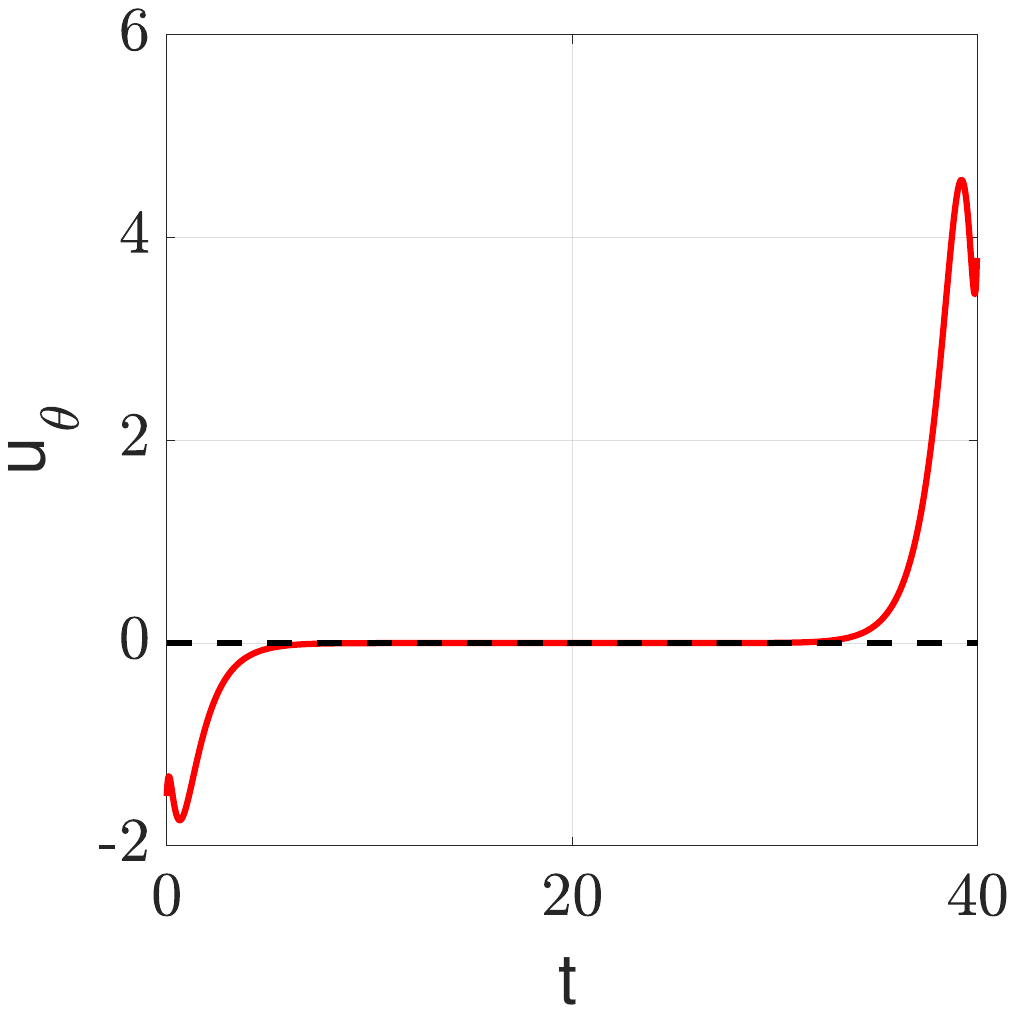}
\vspace{-0.7cm}
\caption{The turnpike property is depicted for the control $u$ and state variables $s, v_s, v_\theta$. The trim depicted in dashed black and the optimal solution in red. Here we observe a partial turnpike because the trim is not constant in $\theta$. Convergence of $\theta$ to trim holds on the whole interval $[0, t_f]$, since the turnpike is linear in this case.}
\label{fig:Kepler-example-2}
\end{figure}

\begin{remark}
Hamiltonian system from PMP associated with the reduced problem is
\begin{equation} \label{eq:PMP.Kepler}
\begin{aligned}
\dot s  &= v_s,   &  \hspace{-1cm} \dot p_s = -v_\theta^2 p_{v_s} - \frac{2k}{s^3}p_{v_s} + \frac{2}{s^2}v_\theta v_s p_{v_\theta} - \frac{u_\theta}{m_2 s^2} p_{v_\theta} + (s - \bar s), \\
\dot v_s &= sv_\theta^2 - \frac{k}{s^2} + \frac{1}{m_2^2} p_{v_s}, & \dot p_{v_s} = -p_s - \frac{2}{s}v_\theta p_{v_\theta} + (v_s - \bar v_s), \\
\dot v_\theta &= -\frac{2}{s}v_\theta v_s + \frac{1}{(m_2 s)^2}p_{v_\theta},  &
\dot p_{v_\theta} = -2 s v_\theta p_{v_s} - \frac{2}{s}v_\theta v_s p_{v_\theta} + (v_\theta - \bar v_\theta).
\end{aligned}
\end{equation}
Critical point of \eqref{eq:PMP.Kepler} is $(\bar s, \bar v_s, \bar v_\theta, \bar p_s, \bar p_{v_s}, \bar p_{v_\theta}) = (\bar s, 0, \sqrt{\frac{k}{\bar s^3}}, 0,0,0)$, assuming $\bar u = (0,0)$. Matrix obtained by linearization of \eqref{eq:PMP.Kepler} around this point has no eigenvalues with zero real part, thus, the hyperbolicity property holds in this case. The exponential trim turnpike property follows from Theorem~\ref{th:trim-turnpike}. This result is illustrated by numerical simulations shown in  Figures~\ref{fig:Kepler-example-1}-\ref{fig:Kepler-example-2}.
\end{remark}

For the implementation, we just consider the same boundary conditions and reference values as in reference \cite{FFOSW22}, i.e.\ $(s_0,\theta_0) = (5,0)$, $s_T = 6.0$ and $\bar s = 4.5$. The corresponding tracked trim reads
\begin{equation} \label{eq:trim.Kepler}
    (\bar s, \bar \theta, \bar v_s, \bar v_\theta) = \left(4.5, \, \bar v_\theta t, \, 0.0,\, \sqrt{\frac{k}{4.5^3}} \right).
\end{equation}
We solve the reduced problem using single shooting direct method in CasADI solver~\cite{Casadi:2019}. Control system is discretized by Runge-Kutta 4 method. The final time is fixed to $T = 40$. The obtained solution is represented in Figures~\ref{fig:Kepler-example-1} and \ref{fig:Kepler-example-2}. The obtained numerical results show that the optimal solution converges to a circular orbit with constant angular and zero radial velocities. The orbit corresponds to the trim \eqref{eq:trim.Kepler}.  

\subsection{Rigid body problem} \label{sec:rb}

We consider rotations of a rigid body in 3-dimensional space. Configuration manifold is given by $Q = SO(3)$.
We assume that the movement of the rigid body is controlled by an external force. The motions can be described by the Lagrange-d'Alembert principle~\cite{RalphMarsden2021}, with the Lagrangian given by
\begin{equation} 
\label{eq:rigidbody-Lagrangian}
L(R, \dot{R}) = \frac{1}{2} \langle R^{-1} \dot{R}, \, I R^{-1} \dot{R}\rangle,
\end{equation}
where $R \in SO(3)$ is an orientation of the rigid body, $\langle \cdot, \cdot \rangle$ stands for the inner product of matrices $I$ represents the inertia matrix. In what follows, we assume that $\langle A, B \rangle = \mathrm{tr}(A^*B)$ for two given matrices $A,B\in SO(3)$, and also that $I$ is diagonal.
The Lagrangian \eqref{eq:rigidbody-Lagrangian} is invariant with respect to the Lie group action of $G = SO(3)$ on itself by left multiplication, which is a non-abelian Lie group. Moreover, we consider an external force $F(R,\dot{R},u)$ equivariant with respect to the action of $G$. 
Lagrange-d'Alembert principle 

leads to the corresponding controlled Euler-Lagrange equations 
\begin{equation} 
\label{eq:EL-Rigidbody}
\frac{d}{dt}\frac{\partial}{\partial \dot R}L - \frac{\partial}{\partial R}L = F(R(t), \dot R(t), u(t)), 
\end{equation}
$u(\cdot)$ being the control function.
Invariance of the Lagrangian and equivariance of the force field induce the definition of the reduced Lagrangian on $TQ/G \equiv \mathfrak{g} $ with $\Omega \in \mathfrak{g}$, given by
\begin{equation} 
\label{eq:Reduced-lagragian-rigidbody}
l(\Omega) = \frac{1}{2} \langle \Omega, I \Omega\rangle,
\end{equation}
and reduced Euler-Lagrange equations 
\begin{equation} \label{eq:Reduced-EL-RB}
\frac{d}{dt}\frac{\partial}{\partial \Omega}l - \mathrm{ad}^*_{\Omega}\frac{\partial}{\partial \Omega}l = f(\Omega, u). \\
\end{equation}
Using identification $\mathfrak{g} = so(3) \equiv \R^3$ and plugging \eqref{eq:Reduced-lagragian-rigidbody} in \eqref{eq:Reduced-EL-RB}, we obtain reduced Euler-Lagrange equations on $\R^3$ in the following form
\begin{equation} 
\label{eq:EL-rigidbody-reduced}
I \dot{\Omega} = I \Omega \times \Omega + f(\Omega, u),  
\end{equation}
where $\Omega \times \Omega$ is a cross product in $\R^3$. The reconstruction equation for $R \in SO(3)$ is given by $\dot R = R \Omega$.
\begin{remark}
To  simplify equations, we consider the reduced force $f(\Omega, u) = u\diff\Omega$, which corresponds to $F(R(t), \dot R(t), u(t)) = Ru \diff R$.  
\end{remark}
Next, we consider the following optimal control problem
\begin{equation} 
\label{eq:rigidbody-pb}
\begin{aligned}
& \min  && J(R, \dot R, u) = \int_0^T F^0(R, \dot R,u) \diff t \\
& \text{s.t.} && 
\begin{cases}
\frac{d}{dt}\frac{\partial}{\partial \dot R}L - \frac{\partial}{\partial R}L = Ru, \\
R^{-1} \dot R(0) = \Omega_0, ~~~ R^{-1} \dot R(T) = \Omega_T.
\end{cases}
\end{aligned}
\end{equation}
Let us set the cost function
$$F^0(R, \dot R, u) = \frac{1}{2}\left(\| R^{-1} \dot R - \Omega_\mathrm{ref}\|_2^2 + \|u - u_\mathrm{ref} \|_2^2 \right),$$ 
where $u_{ref} \in \R^3$, $\Omega_\mathrm{ref} \in \R^3$ are constant. 
It is easy to see that $F^0$ is invariant with respect to $SO(3)$. The induced reduced function $f^0$ on $\R^3 \times \R^3$ takes the form 
%
$$f^0(\Omega, u) = \frac{1}{2} \left(\|\Omega - \Omega_\mathrm{ref}\|_2^2  + \| u - u_{ref}\|_2^2 \right).$$
The resulting reduced OCP is equivalent to \eqref{eq:rigidbody-pb} by Theorem~\ref{th:reduced.OCP.equiv}. It is given by
\begin{equation} 
\label{eq:reduced-rigidbody-pb}
\begin{aligned}
& \text{min}  && J(R, \dot R, u) = \frac{1}{2}\int_0^T \left(\|\Omega - \Omega_\mathrm{ref}\|_2^2  + \| u - u_{ref}\|_2^2 \right) \diff t \\
& \text{s.t.} && 
\begin{cases}
I \dot{\Omega} &= I \Omega \times \Omega + u  \\
\Omega(0) &= \Omega_0, ~~~ \Omega(T) = \Omega_T
\end{cases}
\end{aligned}
\end{equation}
%
%
Static problem associated with \eqref{eq:reduced-rigidbody-pb} is defined by
\begin{equation} 
\label{eq:rigidBody-trim}
\begin{aligned}
& \min  & f^0(\Omega, u) &= \frac{1}{2} \left(\|\Omega - \Omega_\mathrm{ref}\|_2^2  + \| u - u_{ref}\|_2^2\right), \\
&  \text{s.t.}  & 0 &= I \Omega \times \Omega + u.   \\
\end{aligned}
\end{equation}
If $\Omega_{ref}$ is equal to one of the canonical basis vectors of $\R^3$, then \eqref{eq:rigidBody-trim} admits a unique solution $(\bar{\Omega}, \bar u) = (\Omega_\mathrm{ref}, 0)$. The corresponding trim, $\bar R$, is a solution of $\dot R = R \bar{\Omega}$, and the initial condition $\bar R(0)$ of the trim is determined by Theorem~\ref{th:trim-turnpike}. As in the Kepler example above, we can easily verify the hyperbolic property of the reduced system.

\begin{figure}[ht]  
\centering
\def\size{0.33}
\def\Size{0.21}
\def\sizeh{-0.7}
\includegraphics[width=\size\textwidth]{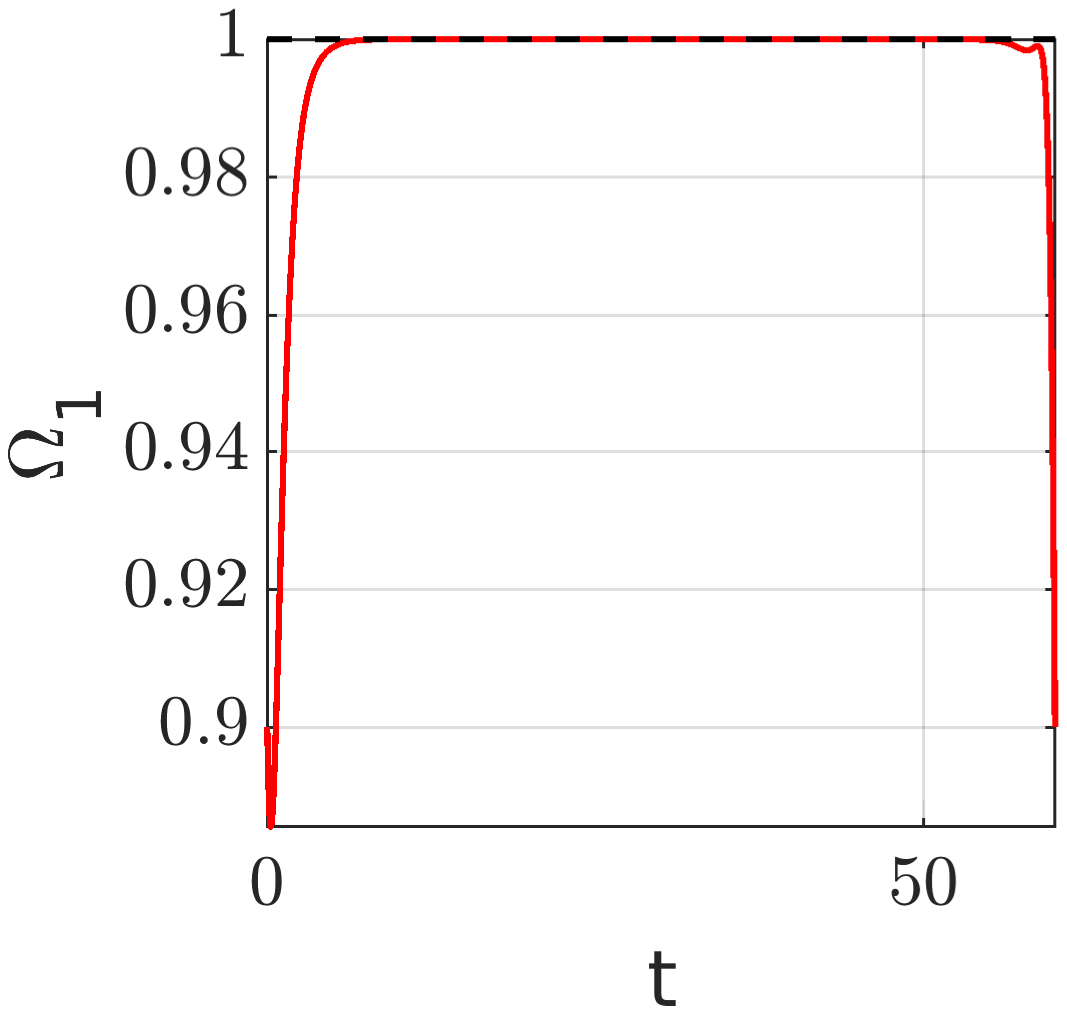}
\hspace{\sizeh cm}
\includegraphics[width=\size\textwidth]{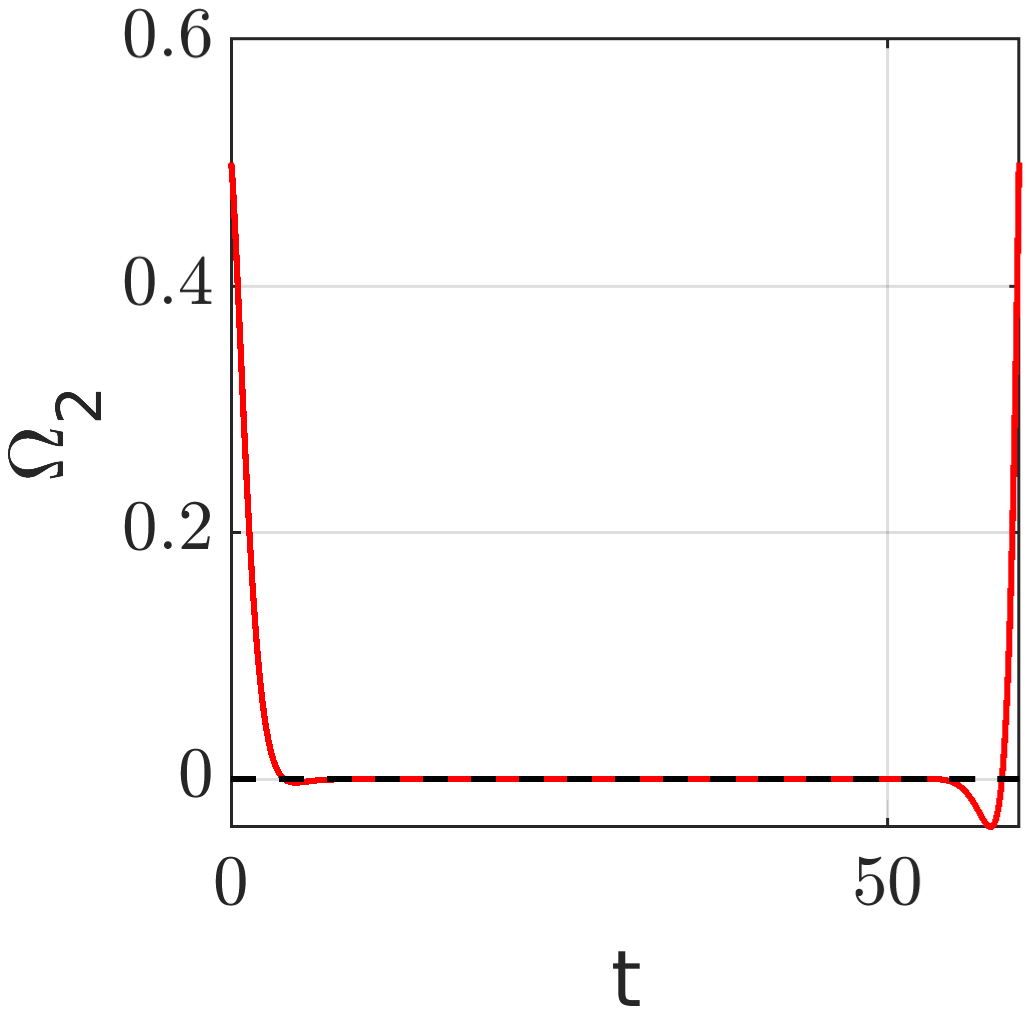}
\hspace{\sizeh cm}
\vspace{-2cm}
\includegraphics[width=\size\textwidth]{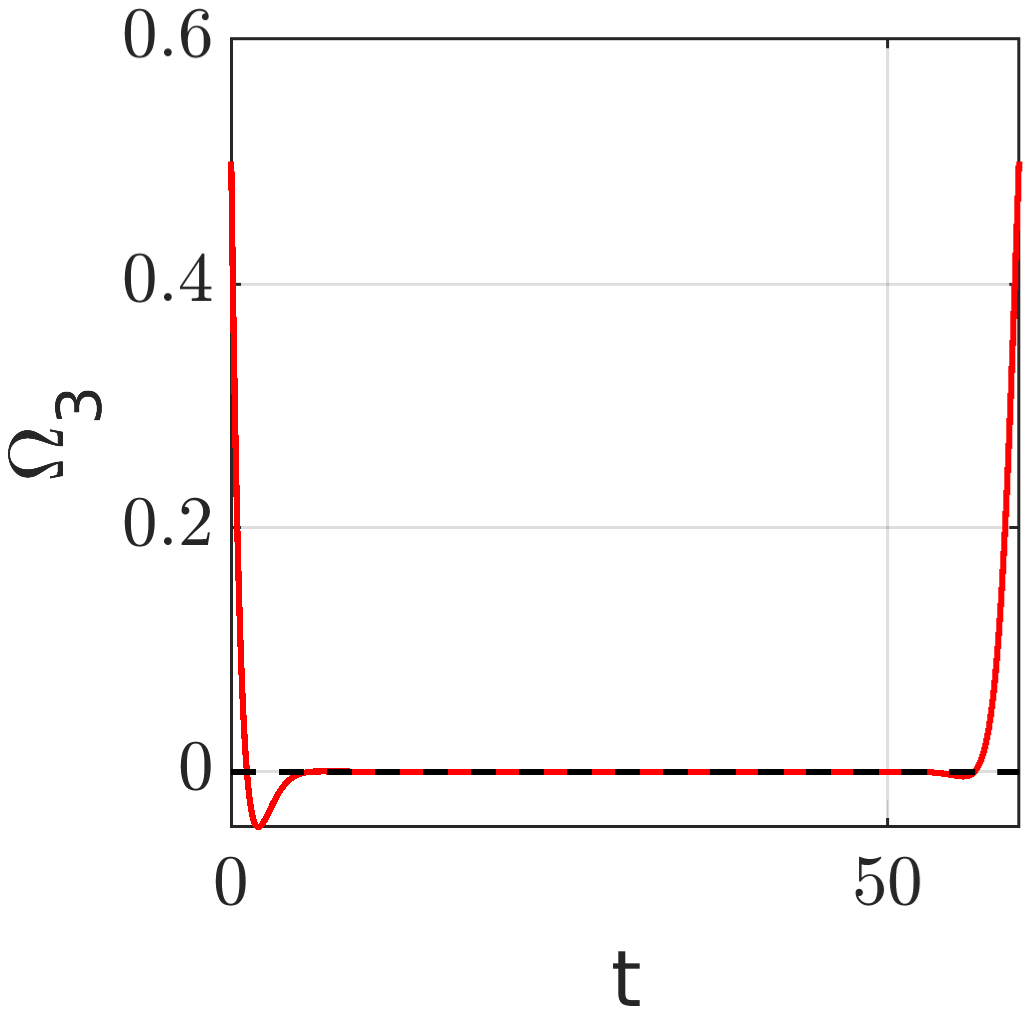}
\vspace{-1cm}
\includegraphics[width=\size\textwidth]{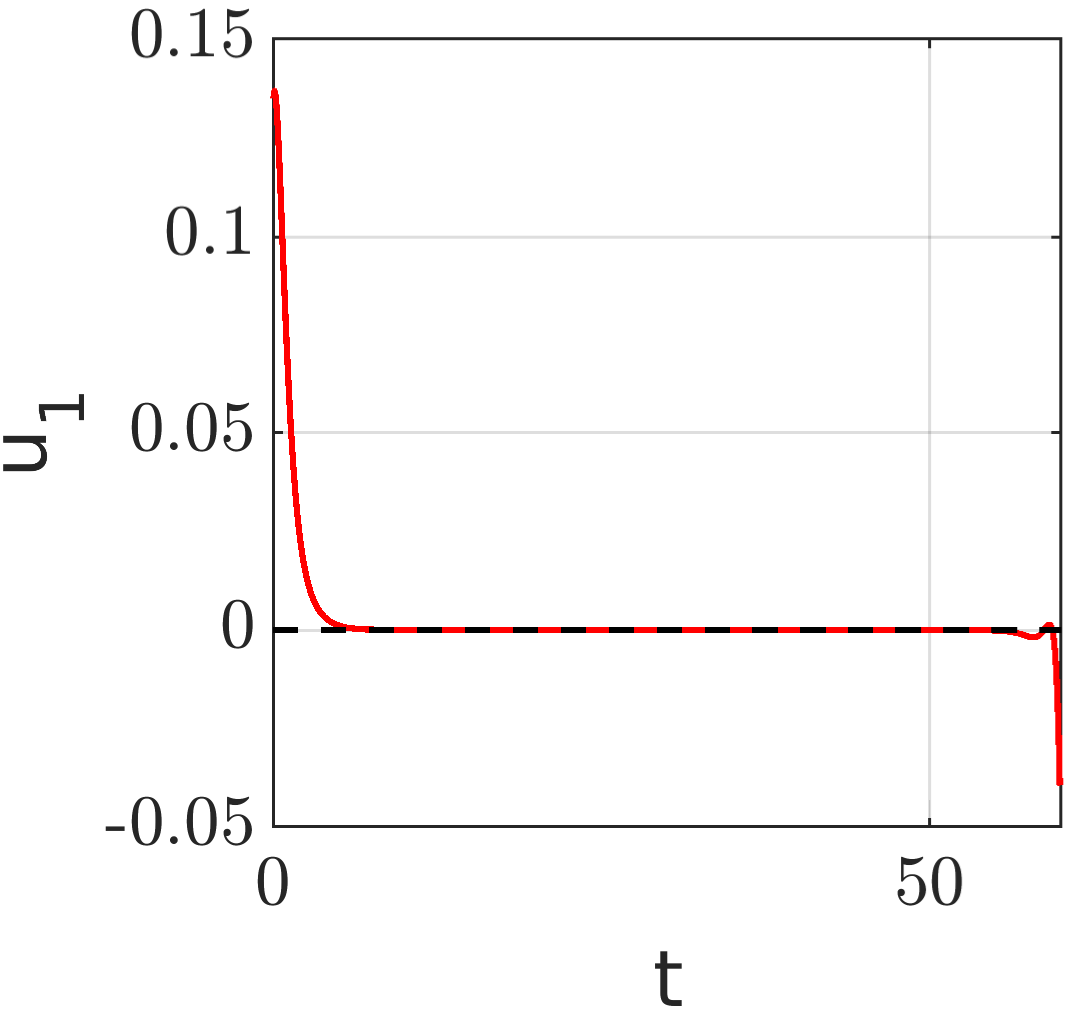}
\hspace{\sizeh cm}
\includegraphics[width=\size\textwidth]{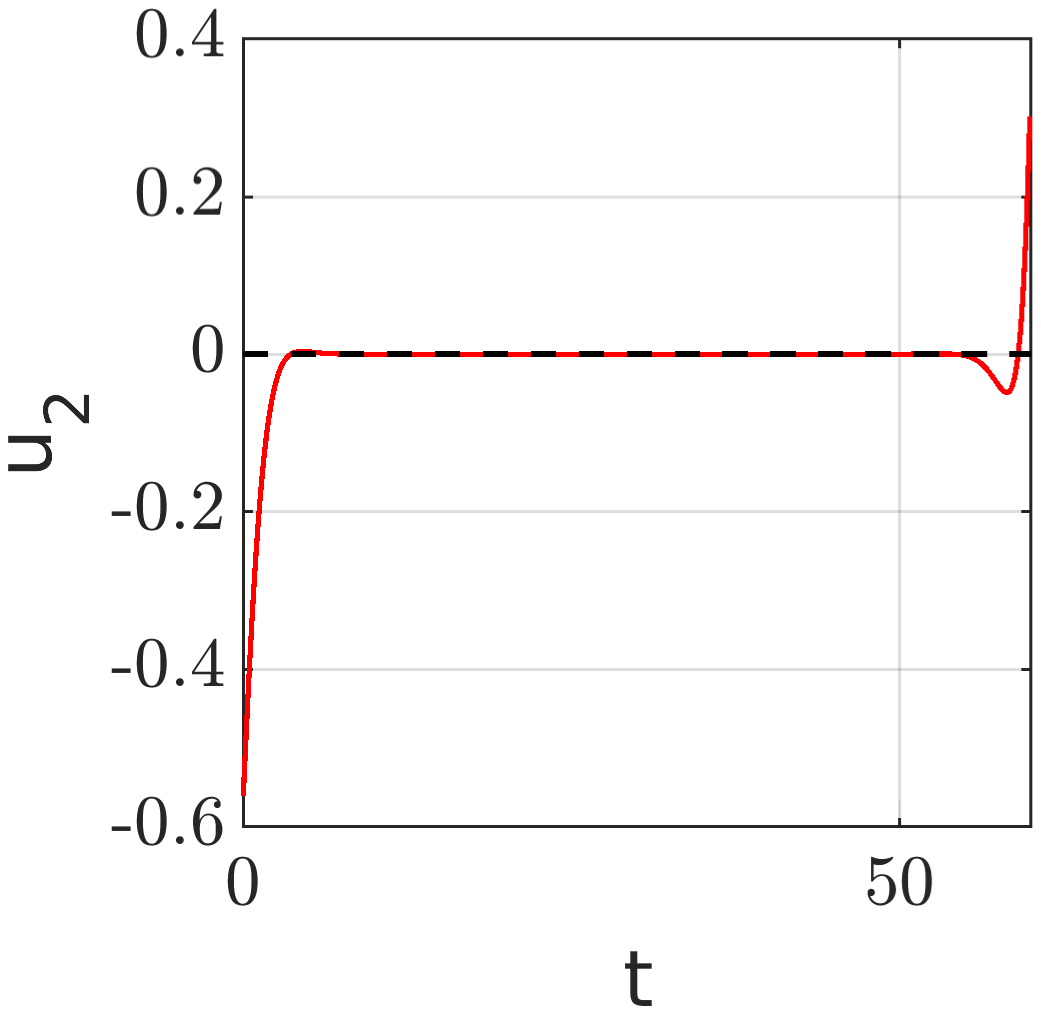}
\hspace{\sizeh cm}
\vspace{-1cm}
\includegraphics[width=\size\textwidth]{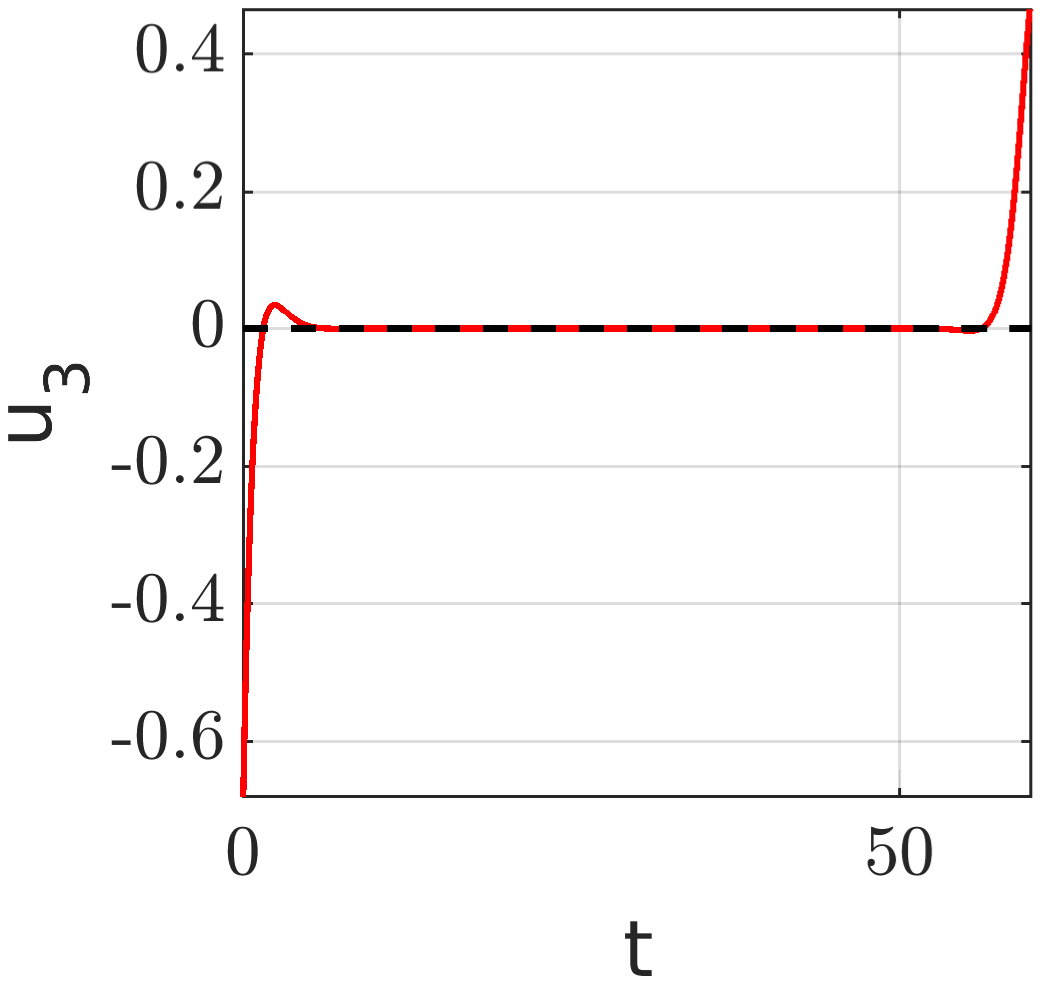}
\vspace{-0cm}
\includegraphics[width=\size\textwidth]{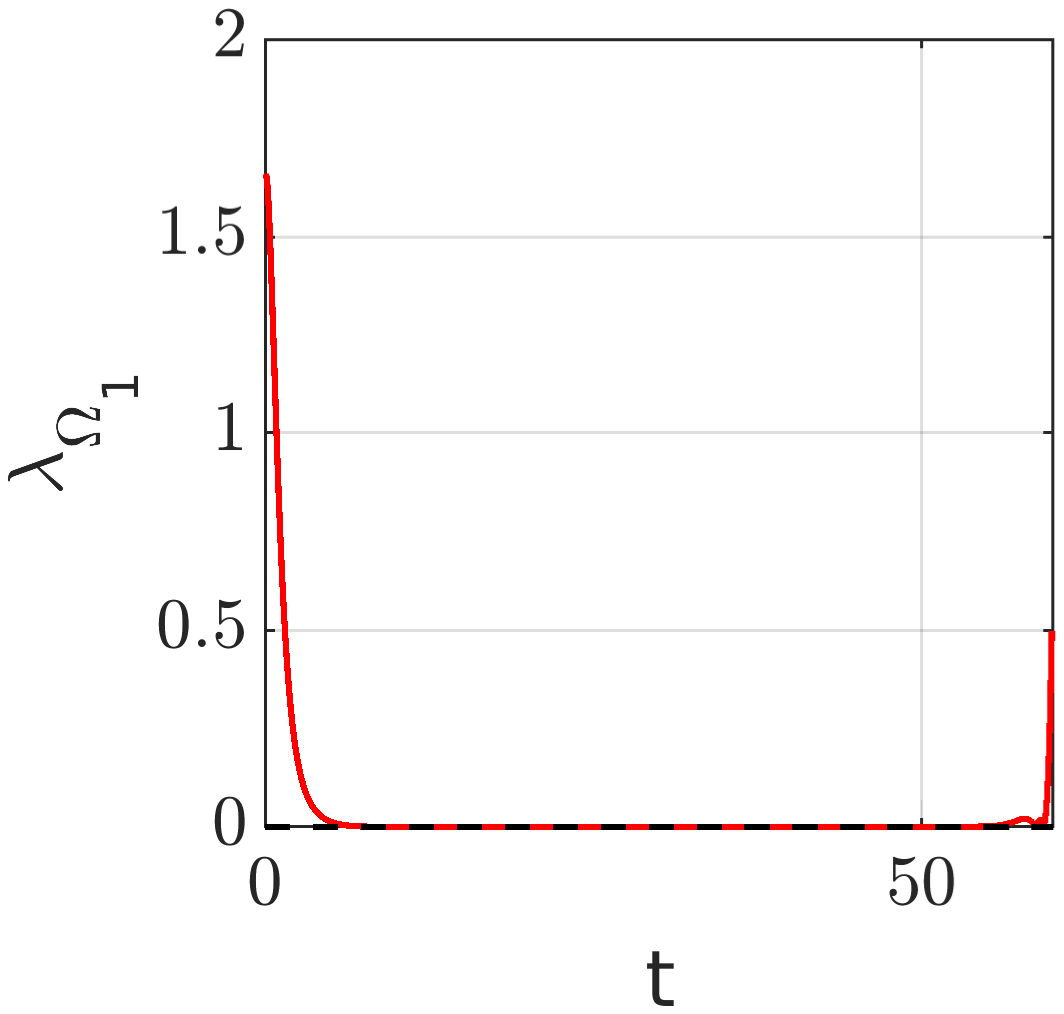}
\hspace{\sizeh cm}
\includegraphics[width=\size\textwidth]{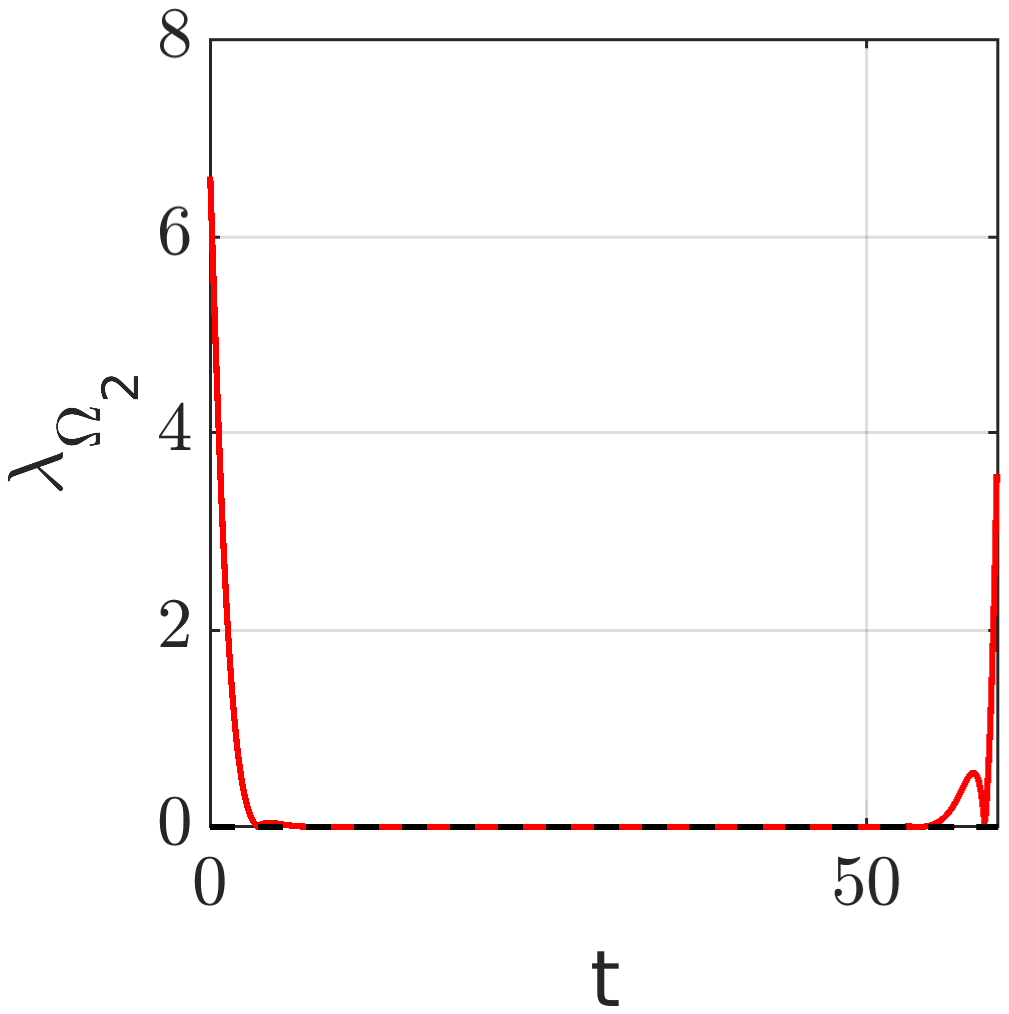}
\hspace{\sizeh cm}
\includegraphics[width=\size\textwidth]{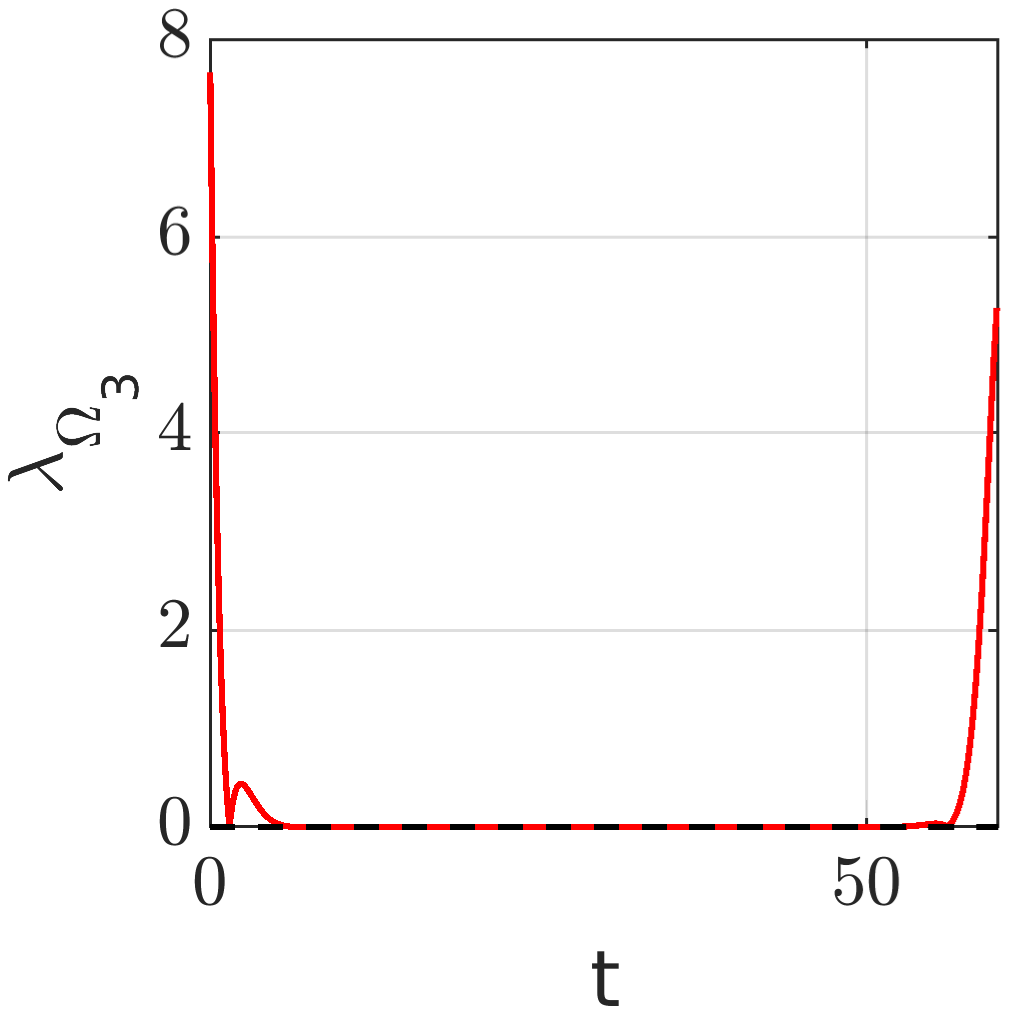}
\vspace{-0.7cm}
\caption{The turnpike property is depicted for the control $u$ and state variables $\Omega_1, \Omega_2,\Omega_3$ and their adjoints $\lambda_1, \lambda_2,\lambda_3$. The trim depicted in dashed black and the optimal solution in red. 
For the representation we consider $\Omega_0 = \Omega_T = (0.9,0.5,0.5)^T$, $T = 60$, $I = \mathrm{diag}(1,5,10)$.}
\label{fig:rigidbody-target}
\end{figure}

\begin{figure}[ht]  
\centering
\includegraphics[width=0.54\textwidth, height=0.5\textwidth]{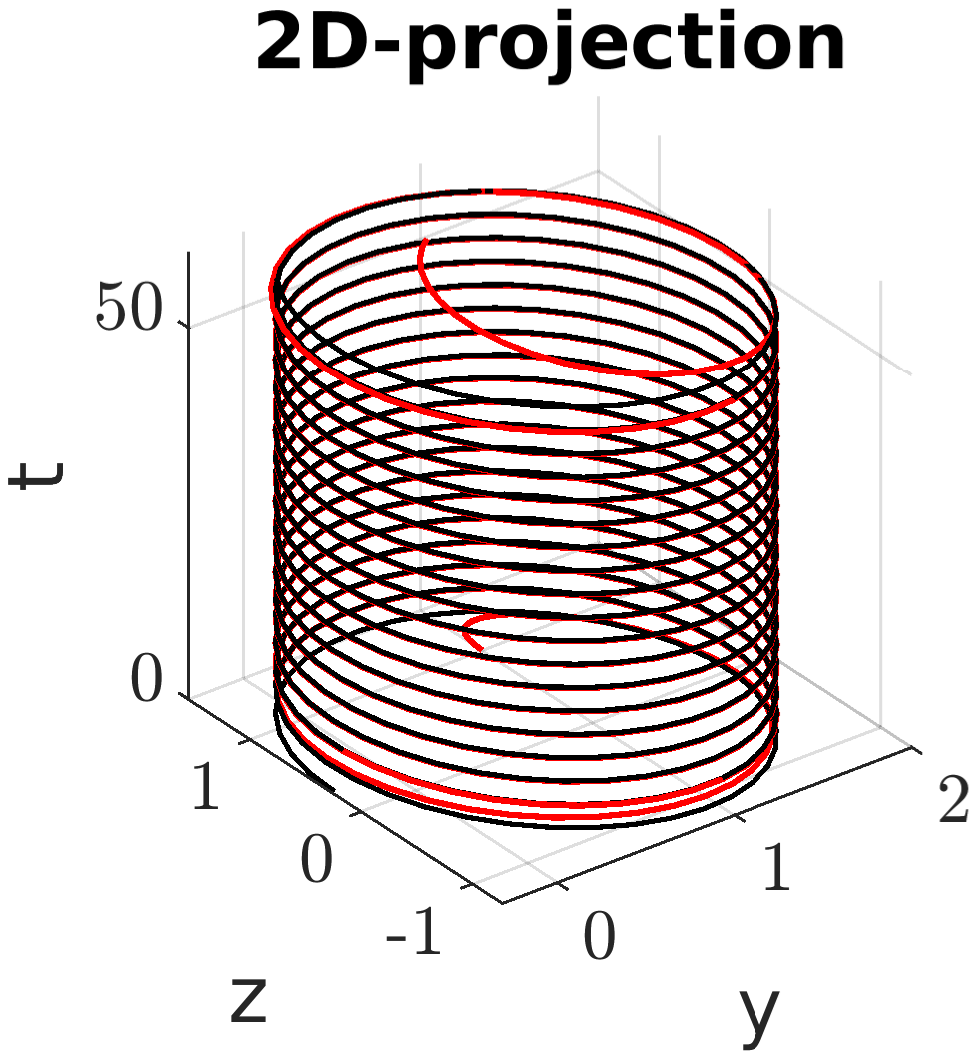}
\hspace{-1.5cm}
\includegraphics[width=0.54\textwidth, height=0.5\textwidth]{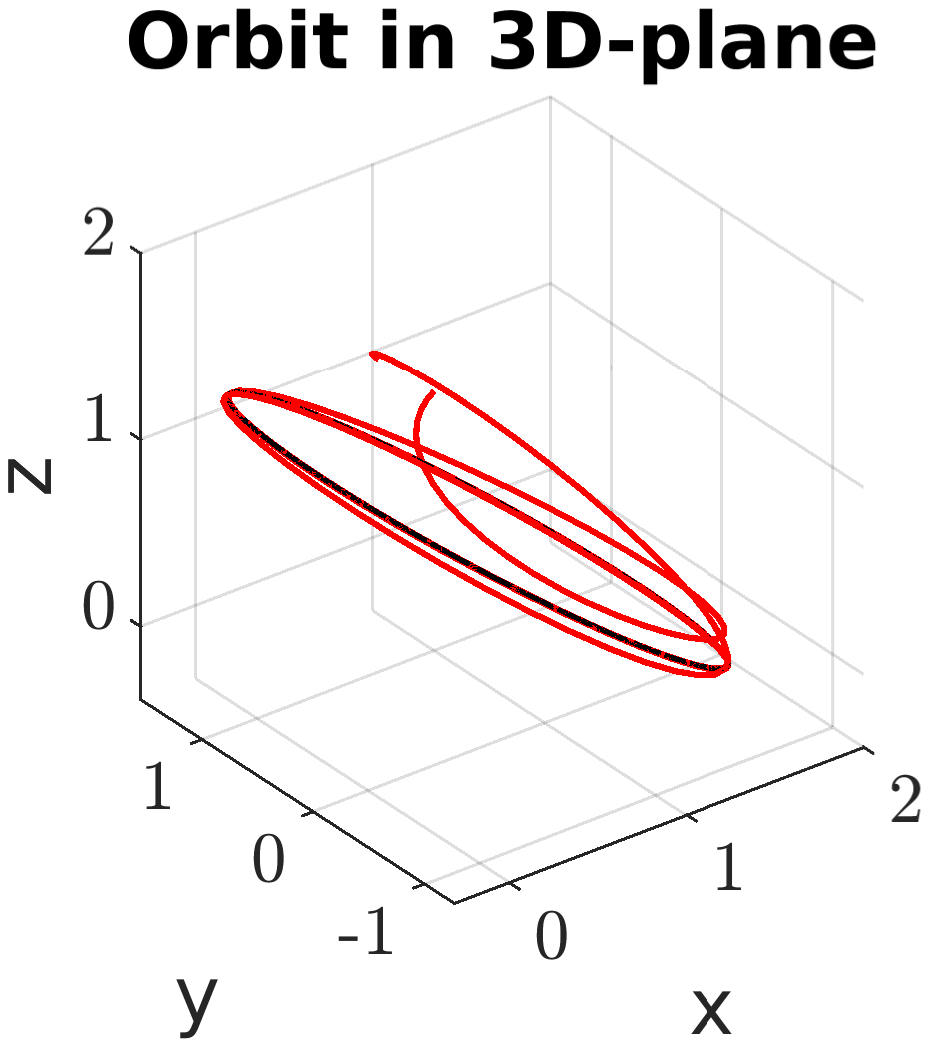}
\vspace{-1.cm}
\caption{Plots of optimal trajectory in red and trim in black. The plot on the left shows the convergence of $(y,z)$ coordinates of the optimal trajectory to $(y,z)$ coordinates of the trim in time, while the right plot shows the convergence in phase space, where we can observe the convergence of geometrical curves.} 
\label{fig:rigidbody-trim}
\end{figure}

To obtain optimal solutions, we use CasADi software and Runge-Kutta 4 discretization of the control system as before. 
To perform the integration of the rotational matrix $R$, we use the quaternion parameterization of $R$ and normalize the quaternion vector after each integration step. 
The numerical optimal solution is depicted together with the trim determined by \eqref{eq:rigidBody-trim} in Figures~\ref{fig:rigidbody-target} and \ref{fig:rigidbody-trim}. As suggested by Theorem~\ref{th:trim-turnpike}, the initial condition for the trim is determined from the optimal solution, obtained numerically in this case. For this we find the value of $R(t)$ at $\lfloor T/2 \rfloor$ and integrate $\dot R = R \bar{\Omega}$ backward in time for $t\in \left[\lfloor T/2 \rfloor, 0 \right]$ to determine the initial condition of the trim at $t = 0$. We observe that the optimal solution converges to that specific trim, which confirms our \emph{trim turnpike property}. As can be seen, the trim corresponds to a rotational motion of the body around one of the principal axes with constant rotational velocity.

\begin{remark}
For the visualization of the rotation in Figure~\ref{fig:rigidbody-trim}, we first compute the solution in quaternions parameterization, then we reconstruct the rotational matrices $R$ and $\bar R$. Finally we visualize the evolution of the vector $(1,1,1)^T$ i.e. by computing $R^T \cdot (1,1,1)^T$ (resp. $\bar R^T \cdot(1,1,1)^T$) at each time step. Notice that the use of quaternions parameterization leads to a change of basis.
\end{remark}


\subsection{Rigid body with rotors problem} \label{sec:rb.with.rotors}
%

\begin{figure}[ht!]
\centering  
\includegraphics[width=0.4\textwidth]{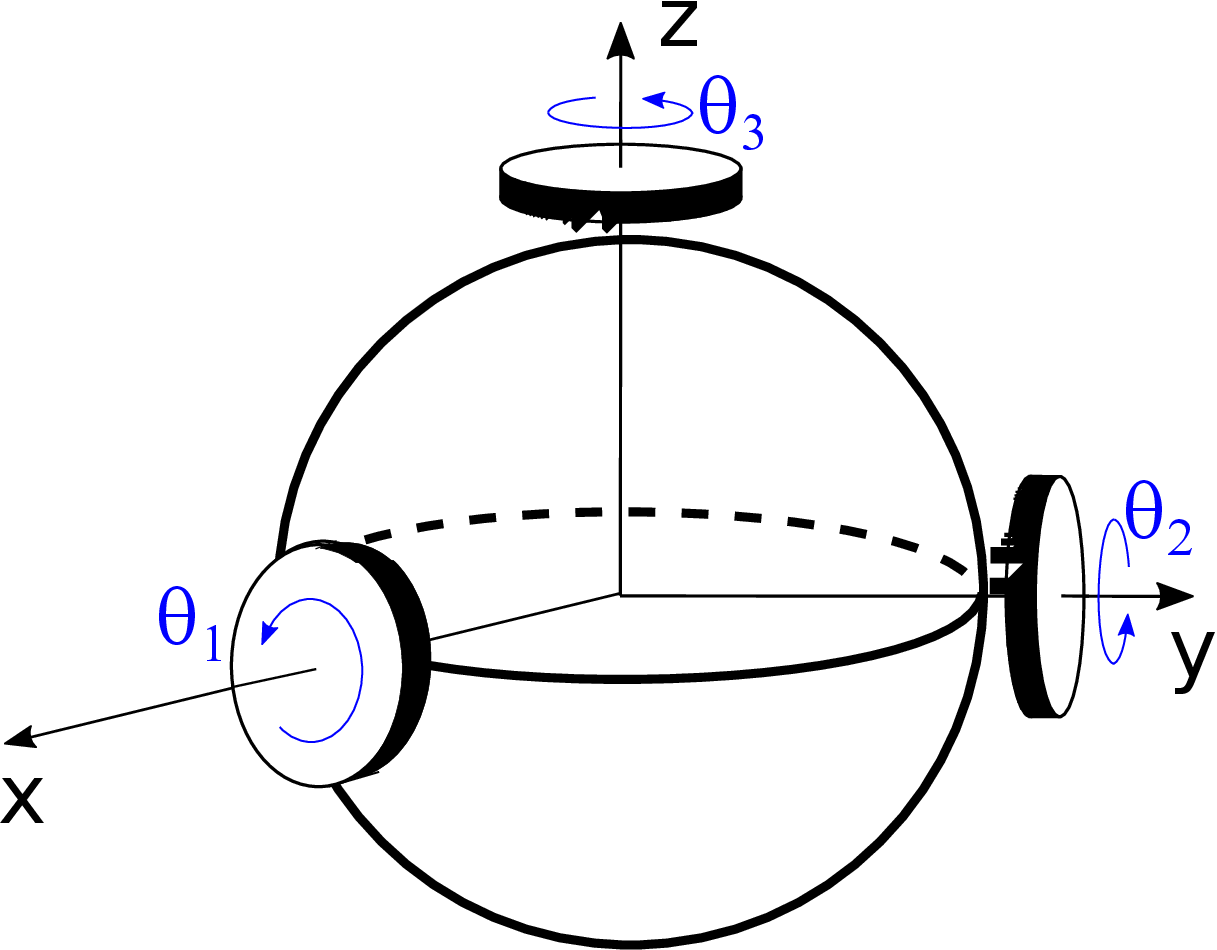}
\vspace{0.3cm}
\caption{Visualization of a rigid body with three rotors, where the body is represented as a ball for simplicity. Rotors spin under a torque force due to an external action. Movement of rotors influence the rotation of the rigid body.}
\label{fig:rigidbody}
\end{figure}

In the last example, we consider rotations of a rigid body with three attached controlled rotors, see Figure~\ref{fig:rigidbody}.
The configuration space is now given by $Q = SO(3)\times S^1 \times S^1 \times S^1$ with the configuration of the rigid body $R \in SO(3)$ and configuration of the $i$th rotor $\theta_i \in S^1$ for $i = 1,2,3$. We set $\theta = (\theta_1, \theta_2, \theta_3)$ for the configuration of all three rotors. 
The equations of motions of the rigid body with rotors can be derived from the Lagrange--d'Alembert principle with the Lagrangian given by
\begin{equation} \label{eq:rigidbody.Lagrangian}
L(R, \dot{R}, \theta, \dot{\theta}) = \frac{1}{2} \langle R^{-1} \dot{R}, I R^{-1} \dot{R}\rangle + \frac{1}{2} \langle R^{-1} \dot{R} + \dot{\theta}, K \left( R^{-1} \dot{R} + \dot{\theta} \right) \rangle,
\end{equation}
where $I,K$ are inertia matrices, $I$ of the rigid body and $K$ of rotors, which we assume to be diagonal. As in Section~\ref{sec:rb}, we consider a group action $G = SO(3)$ on itself by left multiplication. In addition, the external forces are assumed to act only on rotors and the force field has the form $F(R, \dot R, \theta, \dot \theta, u) = u \, d\theta$.
In this case, $F$ is both invariant and orthogonal with respect to the group action by $G$. 

The reduced Lagrangian corresponding to \eqref{eq:rigidbody.Lagrangian} is defined on $TQ/G \approx so(3) \times T(S^1 \times S^1 \times S^1)$ and given by
\begin{equation}
\label{eq:reduced.lagragian.rigidbody}
l(\Omega, \theta, \dot{\theta}) = \frac{1}{2} \langle \Omega, I \Omega\rangle + \frac{1}{2} \langle \Omega + \dot{\theta}, K (\Omega + \dot{\theta} ) \rangle,
\end{equation}
where $\Omega \in SO(3)$ and $(\theta, \dot \theta) \in T(S^1 \times S^1 \times S^1)$. According to the Lagrange--d'Alembert principle, the forced Euler-Lagrange equations on $SO(3)\times T (S^1 \times S^1 \times S^1) \approx \R^3\times \R^3\times (S^1 \times S^1 \times S^1)$ 
are
\begin{equation} 
\label{eq:EL.rigidbody.reduced}
\begin{aligned}     
\dot{\Omega} &= I^{-1}\left( \Pi \times \Omega \right) - I^{-1} u,  \\
\dot{\theta} &= v_\theta, \\
\dot{v}_{\theta} &= -I^{-1}\left( \Pi \times \Omega \right) + \left( K^{-1} + I^{-1} \right)u,
\end{aligned}    
\end{equation}
where $\Pi = \left(I+K\right) \Omega + K v_\theta.$
Since the right-hand side of equations in \eqref{eq:EL.rigidbody.reduced} is independent from $\theta$, it is a cyclic variable of the problem, so that \eqref{eq:EL.rigidbody.reduced} is equivariant with respect to translations in $\theta$. Therefore, we can perform a further reduction of \eqref{eq:EL.rigidbody.reduced}, described in Remark~\ref{rem:cyclic.OCP}. The control system reduced with respect to both $SO(3)$ and $S^1 \times S^1 \times S^1$ is defined on $\R^3 \times \R^3$ and takes the following form
\begin{equation} 
\label{eq:EL.rigidbody.reduced2}
\begin{aligned}     
& \dot{\Omega} = I^{-1}\left( \Pi \times \Omega \right) - I^{-1} u,  \\
& \dot{v}_{\theta} = -I^{-1}\left( \Pi \times \Omega \right) + \left( K^{-1} + I^{-1} \right)u,
\end{aligned}    
\end{equation}
with the reconstruction equations $\dot R = R \Omega$ and $\dot{\theta} = v_\theta$.
Let us consider an optimal control problem of the form \eqref{eq:OCP_original} with the cost function defined by $f^0(\Omega, v_{\theta},u) = \frac{1}{2} \left(   \|\Omega - \Omega_\mathrm{ref}\|_2^2  + \| u - u_{ref}\|_2^2\right)$ and 
subject to the control system constraint given by  controlled Euler-Lagrange equations defined by \eqref{eq:rigidbody.Lagrangian} and force field~$F$. The defined OCP is symmetric with respect to $SO(3) \times S^1 \times S^1 \times S^1$ by construction.  The corresponding reduced problem is defined as follows
\begin{equation}
\label{eq:partial-reduced-rigidbody-pb-2}
\begin{aligned}
& \text{min}  & J(\Omega, v_{\theta}, u) &= \int_0^T f^0(\Omega, v_{\theta},u) \diff t \\
& \text{s.t.} & \dot{\Omega} &= I^{-1}\left( \Pi \times \Omega \right) - I^{-1} u,  \\
&& \dot{v}_{\theta} &= -I^{-1}\left( \Pi \times \Omega \right) + \left( K^{-1} + I^{-1} \right)u,
\end{aligned}
\end{equation}
Trim turnpikes associated with \eqref{eq:partial-reduced-rigidbody-pb-2} is solution of 
\begin{equation} 
\label{eq:rigidBody-full-trim}
\begin{aligned}
& \min  && f^0(\Omega, v_{\theta},u), \\
&  \text{s.t.}  && 0 = I^{-1}\left( \Pi \times \Omega \right) - I^{-1} u,  \\
&	 	 	&& 0 = -I^{-1}\left( \Pi \times \Omega \right) + \left( K^{-1} + I^{-1} \right)u.
\end{aligned}
\end{equation}
In case where $\Omega_\mathrm{ref}$ is equal to one of the canonical basis vectors of $\R^3$, \eqref{eq:rigidBody-full-trim} admits a unique solution $(\bar{\Omega}, \bar{v}_{\theta}, \bar u) = (\Omega_\mathrm{ref},0,0)$. Then, the corresponding trim is a solution of $\dot R = R \bar{\Omega}$, and the initial condition $R(0)$ of the trim is determined by Theorem~\ref{th:trim-turnpike}. However, 
using any formal calculus solver as Maple or Mathematica, one can easily check that the linearization of the state-adjoint system associated with \eqref{eq:partial-reduced-rigidbody-pb-2} around a solution of \eqref{eq:rigidBody-full-trim} is not hyperbolic. More precisely, the matrix obtained by linearization has zero eigenvalues. Nonetheless, simulations show that the optimal control problem exhibits the trim turnpike and the convergence is observed toward the trim defined by Theorem~\ref{th:trim-turnpike}.  

\begin{figure}[ht]  
\centering
\includegraphics[width=0.54\textwidth, height=0.52\textwidth]{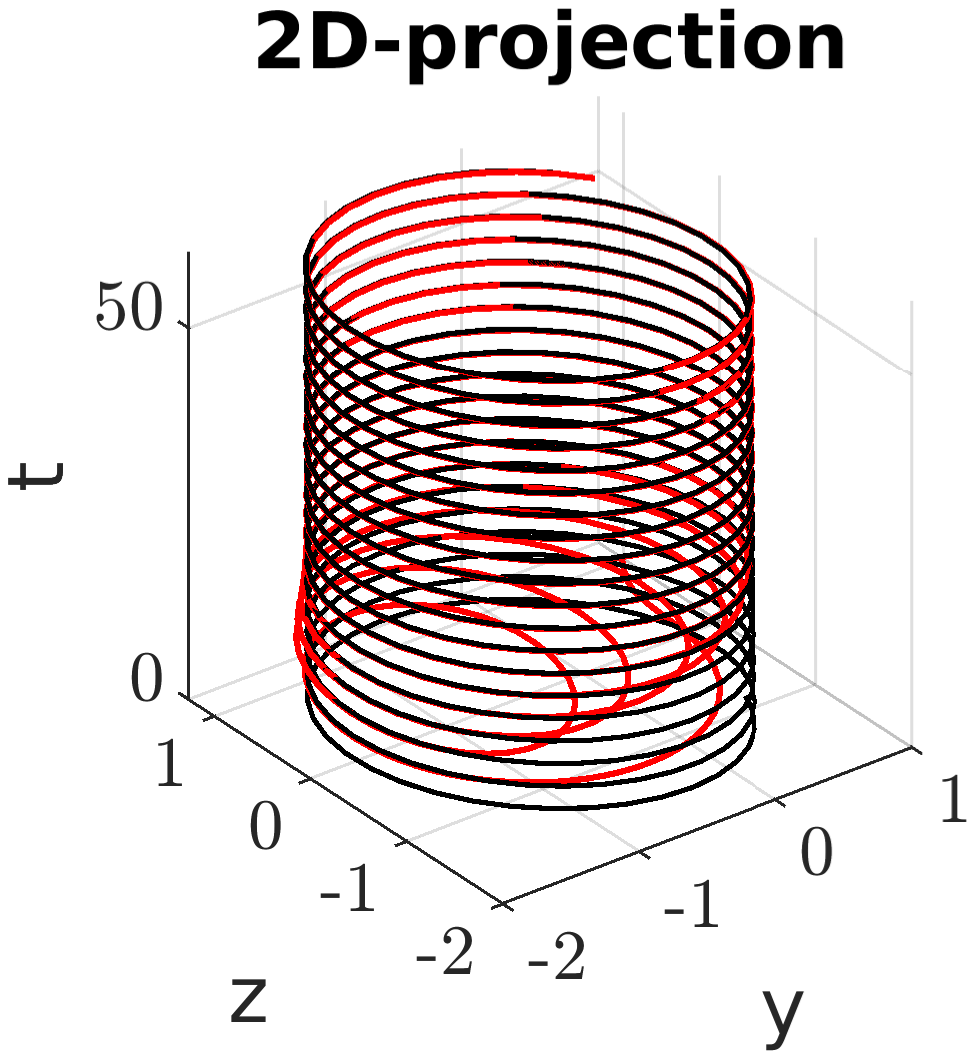}
\hspace{-1.5cm}
\includegraphics[width=0.54\textwidth, height=0.52\textwidth]{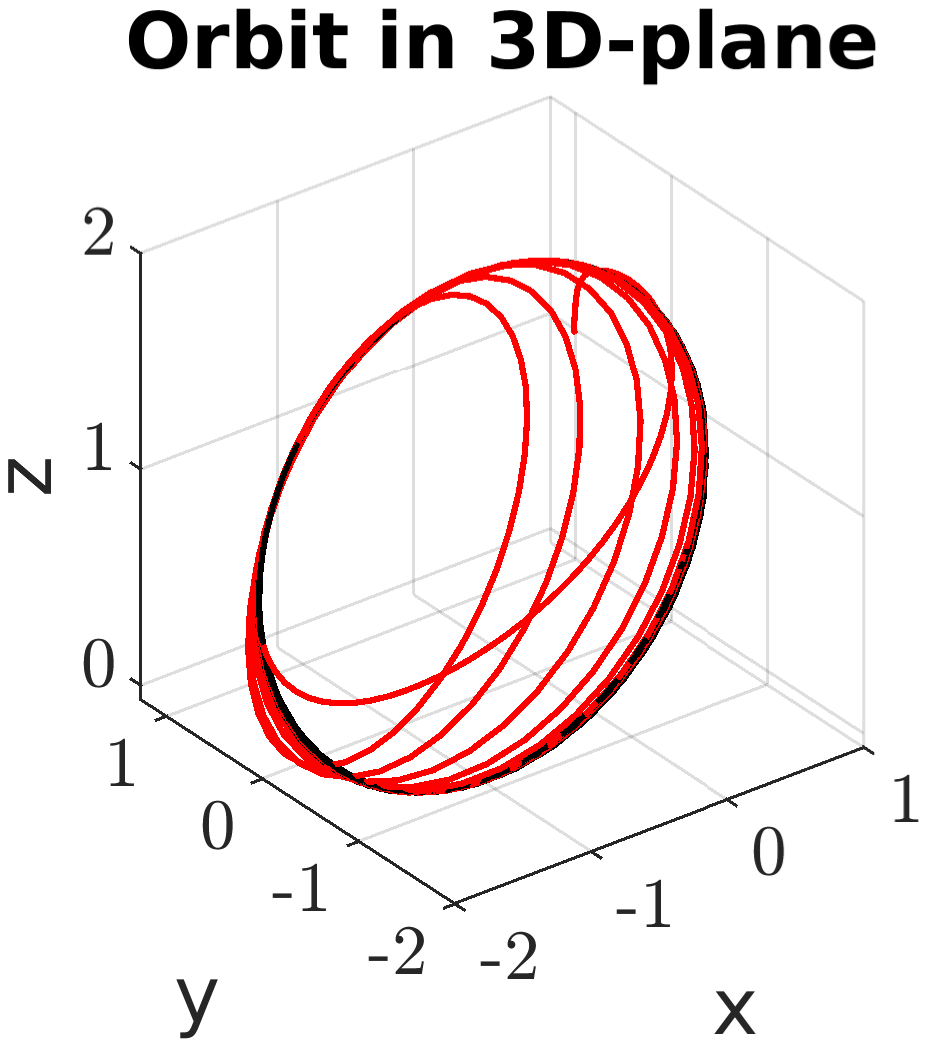}
\vspace{-1.cm}
\caption{Plots of optimal trajectory in red and trim in black. The plot on the left shows the convergence of $(y,z)$ coordinates of optimal trajectory to $(y,z)$ coordinates of the trim over evolving time, 
while the right plot shows the convergence of geometrical curves in phase space.} 
\label{fig:rigidbody.target_e1}
\end{figure}

\begin{figure}[ht]  
\centering
\def\size{0.33}
\def\Size{-0.7}
\def\Sizev{-2.}
\hspace{-0.2cm}
\includegraphics[width=\size\textwidth]{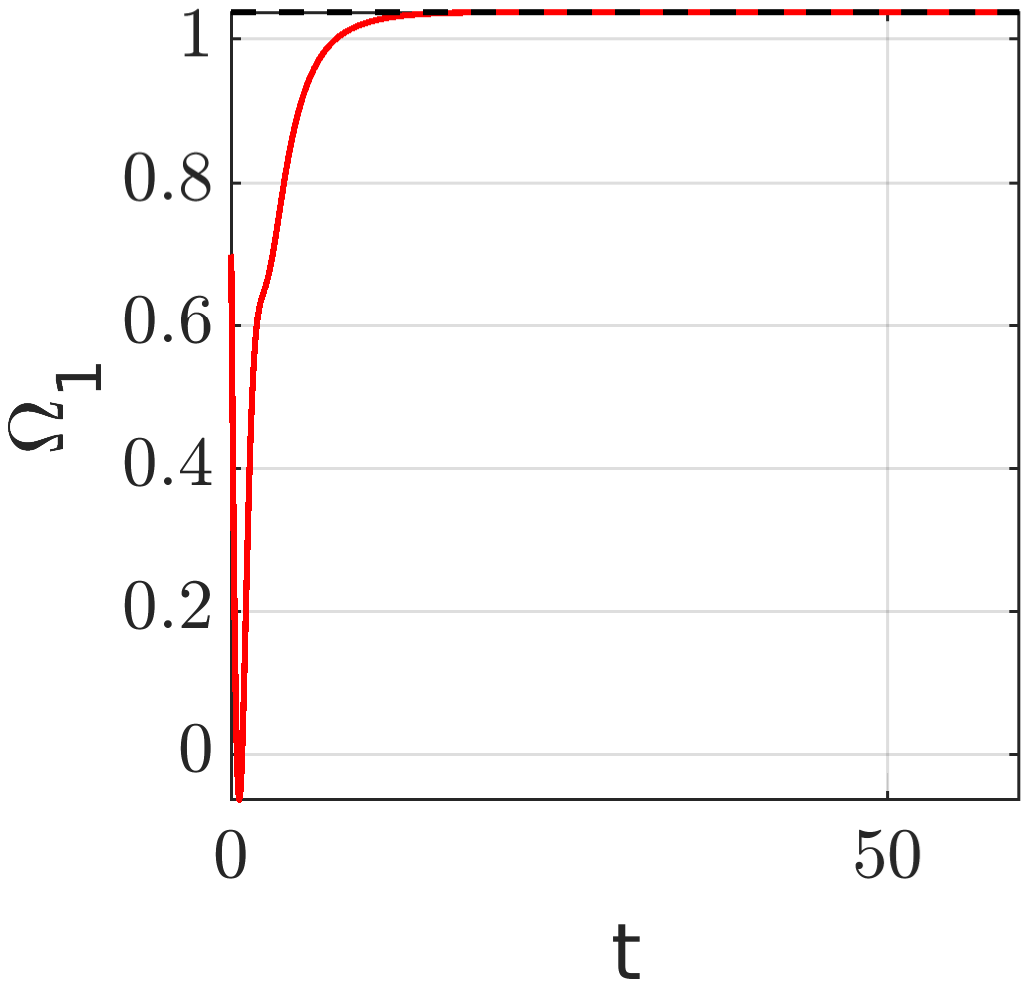}
\hspace{\Size cm}
\includegraphics[width=\size\textwidth]{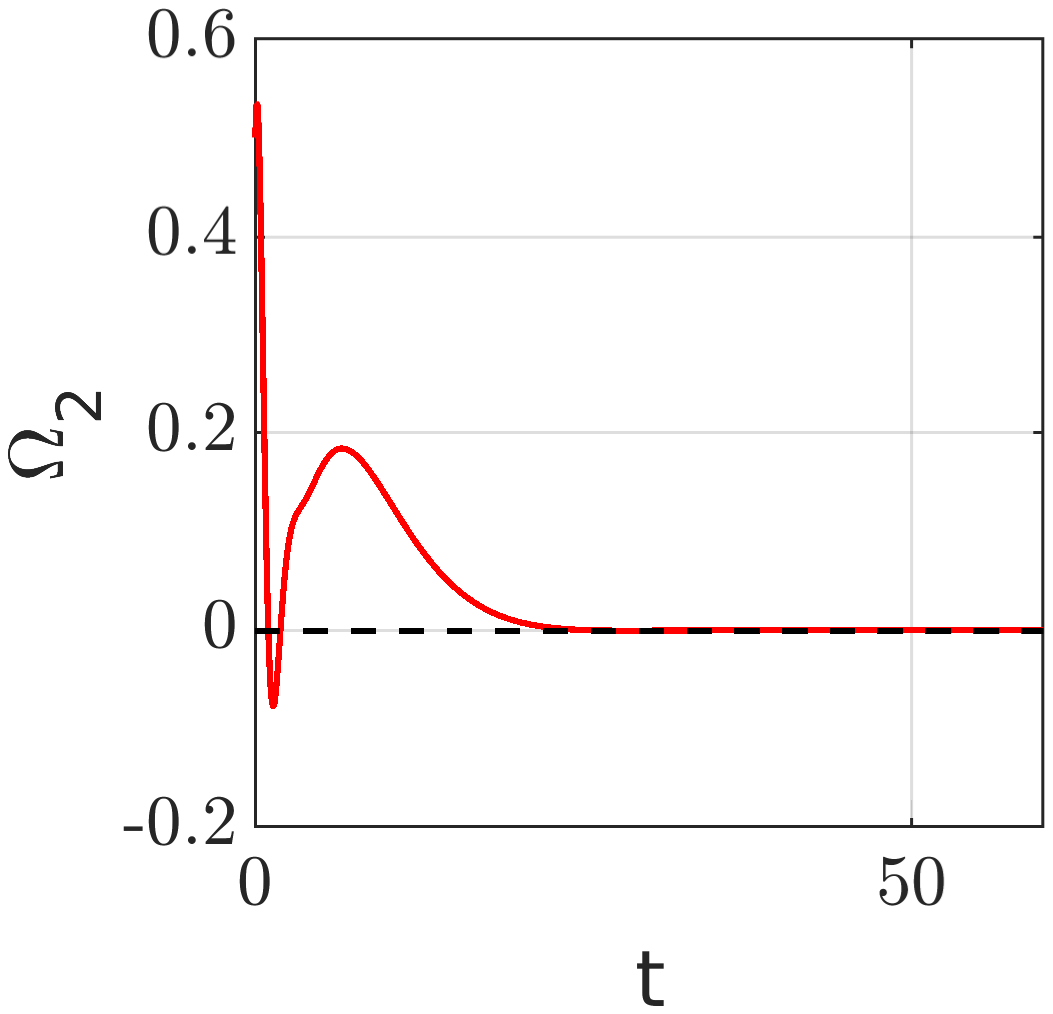}
\hspace{\Size cm}
\vspace{\Sizev cm}
\includegraphics[width=\size\textwidth]{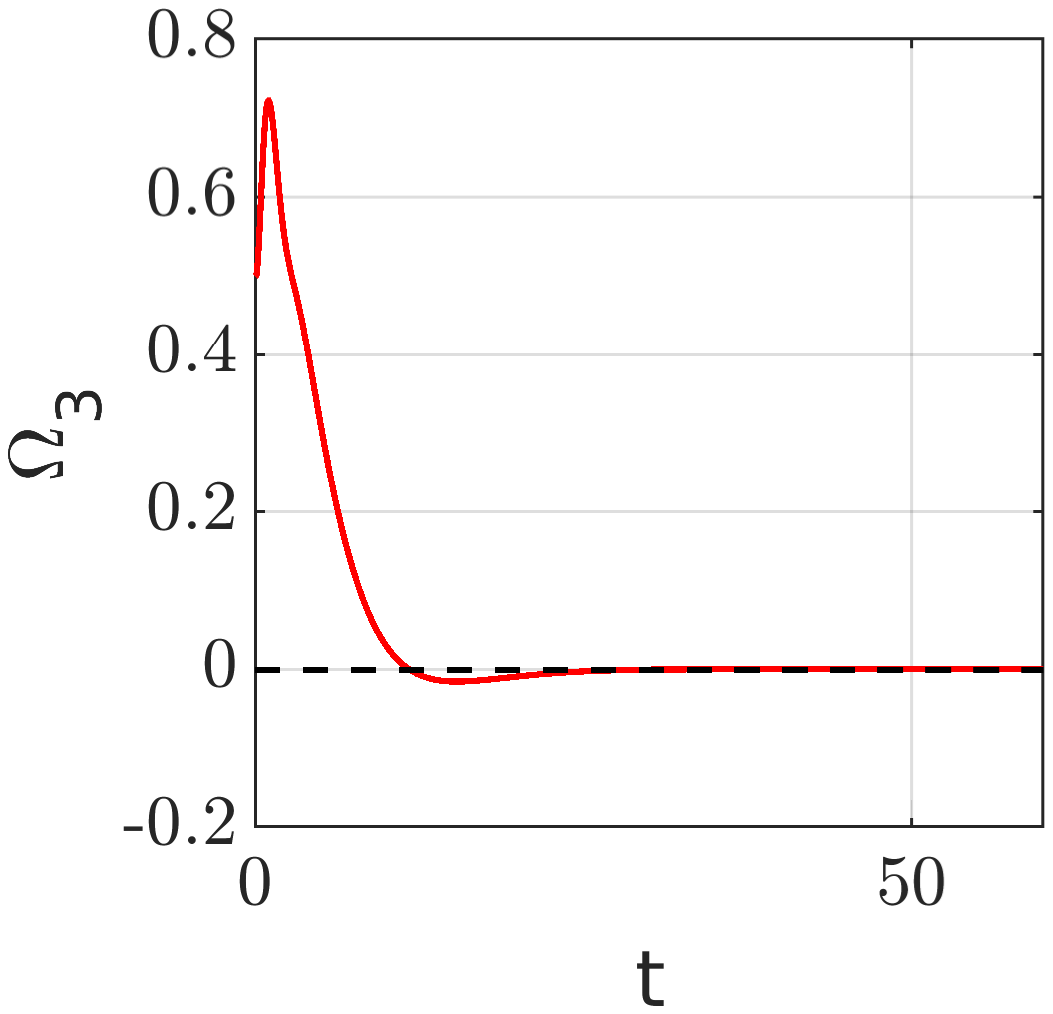}
\includegraphics[width=\size\textwidth]{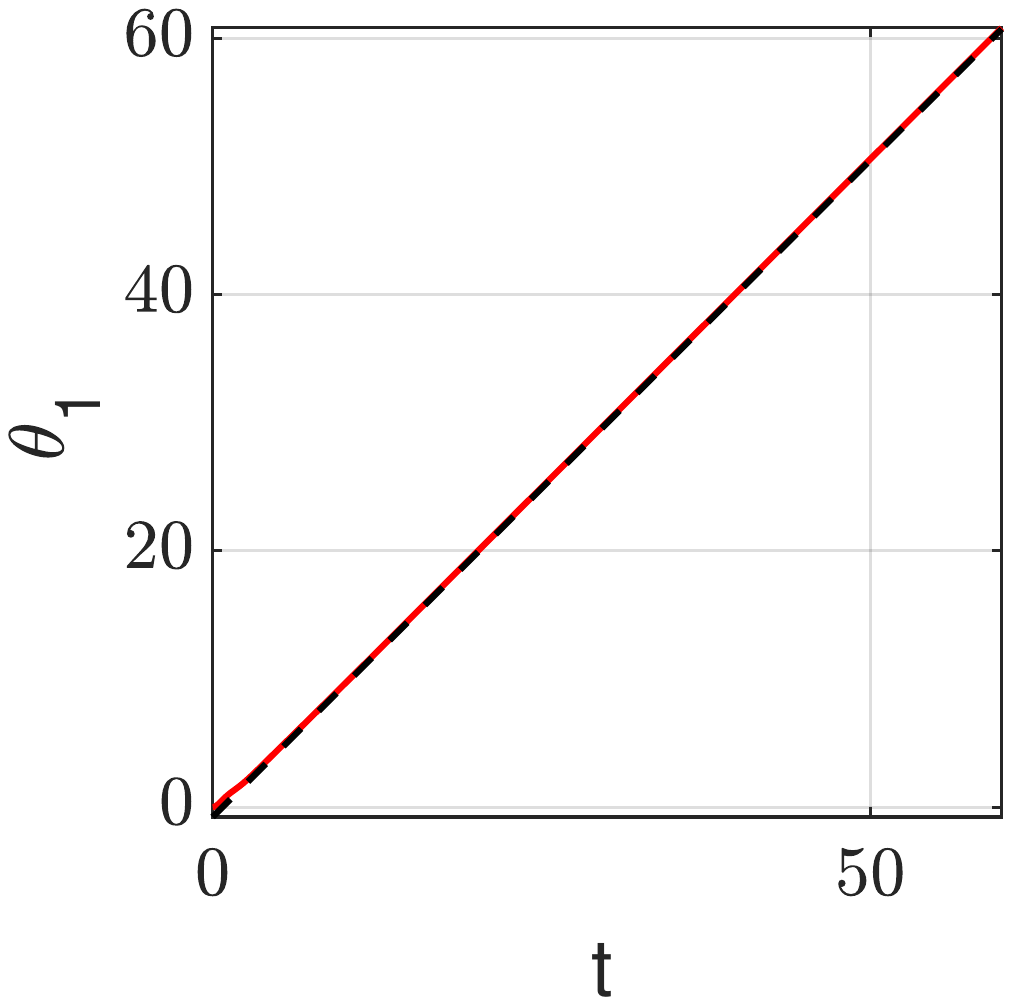}
\hspace{\Size cm}
\includegraphics[width=\size\textwidth]{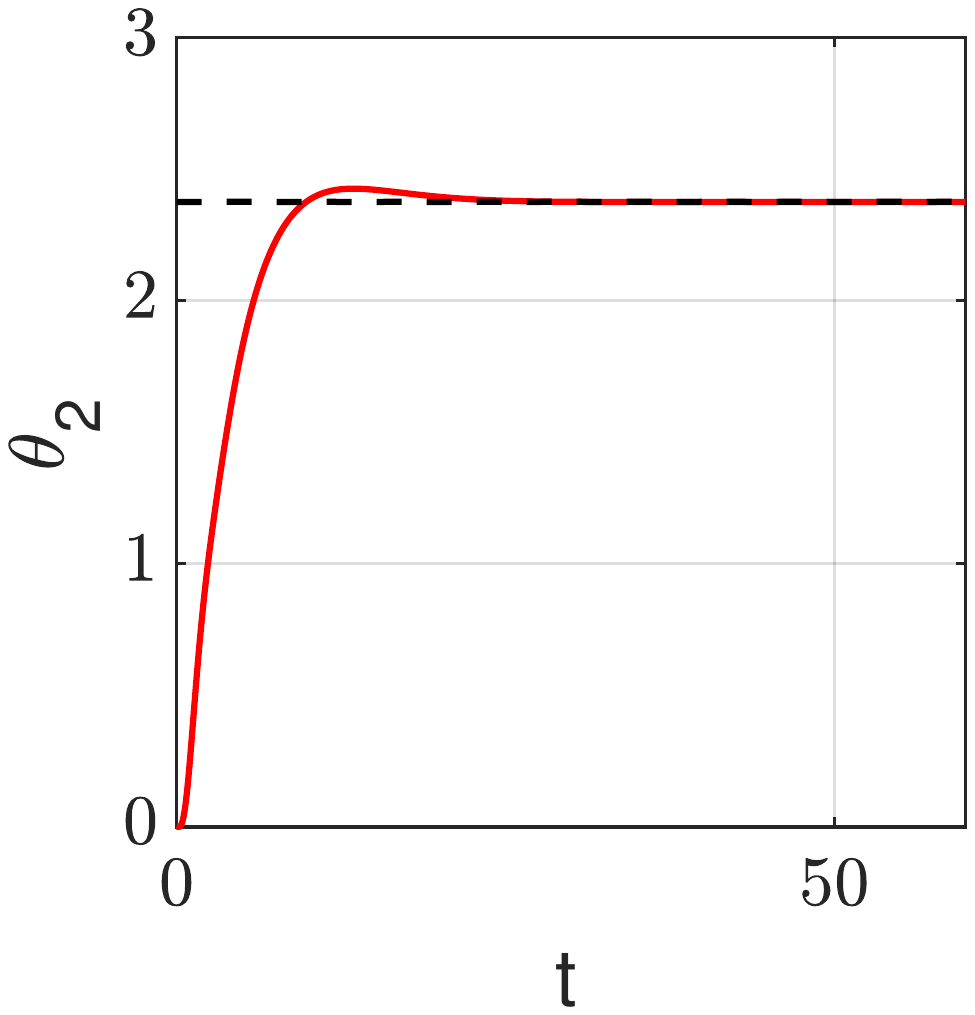}
\hspace{\Size cm}
\vspace{\Sizev cm}
\includegraphics[width=\size\textwidth]{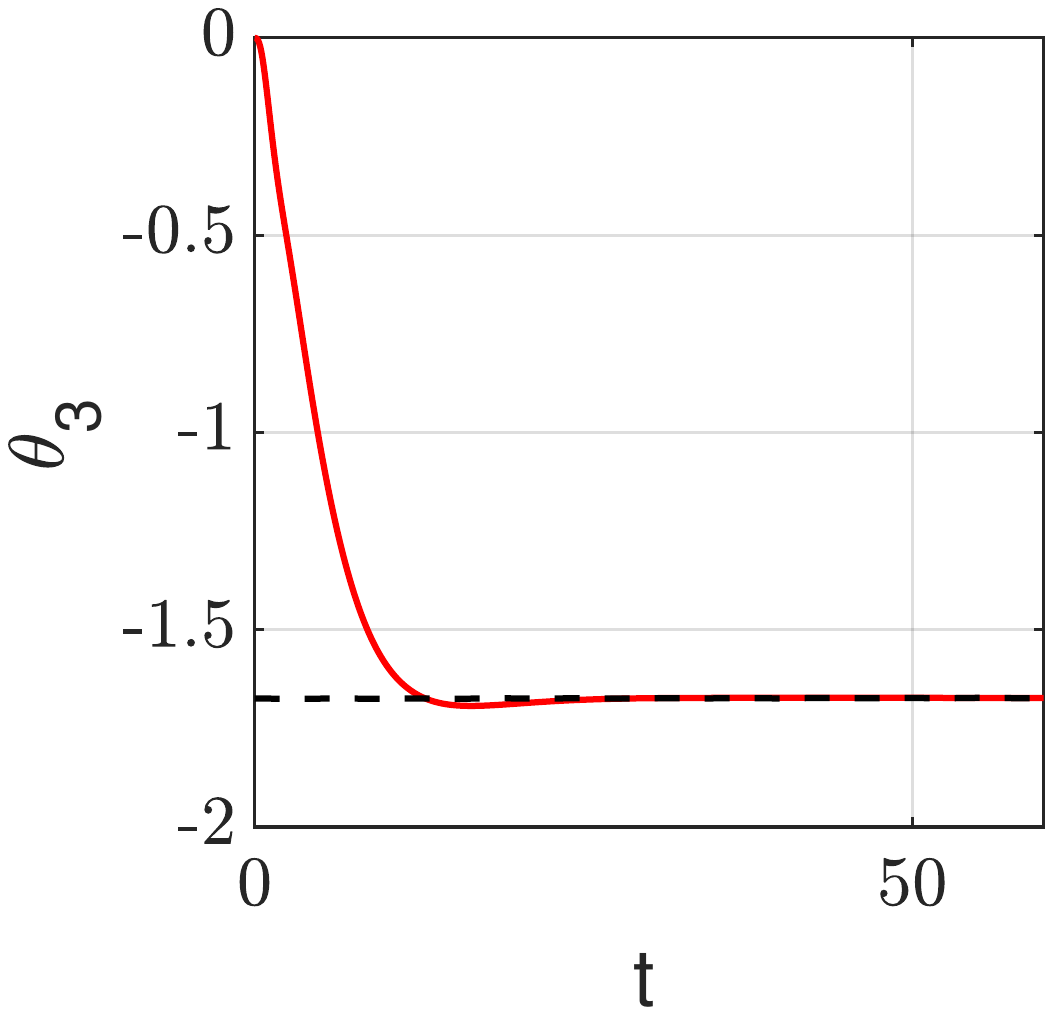}
\vspace{\Sizev cm}
\includegraphics[width=\size\textwidth]{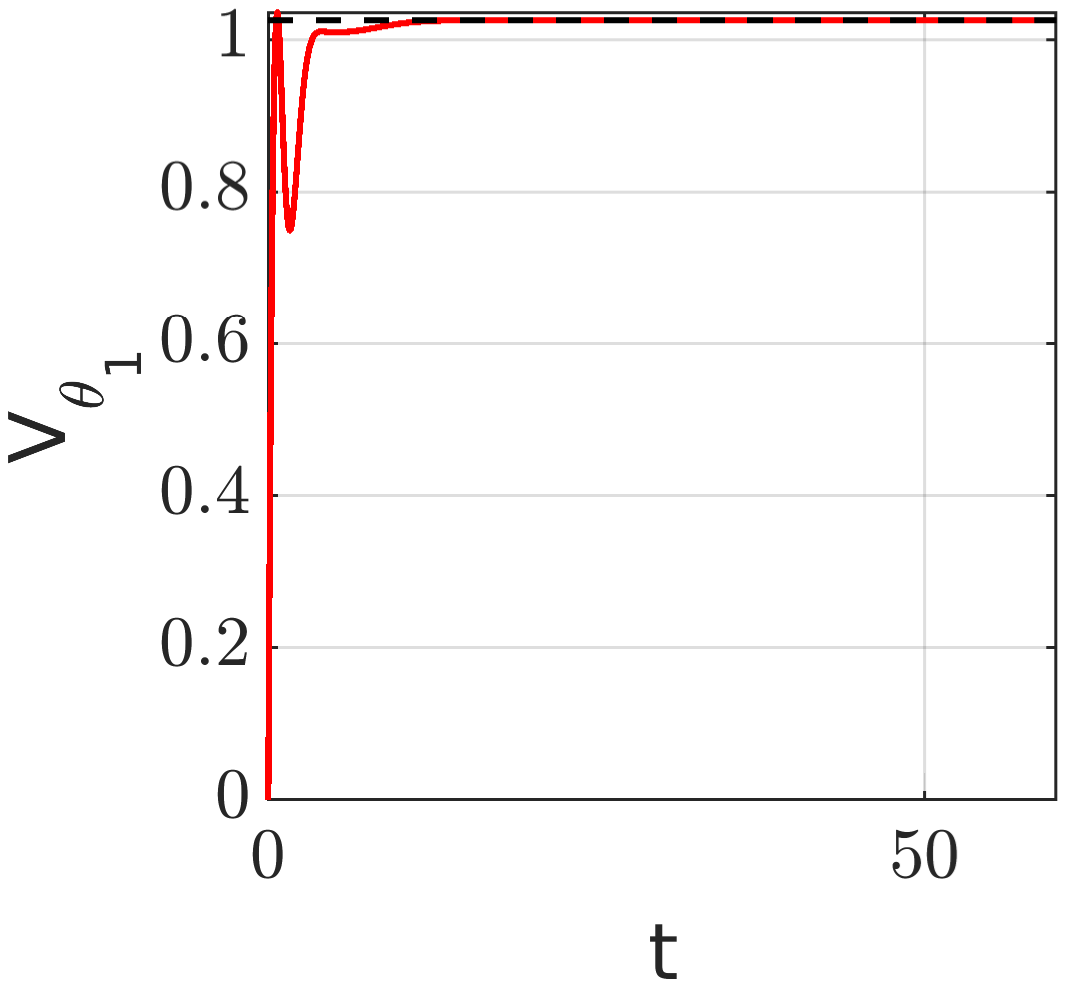}
\hspace{\Size cm}
\includegraphics[width=\size\textwidth]{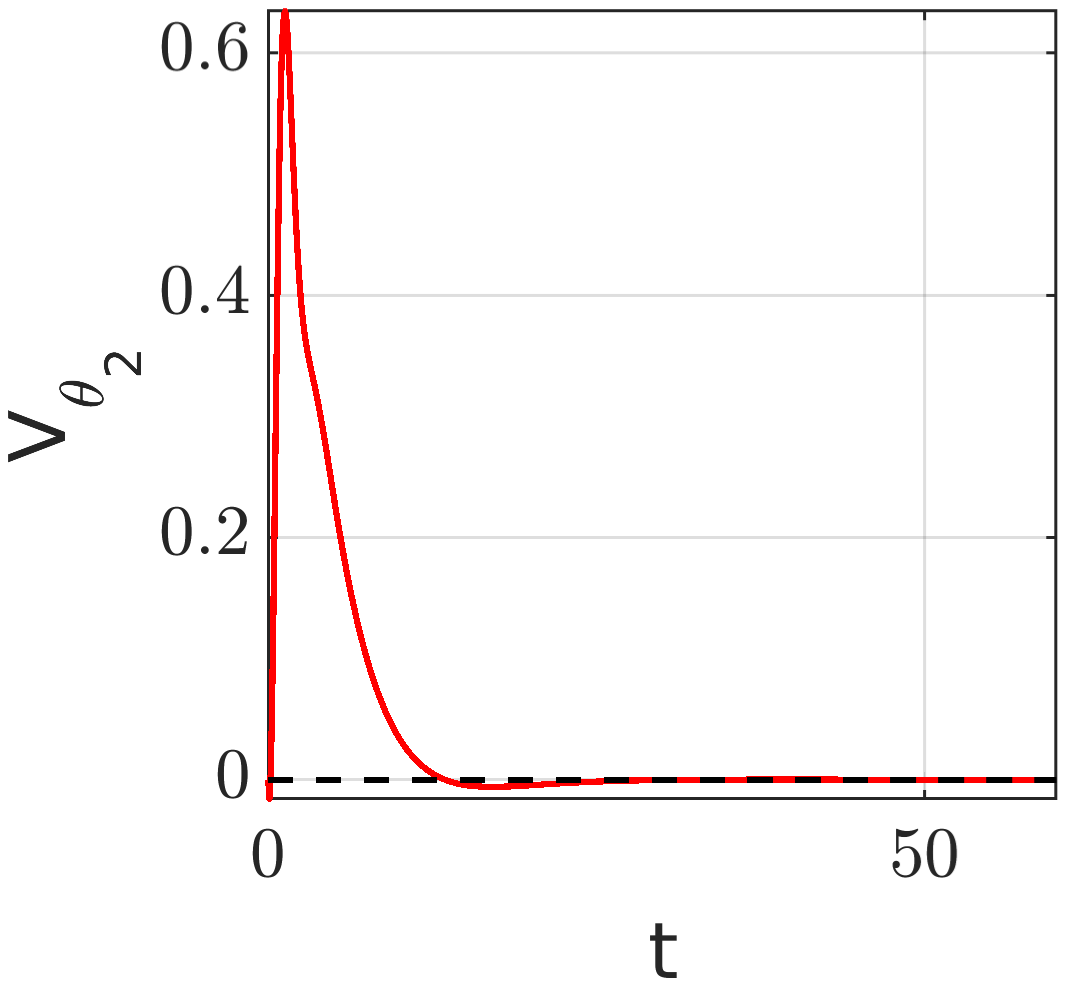}
\hspace{\Size cm}
\includegraphics[width=\size\textwidth]{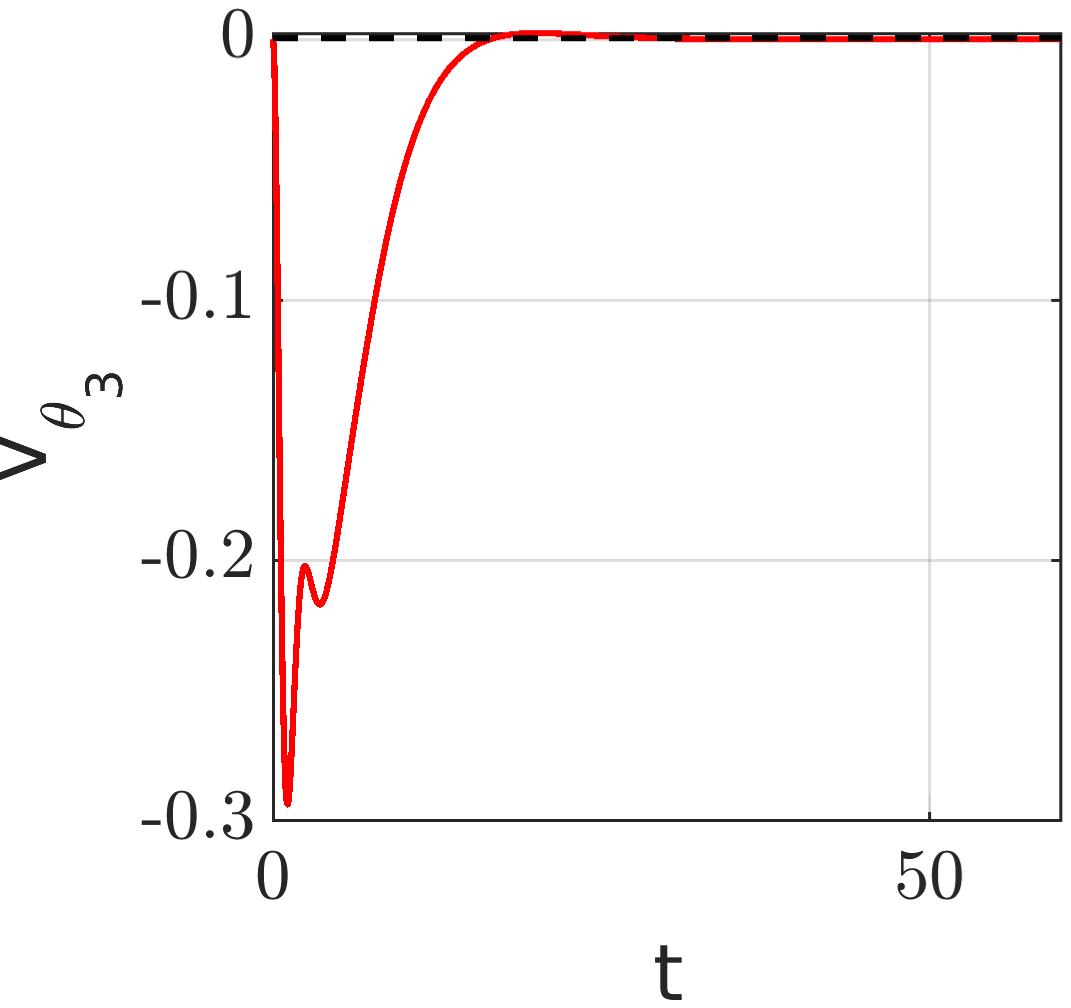}
\vspace{-1.cm}
\includegraphics[width=\size\textwidth]{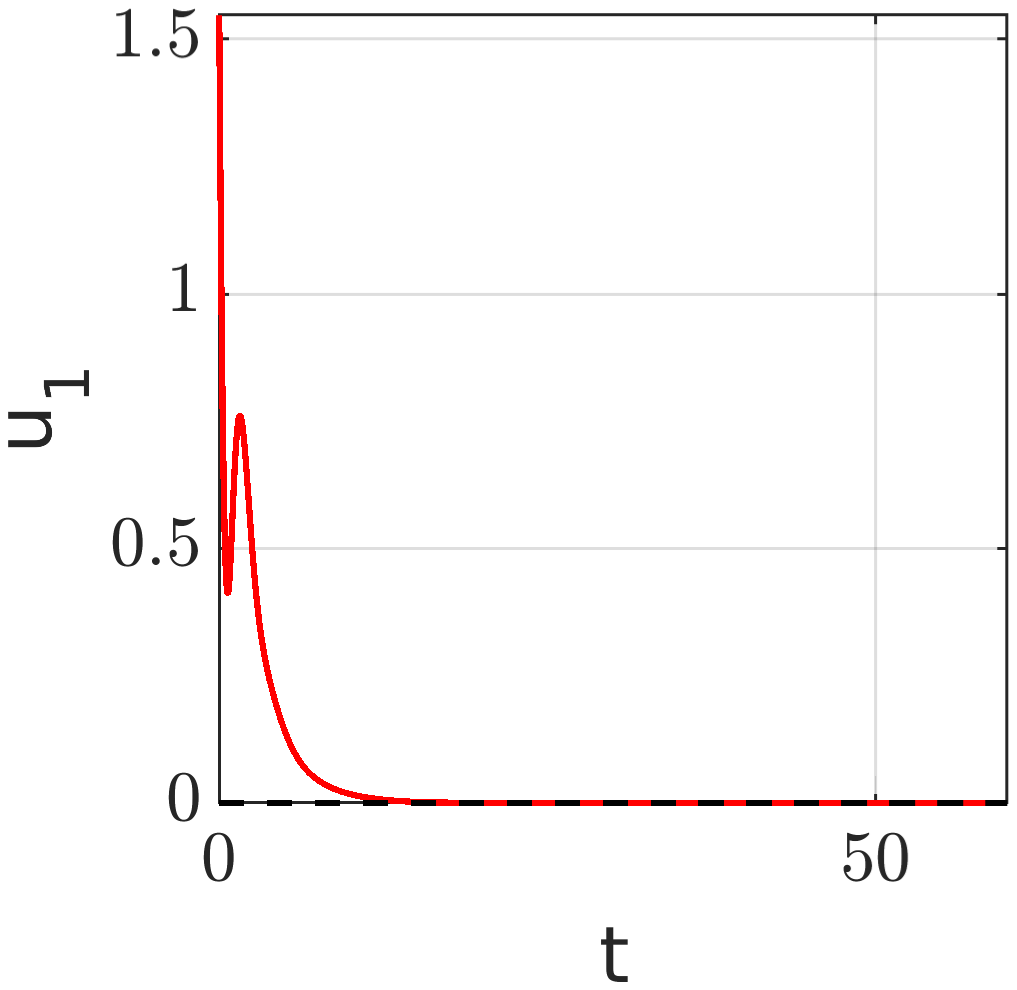}
\hspace{\Size cm}
\includegraphics[width=\size\textwidth]{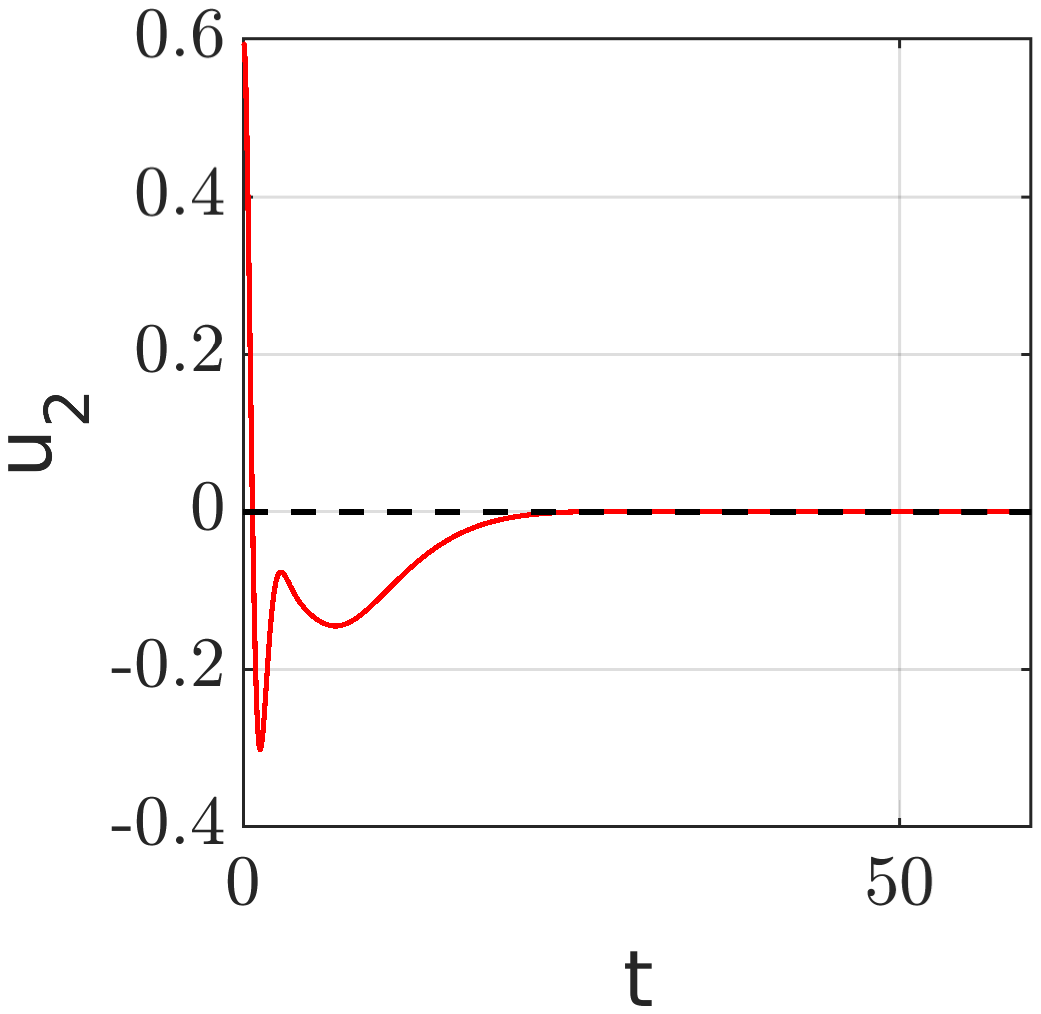}
\hspace{\Size cm}
\includegraphics[width=\size\textwidth]{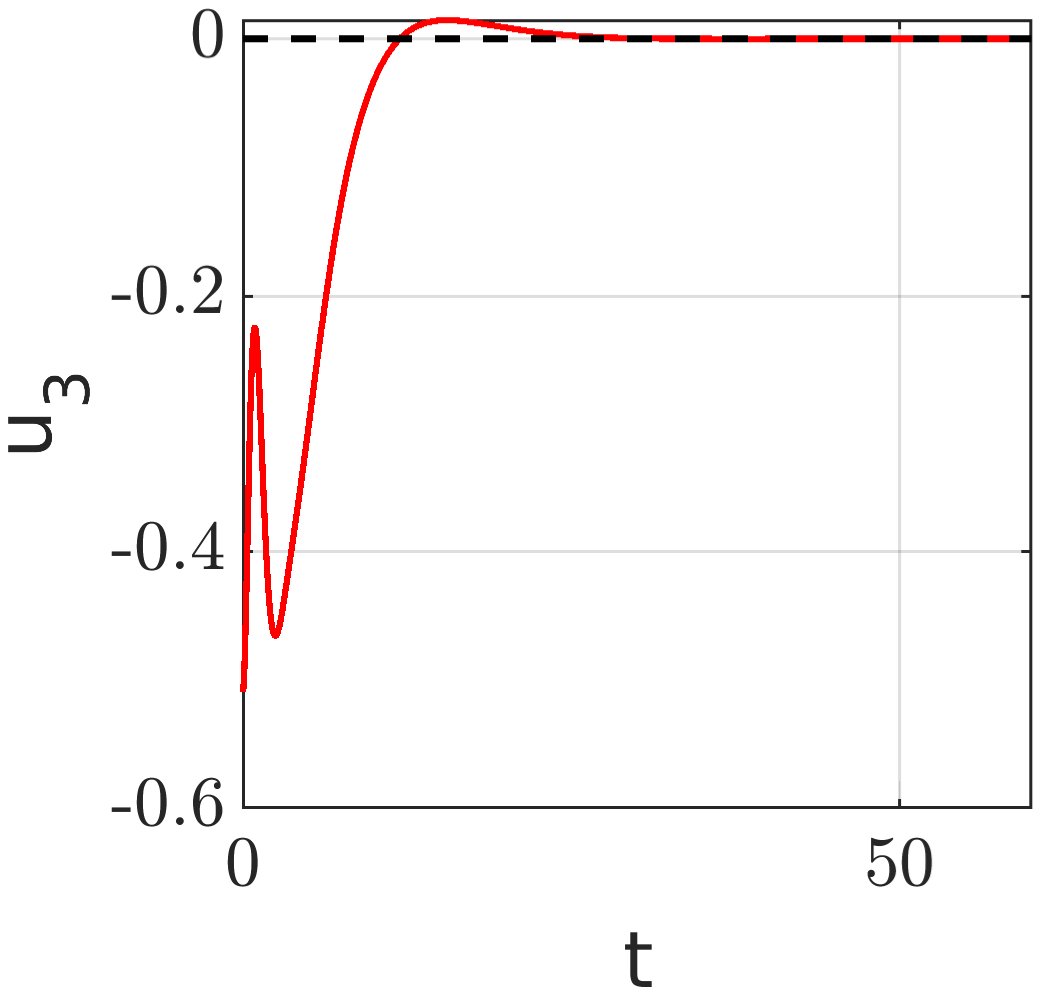}
%
\vspace{-0.2cm}
\caption{The turnpike property is depicted for the control $u$ and state variables $\Omega, \theta, v_\theta$. The trim depicted in dashed black and the optimal solution in red. Here we observe the partial turnpike as the trim is not constant in $\theta_1$.}
\label{fig:rigidbody.all.variables}
\end{figure}

To obtain numerical solutions of \eqref{eq:partial-reduced-rigidbody-pb-2}, we followed the same implementation steps as described in Section~\ref{sec:rb}. 
The obtained results for a fixed initial and free final conditions are shown on Figures~\ref{fig:rigidbody.target_e1} and \ref{fig:rigidbody.all.variables}. We observe a converge of the optimal solution to the trim for all the variables except $\lambda_{\Omega_1}$ and $\lambda_{v_{\theta_1}}$. 
This example shows that the hyperbolicity of the linearization around a static solution is not a necessary condition for trim turnpikes (at least on state variables). But, this has to be investigate more rigorously. 

\section{Conclusion} \label{sec:conclution}

In this paper, a new exponential trim turnpike property has been derived for optimal control problems with symmetries. More precisely, we showed that under certain conditions, optimal solutions converge exponentially to trim primitives, which are orbits traveling with constant velocity. These findings are based on the reduction of the problem with respect to a symmetry, which allowed us to formulate the conditions for the new turnpike to hold and to find the precise trims toward which the solutions converge. In addition, combining mechanical reduction (introduced by Marsden et al. \cite{CMR:2001b, CM:2004}) and optimal control reduction (developed by Ohsawa \cite{Ohsawa13}) we obtained a reduced problem in a simple form for controlled mechanical systems.
Two particular systems were discussed: the Kepler problem and the Rigid body  (without and with attached rotors). Numerical simulations obtained for these non-trivial academic examples validate the obtained theoretical results.
 
In the presented approach, we have considered restricted boundary conditions, which only fix the initial condition on the group. Next step would be to show the exponential trim turnpike in case of general boundary conditions, such as the conditions considered in \cite{Tre:23}, for a more restricted class of optimal control problems. 
In addition, the example of rigid body with rotors shows the trim turnpike in absence of hyperbolicity property, which highlights the importance of the sufficient conditions for trim turnpike. This is another direction of the future research on the trim turnpikes.   

\appendix

\begin{figure}[ht]  
\centering
\def\size{0.33}
\def\Size{-0.7}
\def\Sizev{-2.5}
\includegraphics[width=\size\textwidth]{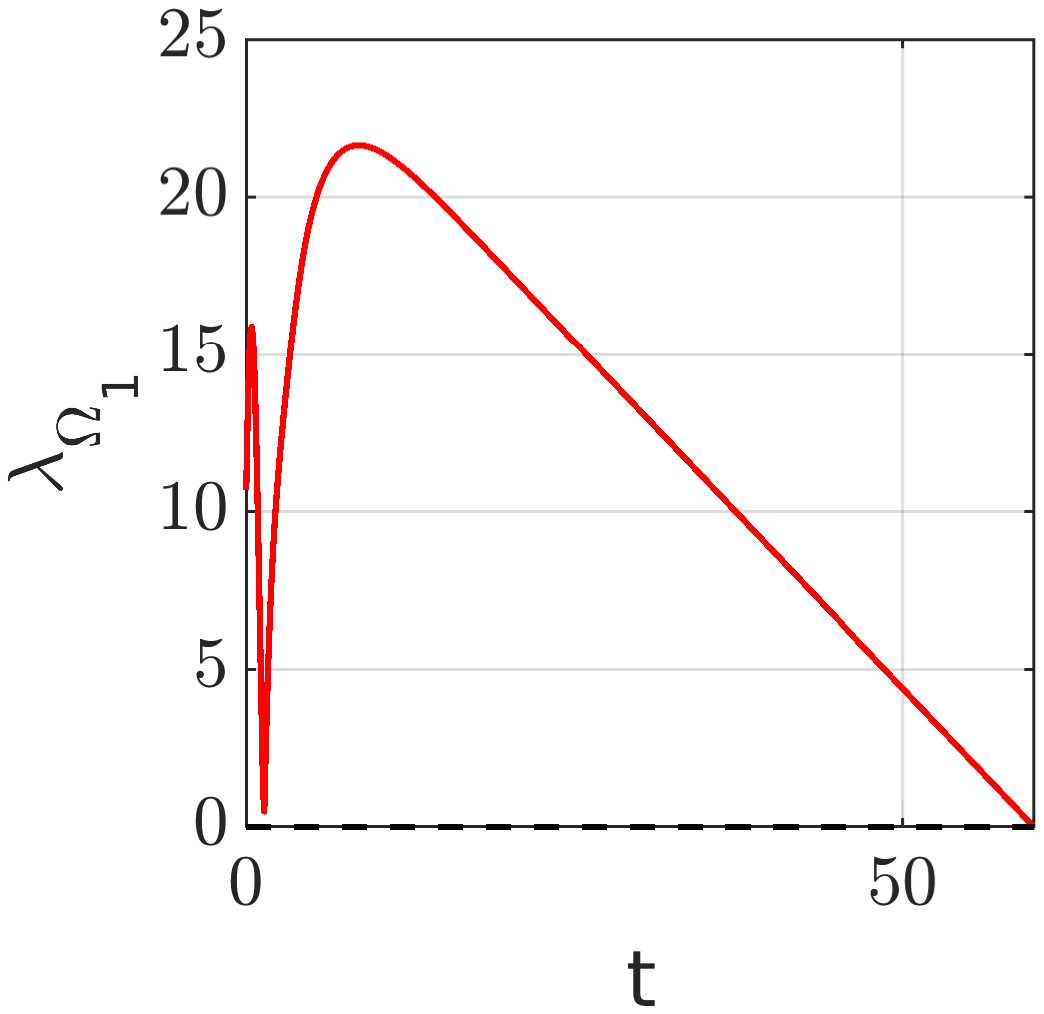}
\hspace{\Size cm}
\includegraphics[width=\size\textwidth]{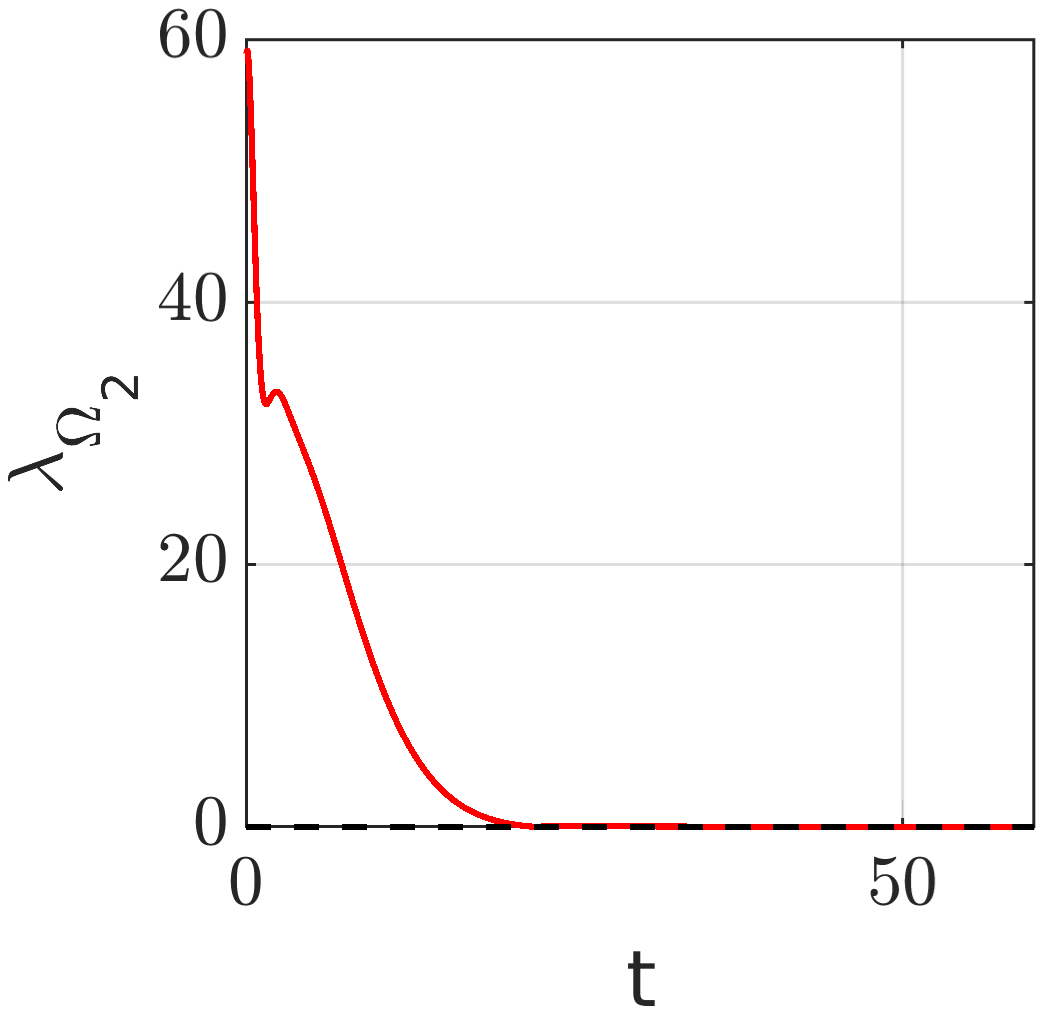}
\hspace{\Size cm}
\vspace{\Sizev cm}
\includegraphics[width=\size\textwidth]{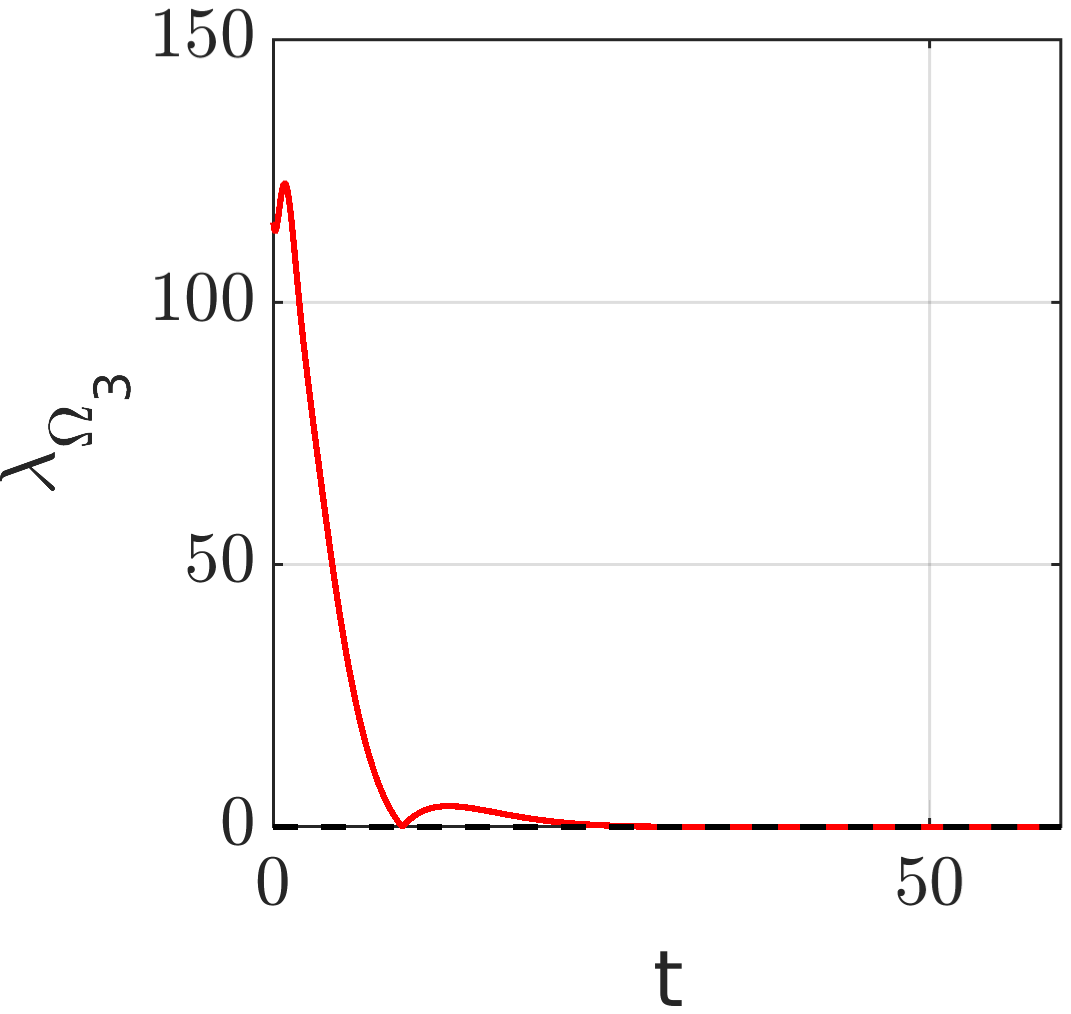}
\hspace{\Size cm}
\includegraphics[width=\size\textwidth]{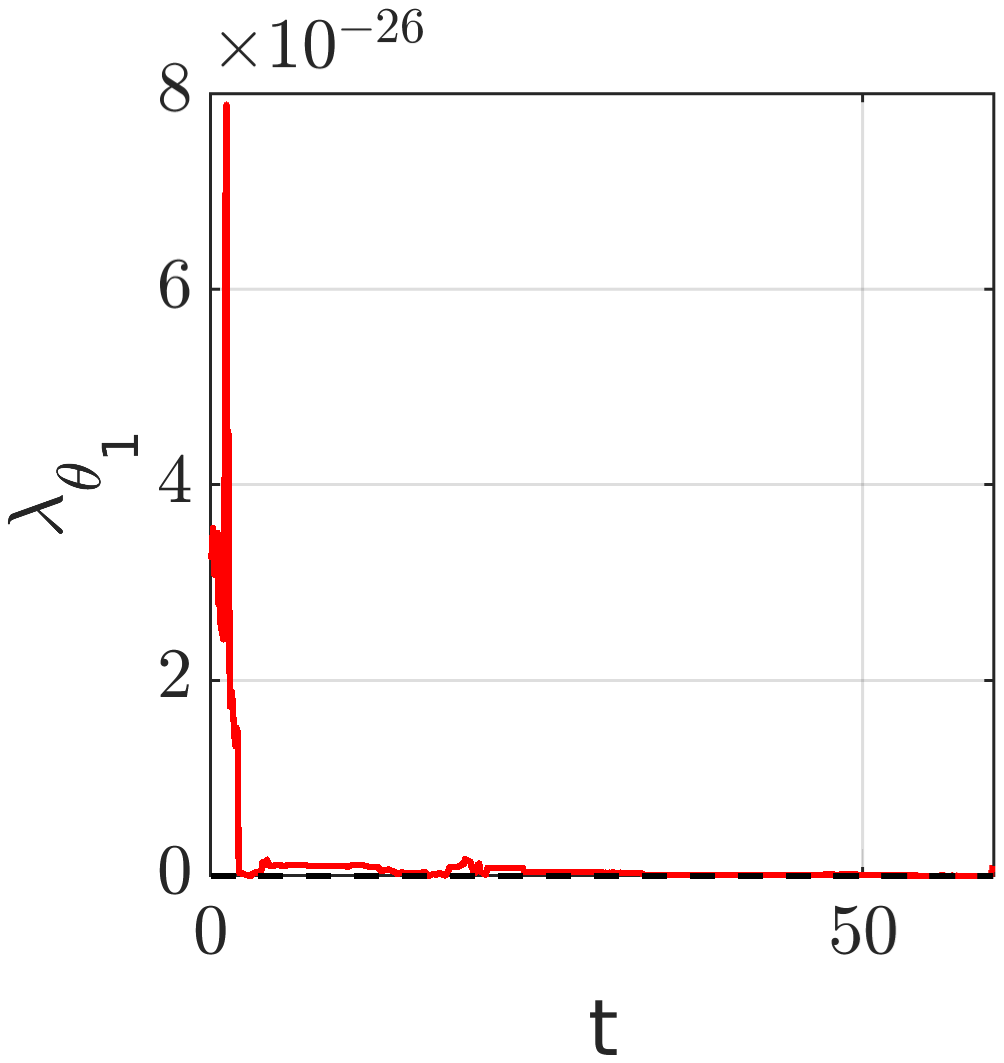}
\hspace{\Size cm}
\includegraphics[width=\size\textwidth]{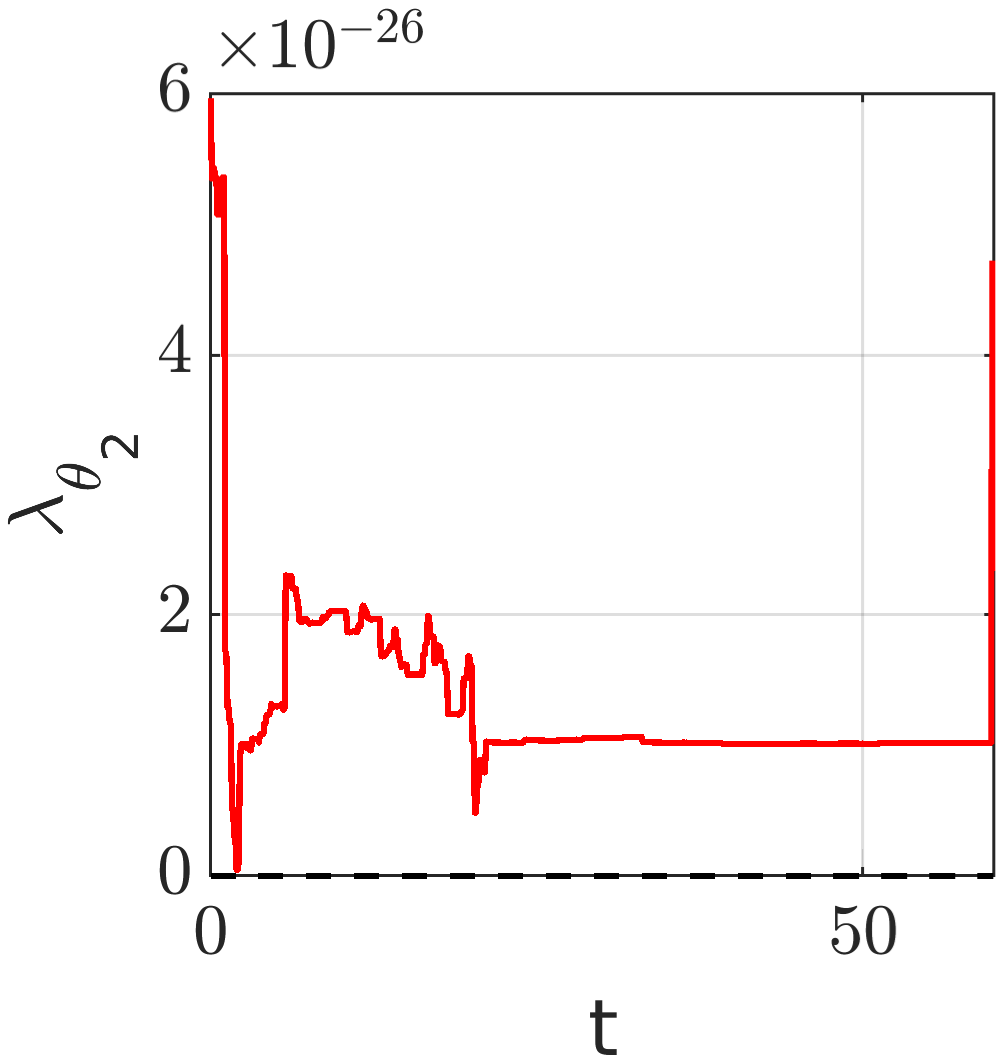}
\vspace{\Sizev cm}
\hspace{\Size cm}
\includegraphics[width=\size\textwidth]{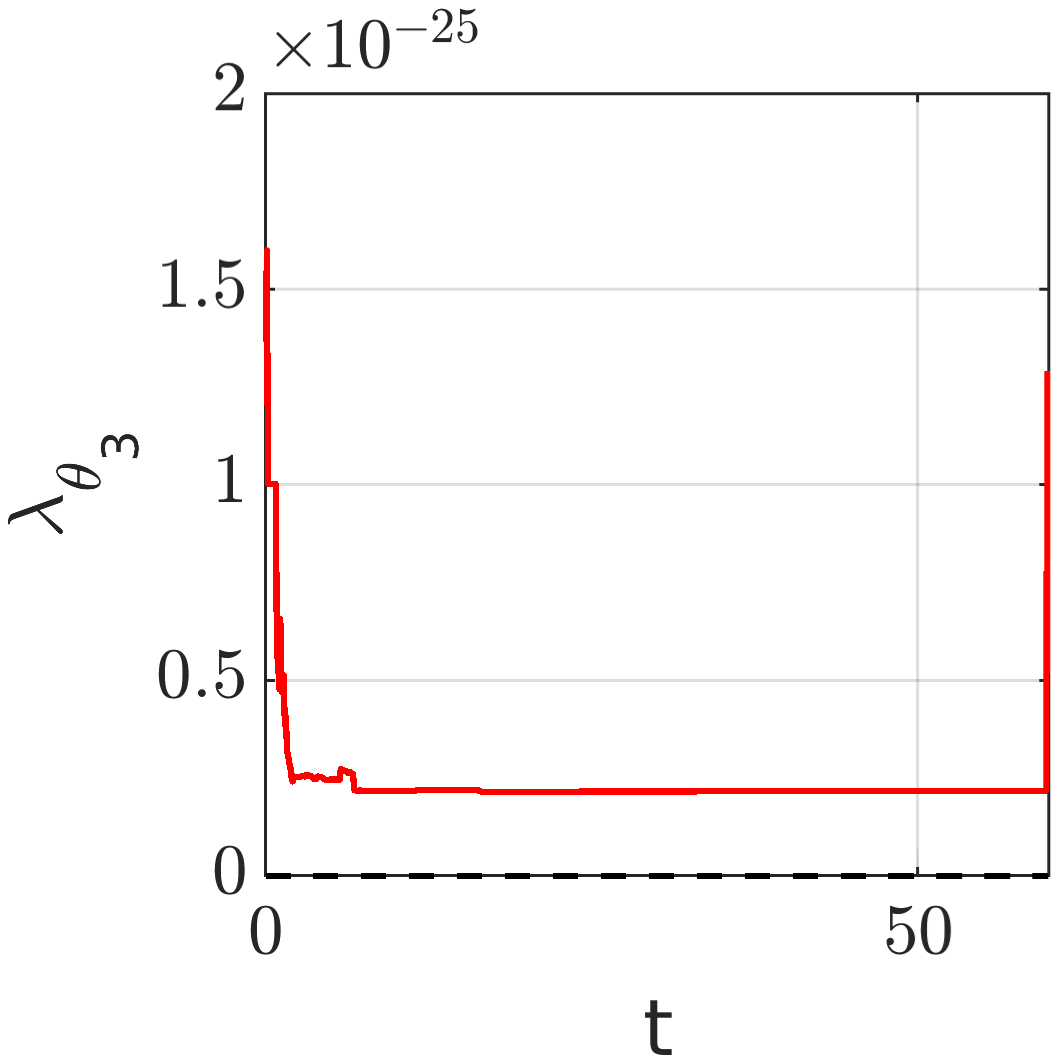}
\includegraphics[width=\size\textwidth]{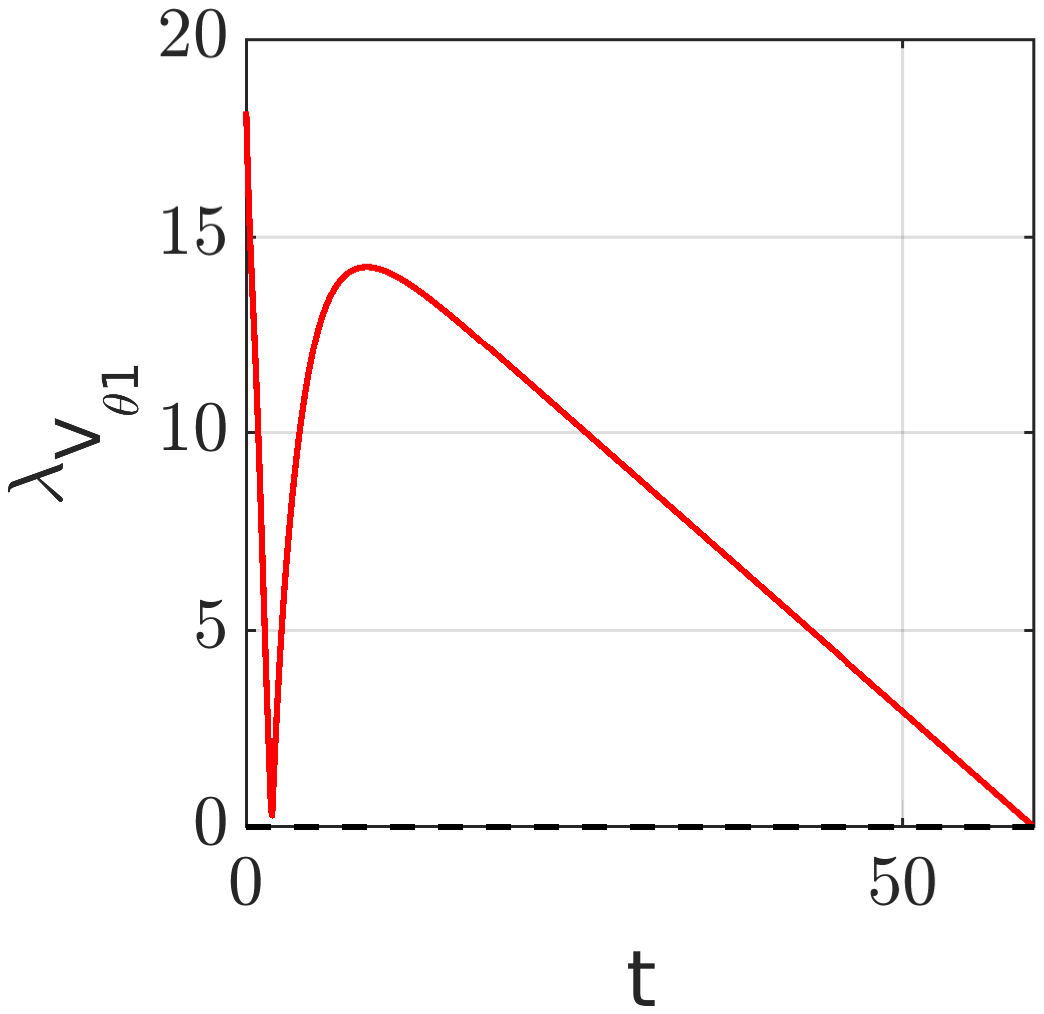}
\hspace{\Size cm}
\includegraphics[width=\size\textwidth]{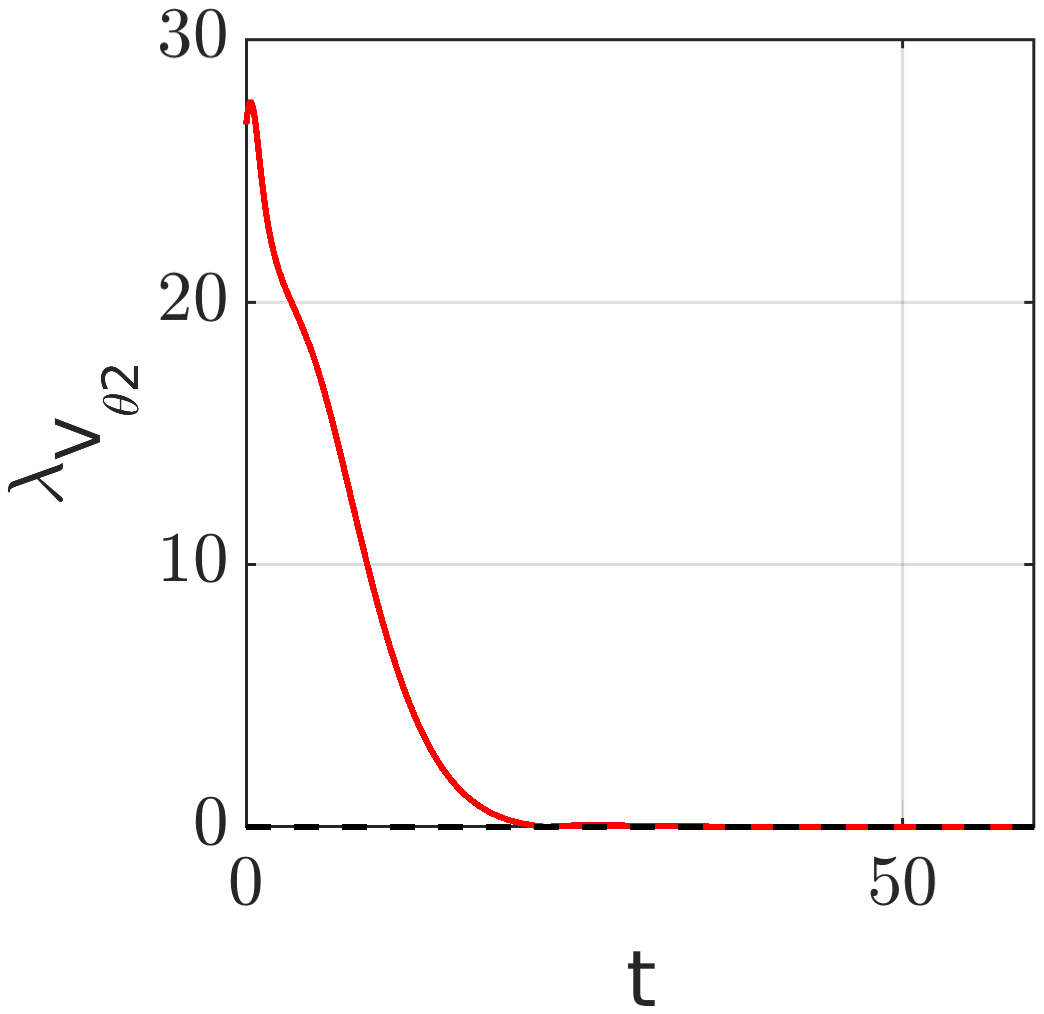}
\hspace{\Size cm}
\includegraphics[width=\size\textwidth]{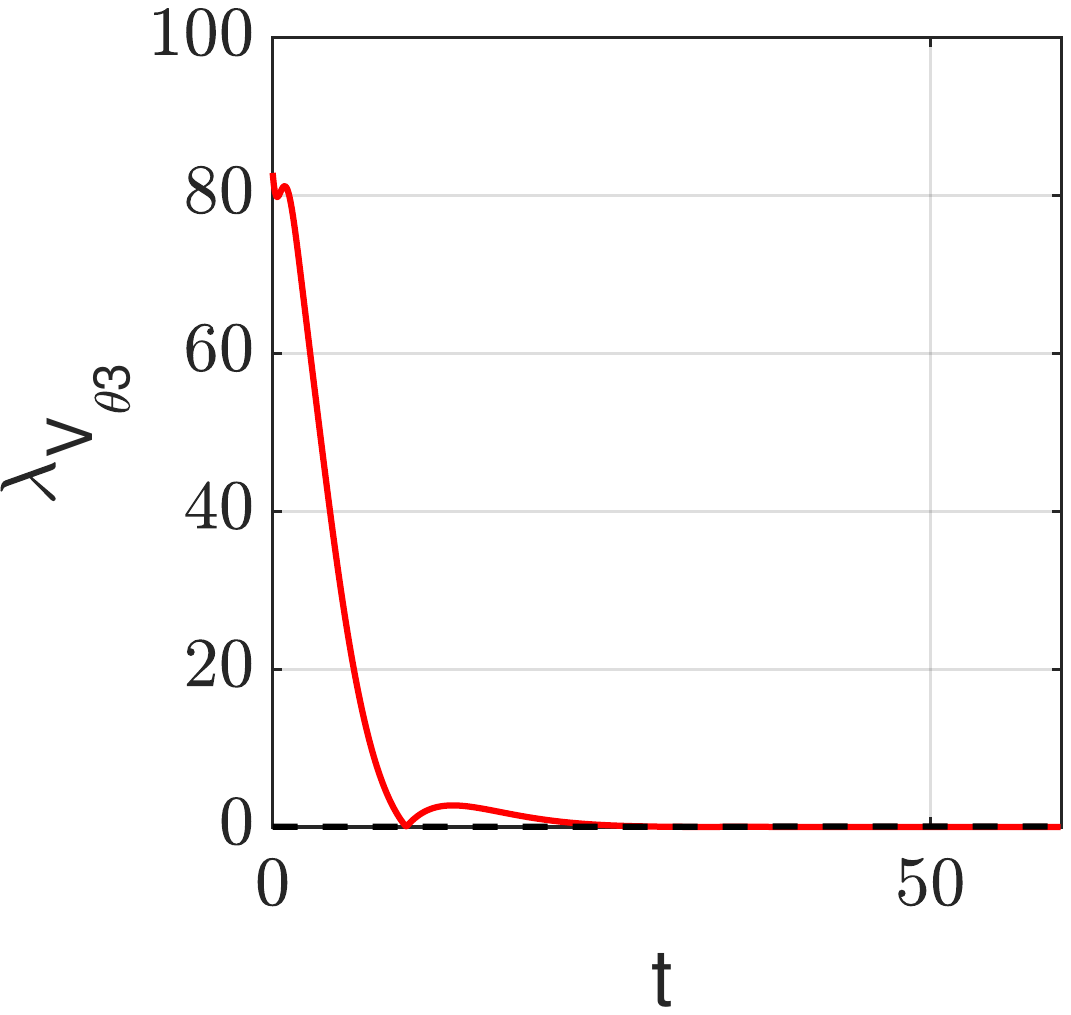}
\vspace{-1.3cm}
\caption{Representation of the adjoint variables for Rigid body problem with rotors. The non-convergence of some adjoint variables, especially $\lambda_{\Omega_1}$ and $\lambda_{V_{\theta1}}$, toward the trim (given here by $0$) is a consequence of the lack of hyperbolicity property.}
\label{fig:rigidbody.adjoint.variables}
\end{figure}

\bibliographystyle{model1-num-names}

\end{document}